\documentclass{article}

\usepackage{arxiv}
\usepackage{graphicx, color}
\usepackage{amsmath,amssymb}
\usepackage[utf8]{inputenc} 
\usepackage[T1]{fontenc}    
\usepackage{hyperref}       
\usepackage{url}            
\usepackage{booktabs}       
\usepackage{amsfonts}       
\usepackage{microtype}      
\usepackage{lipsum}

\title{Exact Solutions of the Cubic-Quintic Duffing Equation Using Leaf Functions}

\author{
  Kazunori Shinohara\thanks{10-3 Takiharu-cho, Minami-ku, Nagoya 457-8530, Daido University, Japan} \\
  Department of Mechanical Systems Engineering\\
  Daido University\\
  10-3 Takiharu-cho, Minami-ku, Nagoya 457-8530, Japan \\
  \texttt{shinohara@06.alumni.u-tokyo.ac.jp} \\
}

\begin{document}
\maketitle

\begin{abstract} 
The exact solutions of both the cubic Duffing equation and cubic-quintic Duffing equation are presented by using only leaf functions. In previous studies, exact solutions of the cubic Duffing equation have been proposed using functions that integrate leaf functions in the phase of trigonometric functions. Because they are not simple, the procedures for transforming the exact solutions are complicated and not convenient. The first derivative of the leaf function $\mathrm{sleaf}_2(l)$ can be derived as the root of $1- (\mathrm{sleaf}_2(l))^4$. This derivative can be factored, and the factors are $1+\mathrm{sleaf}_2(l)$, $1-\mathrm{sleaf}_2(l)$, and $1+ (\mathrm{sleaf}_2(l))^2$. These factors or multiplications of factors are exact solutions to the Duffing equation. Some of these exact solutions are of the same type as the cubic Duffing equation reported in previously. Some of these exact solutions satisfy the exact solutions of the cubic--quintic Duffing equations with high nonlinearity. These rules can also be applied to the leaf function $\mathrm{cleaf}_2(l)$. 
In this study, the relationship between the parameters of these exact solutions and the coefficients of the terms of the Duffing equation is clarified. We numerically analyze these exact solutions, plot the waveform, and discuss the periodicity and amplitude of the waveform. 
\end{abstract}


\section{Introduction}
The Duffing equation is a nonlinear second-order differential equation. The equation consists of the first derivative, the second derivative, the polynomial $x$, and the trigonometric function (the external force term). The behavior of the solution of the Duffing equation easily changes depending on the initial value and the polynomial coefficients, and it is difficult to predict its solution.
To clarify the behavior of the solution, research based mainly on numerical analysis with high-precision calculation is conducted. Different forms of the fractional Duffing equation have been analyzed by using the different parameters of the equation and its fractional derivative orders \cite{Pirmohabbati}. The efficient multistep differential transform method has been applied to obtain accurate approximate solutions \cite{Khatami}. An analytical approximate technique combining both the homotopy perturbation method and variational formulation has also been presented \cite{RAZZAK20162959}. Lai et al. proposed a method of linearized harmonic balance with Newton's method for solving accurate analytical approximations \cite{LAI2009852}.  To obtain an accurate analytical solution for Duffing equations with cubic and quintic nonlinearities, Pirbodaghi et al. proposed the Homotopy Analysis Method and Homotopy Pade technique \cite{PIRBODAGHI}. 
Ganji et al. proposed the complete EBM(Energy Balance Method  solution procedure of the cubic-quintic Duffing oscillator equation \cite{domiri}.
Beléndez et al. obtained an accurate approximate closed-form solutions for the cubic-quintic Duffing oscillator \cite{Bel}
Tsiganov proposed a few exact discretizations of one-dimensional cubic and quintic Duffing oscillators sharing the form of the Hamiltonian and canonical Poisson bracket up to the integer scaling factor \cite{Tsiganov}.

Das et al. proposed some modifications of a domain decomposition method to obtain accurate closed-form approximate solutions of Duffing- and Li\'{e}nard-type nonlinear ordinary differential equations \cite{Das}. 
In the literature \cite{Pirmohabbati}, it was demonstrated that exact solutions of some of these nonlinear differential equations do not exist. Therefore, investigating approximate solutions of these types of equations can play a vital role. However, exact analytical solutions or exact solutions of the Duffing equations have been discussed. To construct a transformation that converts the $n$th power to an $m$th power Duffing-type ordinary differential equation, a systematic analytical approach has been presented \cite{ZAKERI20154607}.
El\'{i}as-Z\'{u}\~{n}iga has attempted to find only one set of exact solutions by using the Jacobi elliptic function \cite{ELIASZUNIGA20132574}. Razzak obtained the exact solution of the cubic--quintic Duffing equation. However, Razzak's solution procedure is complicated and contains a set of algebraic equations with Jacobian elliptic functions that are not easily solved \cite{RAZZAK20162959}. The nonlinear differential equation governing the periodic motion of the one-dimensional, undamped, and unforced cubic-quintic Duffing oscillator is solved exactly\cite{Salas}. 

In this study, by using leaf functions, an exact solution of the cubic--quintic Duffing equation is proposed. The cubic Duffing equation and cubic--quintic Duffing equation \cite{Kovacic} \cite{Kovacic2} are represented by the following equations, respectively: 
\begin{equation}
\frac{\mathrm{d}^2x(t)}{\mathrm{d}t^2} + \alpha_1 x(t)+ \alpha_2 x(t)^3=0 ,\label{1.1}
\end{equation}
\begin{equation}
\frac{\mathrm{d}^2x(t)}{\mathrm{d}t^2} + \alpha_1 x(t)+ \alpha_2 x(t)^3+ \alpha_3 x(t)^5=0. \label{1.2}
\end{equation}

By using leaf functions, exact solutions of the cubic Duffing equation were proposed for a free oscillation system, a divergence system, and a damping system \cite{cmes.2018.02179} \cite{cmes.2019.04472}. Because these exact solutions are expressed by using trigonometric functions for the phases of leaf functions, the expansions of exact solutions become complicated and the usability of applications becomes low.
In this study, by using only leaf functions without trigonometric functions, the exact solutions of the cubic Duffing equation are presented for a free-oscillation system. For example, the following formula is obtained by differentiating the leaf function sleaf$_2(x)$: 
\begin{equation}
\sqrt{1-(\mathrm{sleaf}_2(t))^4}
=\sqrt{ \{1+\mathrm{sleaf}_2(t) \} \{1-\mathrm{sleaf}_2(t) \} \{1+(\mathrm{sleaf}_2(t))^2 \} }.
\label{1.3}
\end{equation}

As shown in the above equation, a factor obtained by factoring the equation is assumed to be the solution $x(t)$. For example, these solutions are as follows:
\begin{equation}
x(t)=\sqrt{ 1+\mathrm{sleaf}_2(t)  },\; \sqrt{ 1-\mathrm{sleaf}_2(t)  }, \;\sqrt{ 1+(\mathrm{sleaf}_2(t))^2  }, \;\sqrt{ 1-(\mathrm{sleaf}_2(t))^2  }. 
\label{1.4}
\end{equation}

A factor obtained by factoring the derivative of the leaf function or a part of the multiplication of these factors satisfies the cubic Duffing equation or the cubic--quintic Duffing equation.
Although the exact solution of the cubic Duffing equation has been discussed, this study is the first to describe the exact solution of the cubic--quintic Duffing equation.
 In this study,  the exact solutions of the cubic--quintic Duffing equation are presented. Furthermore, this equation can be extended to the Duffing equation, which includes the first derivative. The exact solutions that satisfy this extended Duffing equation can also be derived by the leaf function and the elementary function (exponential function). The exact solution can be expressed simply. The waves obtained by numerical analysis with the extended Duffing equation are in perfect agreement with the waves obtained by these exact solutions. There is consistency both numerically and theoretically.
 
In this paper, we first discuss the Duffing equation for the undamped systems. The exact solutions are graphed, and the waveforms of the exact solutions are visualized.
After that, we discuss the damped system through the comparison of the waveforms obtained by the Duffing equation for the undamped systems.

\section{Literature comparison}
Systematic study of trigonometric functions began in Hellenistic mathematics \cite{katz}. Leaf functions, hyperbolic leaf functions, inverse leaf functions and inverse hyperbolic leaf functions based on the basis $n = 1$ represent trigonometric functions, hyperbolic functions, inverse trigonometric functions ( $0 \leqq t < \frac{1}{2} \pi $) and inverse hyperbolic functions, respectively \cite{Kaz_sl} \cite{Kaz_cl} \cite{Kaz_slh} \cite{Kaz_clh} \cite{Abramowitz}.

\begin{equation}
x(t) =\mathrm{sleaf}_{1}(t)=\mathrm{sin}(t) 
\qquad \left(\quad \mathrm{arcsin}(x) = \int_{0}^{x}  \frac{\mathrm{d}u }{ \sqrt{1- u^2} }
=\mathrm{arcsleaf}_1(x)=t  \quad \right)
\label{2.1}
\end{equation}

\begin{equation}
x(t) =\mathrm{cleaf}_{1}(t)=\mathrm{cos}(t)
\qquad \left( \quad \mathrm{arccos}(x) = \int_{x}^{1}  \frac{\mathrm{d}u }{ \sqrt{1- u^2} }
=\mathrm{arccleaf}_1(x) =t  \quad \right)
\label{2.2}
\end{equation}

\begin{equation}
x(t) =\mathrm{sleafh}_{1}(t)=\mathrm{sinh}(t) 
\qquad \left( \quad \mathrm{asinh}(x) = \int_{0}^{x}  \frac{\mathrm{d}u }{ \sqrt{1+ u^2} }
=\mathrm{asleafh}_1(x)=t  \quad \right)
\label{2.3}
\end{equation}

\begin{equation}
x(t) =\mathrm{cleafh}_{1}(t)=\mathrm{cosh}(t) 
\qquad \left( \quad \mathrm{acosh}(x) = \int_{x}^{1}  \frac{\mathrm{d}u }{ \sqrt{u^2-1} }
=\mathrm{acleafh}_1(x) =t \quad \right)
\label{2.4}
\end{equation}

In 1796, Carl Friedrich Gauss presented the lemniscate function\cite{Gauss}. The inverse leaf functions based on the basis $n = 2$ represents inverse functions of the sin and cos lemniscates \cite{roy2017elliptic} (  $0 \leqq t < \frac{1}{2} \pi_2 $) and hyperbolic leaf functions, respectively \cite{Ramanujan}. (See Appendix Q for the constant $\pi_2$)

\begin{equation}
x(t) =\mathrm{sleaf}_{2}(t)=\mathrm{sl}(t) 
\qquad \left( \quad \mathrm{arcsl}(x) = \int_{0}^{x}  \frac{\mathrm{d}u }{ \sqrt{1- u^4} }
=\mathrm{arcsleaf}_2(x)={t}  \quad \right)
\label{2.5}
\end{equation}

\begin{equation}
x(t) =\mathrm{cleaf}_{2}(t)=\mathrm{cl}(t)
\qquad \left( \quad \mathrm{arccl}(x) = \int_{x}^{1} \frac{\mathrm{d}u }{ \sqrt{1- u^4} }
=\mathrm{arccleaf}_2(x) =t  \quad \right)
\label{2.6}
\end{equation}

\begin{equation}
x(t) =\mathrm{sleafh}_{2}(t)=\mathrm{slh}(t)
\qquad \left( \quad \mathrm{arcslh}(x) = \int_{0}^{x}  \frac{\mathrm{d}u }{ \sqrt{1+ u^4} }
=\mathrm{arcsleafh}_2(x)={t} \quad \right)
\label{2.7}
\end{equation}

\begin{equation}
x(t) =\mathrm{cleafh}_{2}(t)=\mathrm{clh}(t)
\qquad \left( \quad \mathrm{arcclh}(x) = \int_{x}^{1} \frac{\mathrm{d}u }{ \sqrt{ u^4-1} }
=\mathrm{arccleafh}_2(x) =t \quad \right)
\label{2.8}
\end{equation}

Historically, there had been no discussion for $n = 3$ or more. No clear description for functions exists. Therefore, the functions were defined as the leaf function and the hyperbolic leaf function \cite{booth} \cite{Dixon} \cite{Akhiezer} \cite{walker} \cite{McKean} \cite{Lawden}.  

\section{Definition of Leaf function}
\subsection{Leaf function; sleaf$_n(t)$}
An ordinary differential equation (ODE) comprises a function raised to the $2n-1$ power and the second derivative of this function. Further, the initial conditions of the ODE are defined.

\begin{equation}
\frac{\mathrm{d}^2x(t) }{\mathrm{d}t^2}=-nx(t)^{2n-1} \label{3.1.1}
\end{equation}
\begin{equation}
x(0)=0 \label{3.1.2}
\end{equation}
\begin{equation}
\frac{\mathrm{d}x(0)}{\mathrm{d}t}=1
\label{3.1.3}
\end{equation}

The variable $n$ represents an integer. In the paper, the variable is named as the basis. The preceding equation is the ODE, motivated in this study. Although the equation (\ref{3.1.1}) is a simple ordinary differential equation, it has a very important meaning because it generates characteristic waves. By numerically analyzing the solution that satisfies this equation, we can obtain regular and periodic waves \cite{Kaz_sl}, \cite{Kaz_cl}. The form of these waves differs from the form of the waves based on trigonometric functions. The function that satisfies this ordinary differential equation is defined as leaf function; sleaf$_n(t)$. By multiplying the derivative $dx/dt$ with respect to Eq. (\ref{3.1.1}), the following equation is obtained.

\begin{equation}
\frac{\mathrm{d}^2x }{\mathrm{d}t^2} \frac{\mathrm{d}x }{\mathrm{d}t}=-nx^{2n-1}\frac{\mathrm{d}x }{\mathrm{d}t} \label{3.1.4}
\end{equation}

The following equation is obtained by integrating both sides of Eq. (\ref{3.1.4}).

\begin{equation}
\frac{1}{2} \left( \frac{\mathrm{d}x }{\mathrm{d}t} \right)^2=-\frac{ 1 }{ 2 } x^{2n} +C\label{3.1.5}
\end{equation}

Using the initial conditions in Eqs. (\ref{3.1.2}) and (\ref{3.1.3}), the constant $C = \frac{1}{2}$ is determined. The following equation is obtained by solving the derivative $dx/dt$ in Eq. (\ref{3.1.5}).

\begin{equation}
\frac{\mathrm{d}x }{\mathrm{d}t} = \pm \sqrt{1- x^{2n}} \label{3.1.6}
\end{equation}

We can create a graph with the horizontal axis as the variable $t$ and the vertical axis as the function $x$. Because function $x$ is a wave with a period, the gradient $dx/dt$ has positive and negative values, and it depends on domain $t$.  The relations between variable $t$ and the variable $x$ are sumarized in Table \ref{tab3.1.1}. The variable $m$ represents an integer.

\begin{table}
\begin{center}
\caption{ Relation between the variable $x$ and $t$ of leaf function; sleaf$_n(t)$ (See Appendix Q for the constant $\pi_n$) }
\label{tab3.1.1}
\begin{tabular}{ccccc}
\hline\noalign{\smallskip}
Number   & Domain of variable $t$   & variable $t$ & Range of  & Derivation  \\

 in Fig. \ref{fig3.1.1}  &   &  & variable $x$ &  d$x$/d$t$ \\
\noalign{\smallskip}\hline\noalign{\smallskip}
(1) &	$(2m-2) \pi_n \leqq t < (2m-\frac{3}{2}) \pi_n$ 
     &  $t=(2m-2) \pi_n  $
     &	$0 \leqq x \leqq 1$ 
     &	$\frac{\mathrm{d}x}{\mathrm{d}t}=\sqrt{1-x^{2n}}  $  \\
     &	
     &  $ + \int_{0}^{x}  \frac{\mathrm{d}u }{ \sqrt{1- u^{2n}}} $
     &	
     &	  \\
(2) &	$(2m-\frac{3}{2}) \pi_n \leqq t < (2m-1) \pi_n$
& $t=(2m-\frac{3}{2}) \pi_n $
&	$0 \leqq x \leqq 1$ 
&	$\frac{\mathrm{d}x}{\mathrm{d}t}=-\sqrt{1-x^{2n}}  $ \\
 &	
& $+ \int_{x}^{1}  \frac{\mathrm{d}u }{ \sqrt{1- u^{2n}}} $
&	 
&	\\
(3) &	$(2m-1) \pi_n \leqq t < (2m-\frac{1}{2}) \pi_n$
& $t=(2m-1) \pi_n  $
&	$-1 \leqq x \leqq 0$
&	$\frac{\mathrm{d}x}{\mathrm{d}t}=-\sqrt{1-x^{2n}}  $ \\
  &	
& $ + \int_{x}^{0}  \frac{\mathrm{d}u }{ \sqrt{1- u^{2n}}} $
&	
&	 \\
(4) &	$(2m-\frac{1}{2}) \pi_n \leqq t < 2m \pi_n$
& $t=(2m-\frac{1}{2}) \pi_n  $
&	$-1 \leqq x \leqq 0$ 
&	$\frac{\mathrm{d}x}{\mathrm{d}t}=\sqrt{1-x^{2n}}  $ \\
 &	
& $ + \int_{-1}^{x}  \frac{\mathrm{d}u }{ \sqrt{1- u^{2n}}} $
&	
&	 \\
\noalign{\smallskip}\hline
\end{tabular}
\end{center}
\end{table}

\begin{figure*}[tb]
\begin{center}
\includegraphics[width=0.75 \textwidth]{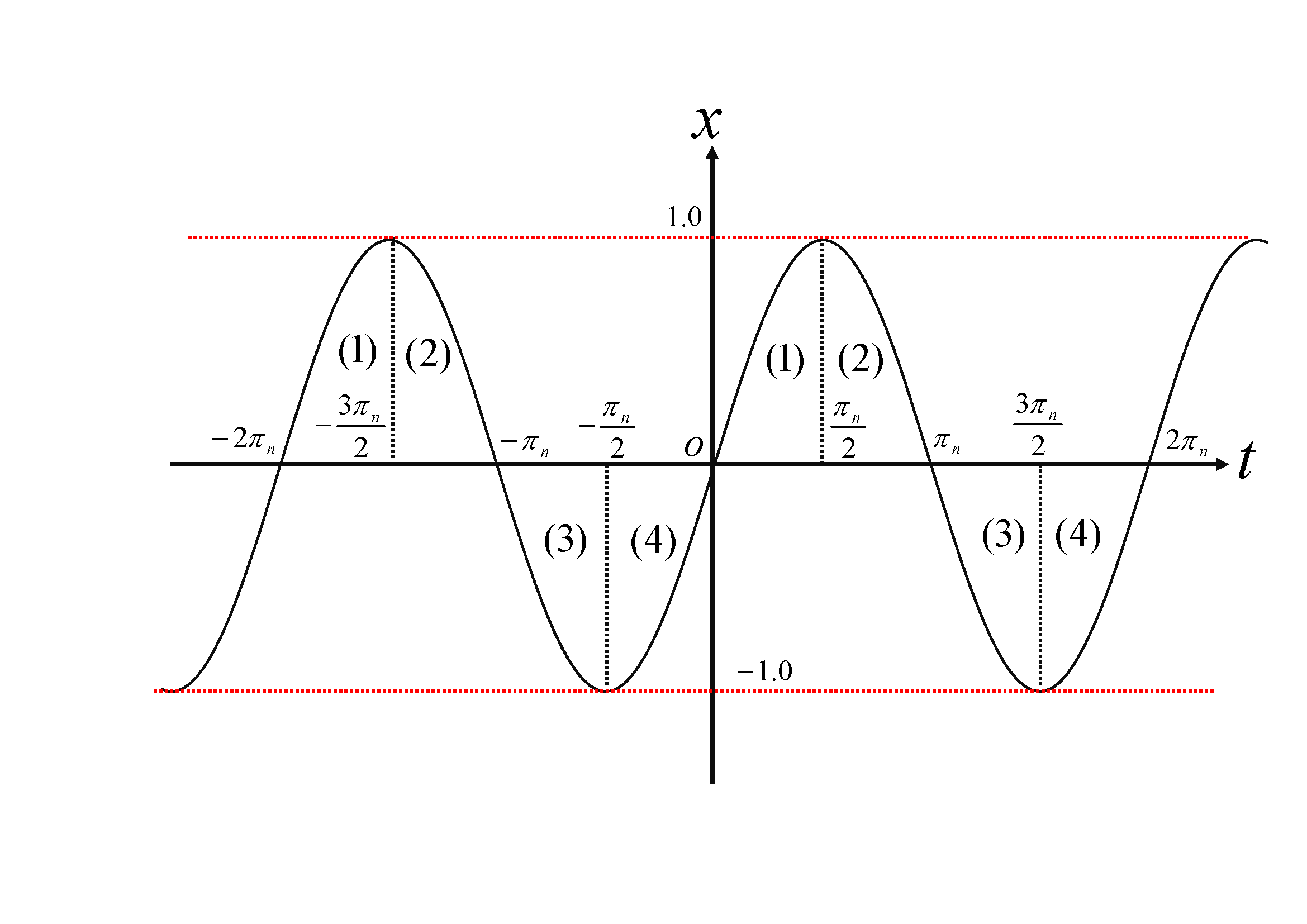}
\caption{Diagram of wave with respect to leaf function; sleaf$_n(t)$  (See Appendix Q for the constant $\pi_n$) }
\label{fig3.1.1}      
\end{center}
\end{figure*}

\subsection{Leaf function; cleaf$_n(t)$}
Another function satisfies the following ordinary differential equations and initial conditions.

\begin{equation}
\frac{\mathrm{d}^2x(t) }{\mathrm{d}t^2}=-nx(t)^{2n-1} \label{3.2.1}
\end{equation}
\begin{equation}
x(0)=1 \label{3.2.2}
\end{equation}
\begin{equation}
\frac{\mathrm{d}x(0)}{\mathrm{d}t}=0
\label{3.2.3}
\end{equation}

The function that satisfies this ordinary differential equation is defined as leaf function; cleaf$_n(t)$. The relations between variable $t$ and the variable $x$ are sumarized in Table \ref{tab3.2.1}.

\begin{table*}
\begin{center}
\caption{ Relation between the variable $x$ and $t$ of leaf function; cleaf$_n(t)$  (See Appendix Q for the constant $\pi_n$)}
\label{tab3.2.1}
\begin{tabular}{ccccc}
\hline\noalign{\smallskip}
Number  & Domain of variable $t$   & variable $t$ & Range of  & Derivation  \\
in Fig. \ref{fig3.2.1}  &    & variable $t$ &  variable $x$ & d$x$/d$t$ \\
\noalign{\smallskip}\hline\noalign{\smallskip}
(1) &	$(2m-2) \pi_n \leqq t < (2m-\frac{3}{2}) \pi_n$ 
     &  $t=(2m-2) \pi_n  $
     &	$0 \leqq x \leqq 1$ 
     &	$\frac{\mathrm{d}x}{\mathrm{d}t}=-\sqrt{1-x^{2n}}  $  \\
& 
     &  $ - \int_{1}^{x}  \frac{\mathrm{d}u }{ \sqrt{1- u^{2n}}} $
     &	 
     &	  \\
(2) &	$(2m-\frac{3}{2}) \pi_n \leqq t < (2m-1) \pi_n$
& $t=(2m-\frac{3}{2}) \pi_n  $
&	$-1 \leqq x \leqq 0$ 
&	$\frac{\mathrm{d}x}{\mathrm{d}t}=-\sqrt{1-x^{2n}}  $ \\
 &	
& $ - \int_{0}^{x}  \frac{\mathrm{d}u }{ \sqrt{1- u^{2n}}} $
&	
&	 \\
(3) &	$(2m-1) \pi_n \leqq t < (2m-\frac{1}{2}) \pi_n$
& $t=(2m-1) \pi_n  $
&	$-1 \leqq x \leqq 0$
&	$\frac{\mathrm{d}x}{\mathrm{d}t}=\sqrt{1-x^{2n}}  $ \\
  &	
& $ + \int_{-1}^{x}  \frac{\mathrm{d}u }{ \sqrt{1-u^{2n}}} $
&
&	\\
(4) &	$(2m-\frac{1}{2}) \pi_n \leqq t < 2m \pi_n$
& $t=(2m-\frac{1}{2}) \pi_n  $
&	$0 \leqq x \leqq 1$ 
&	$\frac{\mathrm{d}x}{\mathrm{d}t}=\sqrt{1-x^{2n}}  $ \\
 &	
& $ + \int_{0}^{x}  \frac{\mathrm{d}u }{ \sqrt{1- u^{2n}}} $
&	
&	 \\
\noalign{\smallskip}\hline
\end{tabular}
\end{center}
\end{table*}

\begin{figure*}[tb]
\begin{center}
\includegraphics[width=0.75 \textwidth]{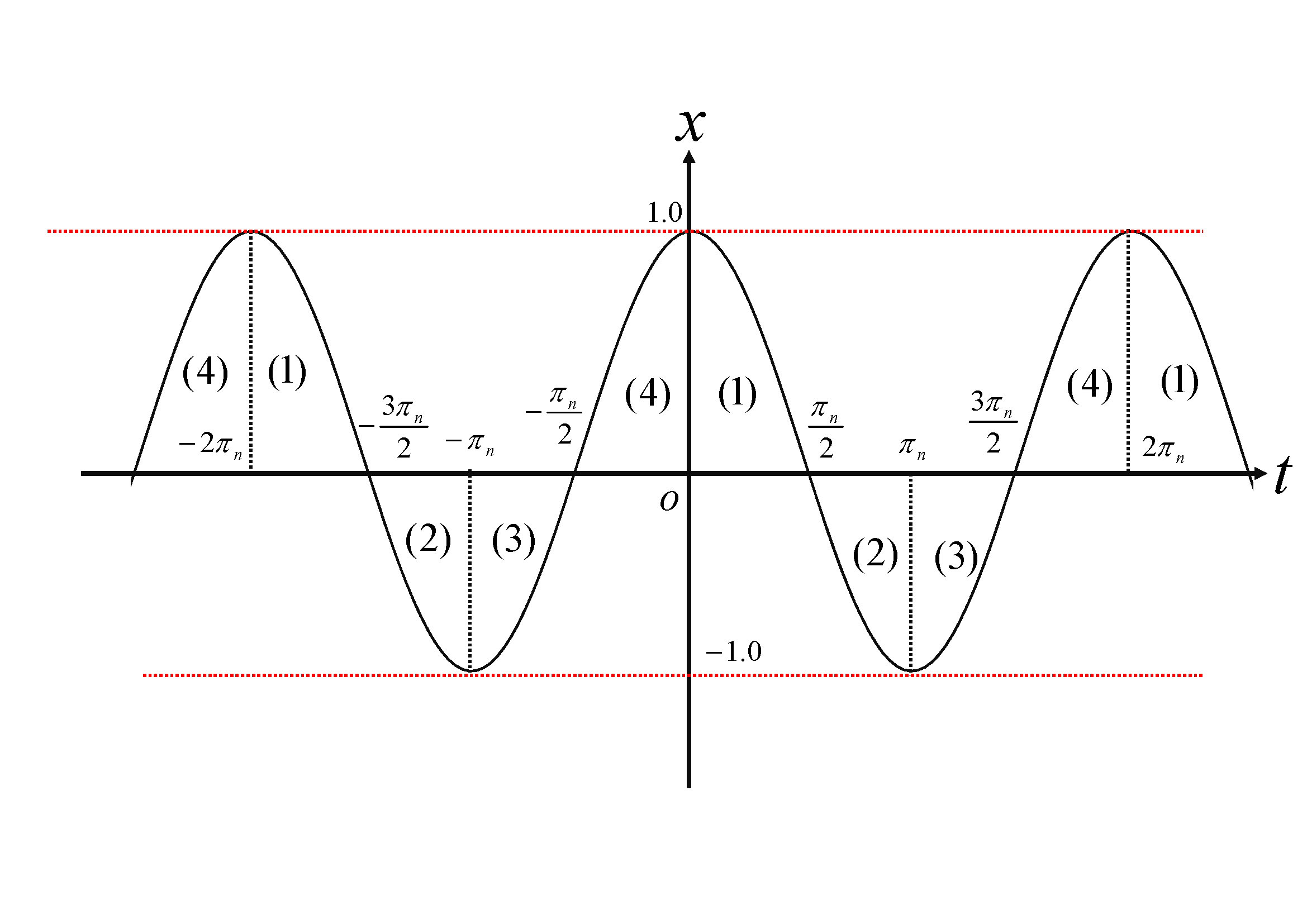}
\caption{Diagram of wave with respect to leaf function; cleaf$_n(l)$  (See Appendix Q for the constant $\pi_n$)}
\label{fig3.2.1}      
\end{center}
\end{figure*}

\subsection{Hyperbolic Leaf function; sleafh$_n(t)$}
On the right side of Eq. (\ref{3.1.1}) and Eq. (\ref{3.2.1}), if the minus changes to plus, the waves disappear with respect to any arbitrary $n$. The variable $x(t)$ increases monotonically as the variable $t$ increases. These ordinary differential equations are as follows:
\begin{equation}
\frac{\mathrm{d}^2x(t)}{\mathrm{d}t^2}=nx(t)^{2n-1} \label{3.3.1}
\end{equation}
\begin{equation}
x(0)=0
\label{3.3.2}
\end{equation}
\begin{equation}
\frac{\mathrm{d}x(0)}{\mathrm{d}t}=1
\label{3.3.3}
\end{equation}
In this study, a function that satisfies the above equations is defined as $\mathrm{sleafh}_{n}(t)$ \cite{Kaz_slh}. The relations between variable $t$ and the variable $x$ are sumarized in Table \ref{tab3.3.1}. 

\begin{table*}
\begin{center}
\caption{ Relation between the variable $x$ and $t$ of leaf function; sleafh$_n(t)$  (See Appendix R for the constant $\zeta_n$)}
\label{tab3.3.1}
\begin{tabular}{ccccc}
\hline\noalign{\smallskip}
Number   & Domain of variable $t$   & variable $t$ & Range of  & Derivation  \\
in Fig. \ref{fig3.3.1}  &    & variable $t$ &  variable $x$ & d$x$/d$t$ \\
\noalign{\smallskip}\hline\noalign{\smallskip}
(1) &	$ (2m-1)\zeta_n < t < (2m+1)\zeta_n$ 
     &  $t=2m \zeta_n  $
     &	$-\infty < x < \infty $ 
     &	$\frac{\mathrm{d}x}{\mathrm{d}t}=\sqrt{1+t^{2n}}  $  \\
&	
     &  $ + \int_{0}^{x}  \frac{\mathrm{d}u }{ \sqrt{ 1+u^{2n}}} $
     &	
     &	 \\
\noalign{\smallskip}\hline
\end{tabular}
\end{center}
\end{table*}

\begin{figure*}[tb]
\begin{center}
\includegraphics[width=0.75 \textwidth]{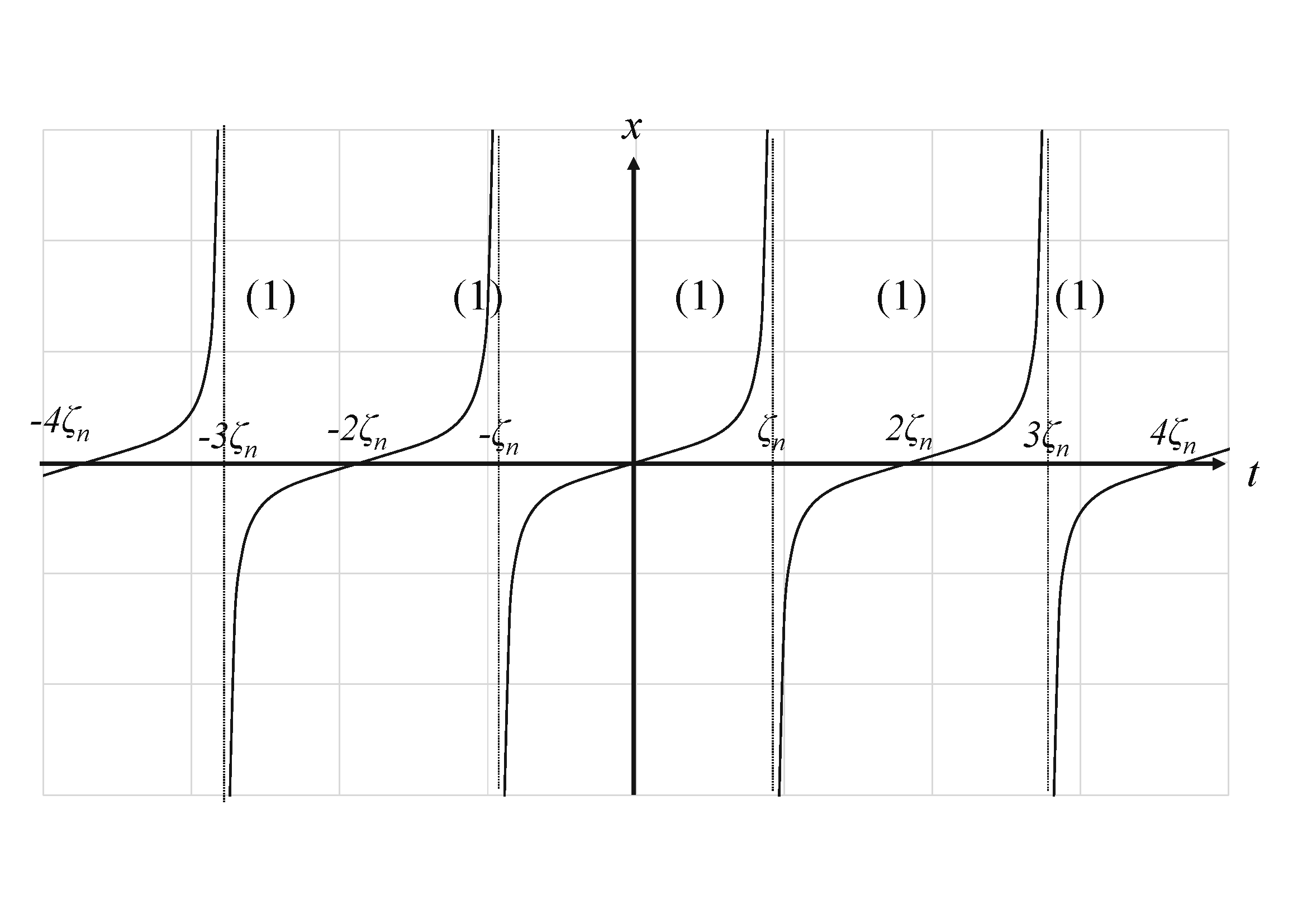}
\caption{Diagram of wave with respect to hyperbolic leaf function; sleafh$_n(t)$ (See Appendix R for the constant $\zeta_n$)}
\label{fig3.3.1}      
\end{center}
\end{figure*}

\subsection{Hyperbolic Leaf function; cleafh$_n(t)$}
Another function satisfies the following ordinary differential equations and initial conditions.

\begin{equation}
\frac{\mathrm{d}^2x(t)}{\mathrm{d}t^2}=nx(t)^{2n-1}
\label{3.4.1}
\end{equation}
\begin{equation}
x(0)=1
\label{3.4.2}
\end{equation}
\begin{equation}
\frac{\mathrm{d}x(0)}{\mathrm{d}t}=0
\label{3.4.3}
\end{equation}

In this study, a function that satisfies the above equations is defined as $\mathrm{cleafh}_{n}(t)$ \cite{Kaz_clh}. The relations between variable $t$ and the variable $x$ are sumarized in Table \ref{tab3.4.1}. 

\begin{table*}
\begin{center}
\caption{ Relation between the variable $x$ and $t$ of leaf function; cleafh$_n(t)$ (See Appendix P for the constant $\eta_n$)}
\label{tab3.4.1}
\begin{tabular}{ccccc}
\hline\noalign{\smallskip}
Number  & Domain of variable $t$   & variable $t$ & Range of & Derivation  \\
in Fig. \ref{fig3.4.1}  &    & variable $t$ &  variable $x$ &  d$x$/d$t$ \\
\noalign{\smallskip}\hline\noalign{\smallskip}
(1) &	$(4m-1) \eta_n < t \leqq 4m \eta_n$ 
     &  $t=4m \eta_n $
     &	$1 \leqq x $ 
     &	$\frac{\mathrm{d}x}{\mathrm{d}t}=-\sqrt{x^{2n}-1}  $  \\
&	 
     &  $- \int_{1}^{x}  \frac{\mathrm{d}u }{ \sqrt{ u^{2n}-1}} $
     &	
     &	    \\
(2) &	$4m \eta_n \leqq t < (4m+1) \eta_n$ 
& $t=4m \eta_n  $
&	$1 \leqq x $ 
&	$\frac{\mathrm{d}x}{\mathrm{d}t}=\sqrt{x^{2n}-1}  $ \\
 &	
& $+ \int_{1}^{x}  \frac{\mathrm{d}u }{ \sqrt{ u^{2n}-1}} $
&	 
&	 \\
(3) & $(4m+1) \eta_n < t \leqq (4m+2) \eta_n$ 
& $t=(4m+2) \eta_n  $
&	$ x \leqq -1$
&	$\frac{\mathrm{d}x}{\mathrm{d}t}=\sqrt{x^{2n}-1}  $ \\
& 
&  $- \int_{x}^{-1}  \frac{\mathrm{d}u }{ \sqrt{ u^{2n}-1}} $
&	
&	 \\
(4) &	$(4m+2) \eta_n \leqq t < (4m+3) \eta_n$ 
& $t=(4m+2) \eta_n  $
&	$ x \leqq -1$
&	$\frac{\mathrm{d}x}{\mathrm{d}t}=-\sqrt{x^{2n}-1}  $ \\
 &	 
& $ + \int_{x}^{-1}  \frac{\mathrm{d}u }{ \sqrt{ u^{2n}-1}} $
&	
&	 \\
\noalign{\smallskip}\hline
\end{tabular}
\end{center}
\end{table*}

\begin{figure*}[tb]
\begin{center}
\includegraphics[width=0.75 \textwidth]{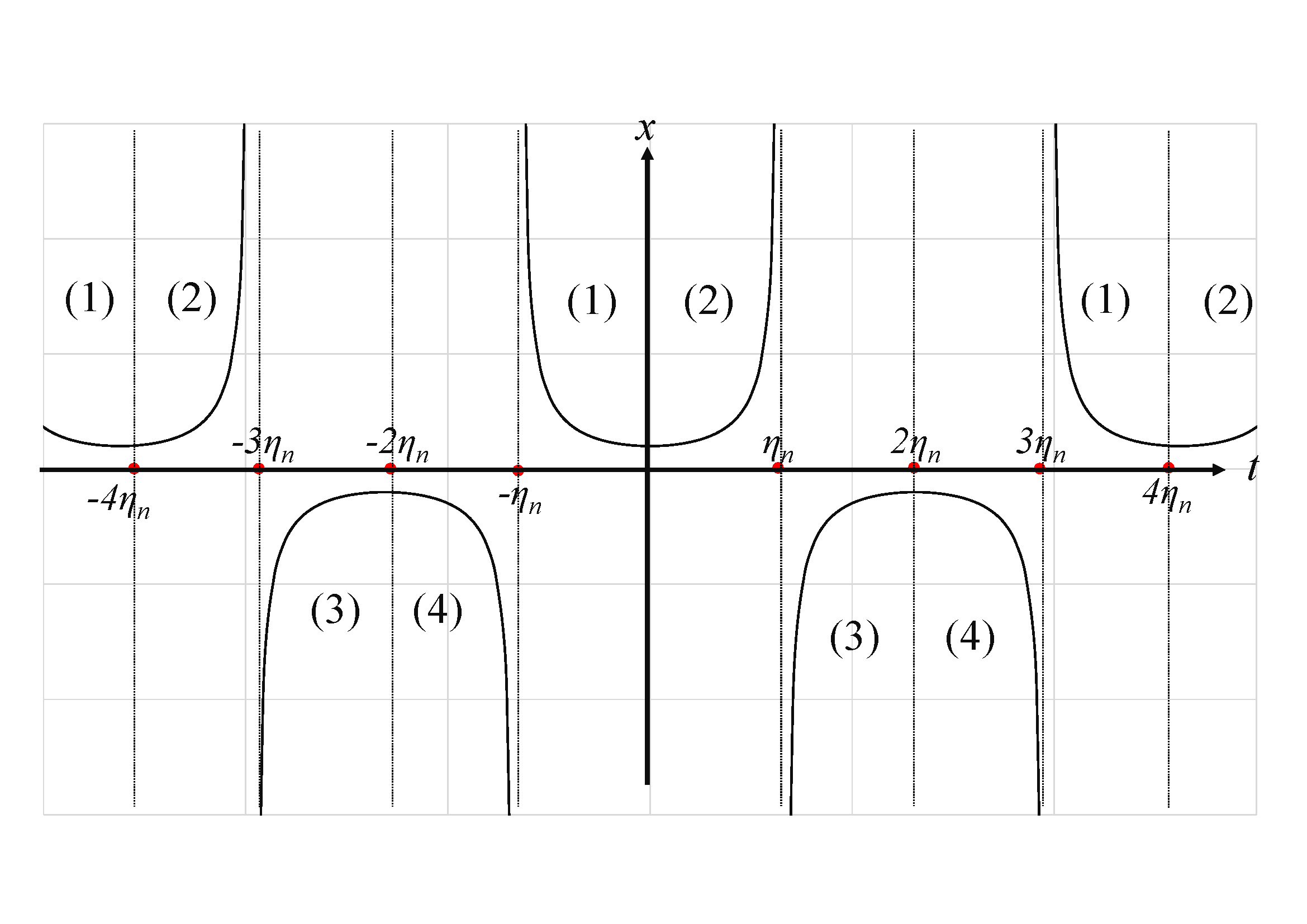}
\caption{Diagram of wave with respect to hyperbolic leaf function; cleafh$_n(t)$ (See Appendix P for the constant $\eta_n$)}
\label{fig3.4.1}      
\end{center}
\end{figure*}


\section{Exact solutions of the cubic Duffing equation}
The exact solutions of the cubic Duffing equation are obtained by using the leaf functions. 
There are eight types of exact solutions.
\subsection{Exact solution 1}
The exact solution is 
\begin{equation}
 x(t)=A \sqrt{1+(\mathrm{cleaf}_2(\omega t))^2}
\label{4.1}
\end{equation}
The ordinary differential equation is 
\begin{equation}
\frac{\mathrm{d}^2 x(t) }{\mathrm{d}t^2} -3 \omega^2 x(t) + 2 \left( \frac{\omega}{A} \right)^2 x(t)^3=0
\label{4.2}
\end{equation}

\subsection{Exact solution 2}
The exact solution is 
\begin{equation}
 x(t)=A \sqrt{1-(\mathrm{cleaf}_2( \omega t))^2}
\label{4.3}
\end{equation}
The ordinary differential equation is 
\begin{equation}
\frac{\mathrm{d}^2 x(t) }{\mathrm{d}t^2} +3 \omega^2 x(t) - 2 \left( \frac{\omega}{A} \right)^2 x(t)^3=0
\label{4.4}
\end{equation}

\subsection{Exact solution 3}
The exact solution is 
\begin{equation}
x(t)=A \sqrt{1+(\mathrm{sleaf}_2( \omega t))^2}
\label{4.5}
\end{equation}

The ordinary differential equation that satisfies the above exact solution is the same as Eq. (\ref{4.2}). 

\subsection{Exact solution 4}
The exact solution is 
\begin{equation}
 x(t)=A \sqrt{1 - (\mathrm{sleaf}_2( \omega t))^2}
\label{4.6}
\end{equation}

The ordinary differential equation that satisfies the above exact solution is the same as Eq. (\ref{4.4}). 

\subsection{Exact solution 5}
The exact solution is 
\begin{equation}
 x(t)=A \sqrt{1 + (\mathrm{sleaf}_2( \omega t))^2} + A \sqrt{1 + (\mathrm{cleaf}_2( \omega t))^2}
\label{4.7}
\end{equation}
The ordinary differential equation is 
\begin{equation}
\frac{\mathrm{d}^2 x(t) }{\mathrm{d}t^2} - 3 \omega^2 ( 2\sqrt{2}+1 ) x(t) + 2 \left( \frac{\omega}{A} \right)^2 x(t)^3=0 \label{4.8}
\end{equation}

\subsection{Exact solution 6}
The exact solution is 
\begin{equation}
 x(t)=A \sqrt{1 + (\mathrm{cleaf}_2( \omega t))^2} -
 A \sqrt{1 + (\mathrm{sleaf}_2( \omega t))^2} 
\label{4.9}
\end{equation}
The ordinary differential equation is 
\begin{equation}
\frac{\mathrm{d}^2 x(t) }{\mathrm{d}t^2} + 3 \omega^2 ( 2\sqrt{2} - 1 ) x(t) + 2 \left( \frac{\omega}{A} \right)^2 x(t)^3=0 \label{4.10}
\end{equation}

\subsection{Exact solution 7}
The exact solution is 
\begin{equation}
 x(t)=A \sqrt{ (\mathrm{cleafh}_2( \omega t))^2+1}
\label{4.11}
\end{equation}

The ordinary differential equation that satisfies the above exact solution is the same as Eq. (\ref{4.4}). 

\subsection{Exact solution 8}
The exact solution is 
\begin{equation}
 x(t)=A \sqrt{ (\mathrm{cleafh}_2( \omega t))^2-1}.
\label{4.12}
\end{equation}
The ordinary differential equation is 
\begin{equation}
\frac{\mathrm{d}^2 x(t) }{\mathrm{d}t^2} -3 \omega^2 x(t) - 2 \left( \frac{\omega}{A} \right)^2 x(t)^3=0
\label{4.13}.
\end{equation}


\section{Exact solutions of the cubic--quintic Duffing equation}
The exact solutions of the cubic--quintic Duffing equation are obtained by using the leaf functions. 
There are six types of exact solutions.
\subsection{Exact solution 9}
The exact solution is 
\begin{equation}
 x(t)=A \sqrt{1+ \mathrm{sleaf}_2(\omega t) }.
\label{5.1}
\end{equation}
The ordinary differential equation is 
\begin{equation}
\frac{\mathrm{d}^2 x(t) }{\mathrm{d}t^2} + \frac{3}{2} \omega^2 x(t) - 2  \frac{\omega^2}{A^2}  x(t)^3 + \frac{3}{4}  \frac{\omega^2}{A^4}  x(t)^5=0
\label{5.2}.
\end{equation}

\subsection{Exact solution 10}
The exact solution is 
\begin{equation}
 x(t)=A \sqrt{1- \mathrm{sleaf}_2(\omega t) }.
\label{5.3}
\end{equation}

The ordinary differential equation that satisfies the above exact solution is the same as Eq. (\ref{5.2}).

\subsection{Exact solution 11}
The exact solution is 
\begin{equation}
 x(t)=A \sqrt{1+ \mathrm{cleaf}_2(\omega t) }.
\label{5.4}
\end{equation}
The ordinary differential equation that satisfies the above exact solution is the same as Eq. (\ref{5.2}).

\subsection{Exact solution 12}
The exact solution is 
\begin{equation}
 x(t)=A \sqrt{1- \mathrm{cleaf}_2(\omega t) .}
\label{5.5}
\end{equation}
The ordinary differential equation that satisfies the above exact solution is the same as Eq. (\ref{5.2}).

\subsection{Exact solution 13}
The exact solution is 
\begin{equation}
 x(t)=A \sqrt{\mathrm{cleafh}_2(\omega t)+1 }.
\label{5.6}
\end{equation}
The ordinary differential equation is 
\begin{equation}
\frac{\mathrm{d}^2 x(t) }{\mathrm{d}t^2} - \frac{3}{2} \omega^2 x(t) + 2  \frac{\omega^2}{A^2}  x(t)^3 - \frac{3}{4}  \frac{\omega^2}{A^4}  x(t)^5=0
\label{5.7}.
\end{equation}

\subsection{Exact solution 14}
The exact solution is 
\begin{equation}
 x(t)=A \sqrt{\mathrm{cleafh}_2(\omega t)-1 }.
\label{5.8}
\end{equation}
The ordinary differential equation is 
\begin{equation}
\frac{\mathrm{d}^2 x(t) }{\mathrm{d}t^2} - \frac{3}{2} \omega^2 x(t) - 2  \frac{\omega^2}{A^2}  x(t)^3 - \frac{3}{4}  \frac{\omega^2}{A^4}  x(t)^5=0
\label{5.9}.
\end{equation}

\section{Extended exact solutions of cubic Duffing equatoin }
\subsection{Theory of Transformation}
The exact solution under the free vibration is proposed in Section 5. These solutions can be extended to the exact solution under damped vibration \cite{Gallegos}.
A suitable chage of variables can leads to a simplification of the problem. The damped duffing equation can be decribed as follows:

\begin{equation}
\frac{\mathrm{d}^2y(t)}{\mathrm{d}t^2} + \beta(t) \frac{\mathrm{d}y(t)}{\mathrm{d}t}+\alpha_1(t) y(t)+ h(t,y(t))=0. \label{6.1}
\end{equation}

The function $h(t,y(t))$ adopt the cubic form or the cubic-quintic form. 

\begin{equation}
h(t,y(t))=\alpha_2(t) y(t)^3.
\label{6.2}
\end{equation}

\begin{equation}
h(t,y(t))=\alpha_2(t) y(t)^3+\alpha_3(t) y(t)^5.
\label{6.3}
\end{equation}

The damped duffing equation can be transformed to the undamped duffing equation by introducing the following function.

\begin{equation}
x(t)= y(t) \mathrm{exp} \left(\frac{1}{2} \int_{c}^{t} \beta(u) du \right)
\label{6.4}
\end{equation}

The variable $c$ represents the constant. Using the above equation and Eq. (\ref{6.1}), the following equation can be derived.

\begin{equation}
\begin{split}
&\frac{\mathrm{d}^2 x(t)}{\mathrm{d}t^2} + 
\left( \alpha_1(t)-\frac{1}{4} \beta(t)^2
-\frac{1}{2} \frac{\mathrm{d} \beta(t)}{\mathrm{d}t}\right)x(t) \\
&+\mathrm{exp} \left(\frac{1}{2} \int_{c}^{t} \beta(u) du \right)
\cdot h\left( t,x(t) \mathrm{exp} \left( -\frac{1}{2} \int_{c}^{t} \beta(u) du \right) \right)=0.
\label{6.5}
\end{split}
\end{equation}

The created equation is simpler than the original one in the sense that the term of damping is not present any more \cite{Gallegos}. 
By comparing the coefficient of the term of the above equation with the coefficient of the term of the undamped Duffing equation, it becomes possible to create a new exact solution that satisfies the damped Duffing equation by leaf functions.


\subsection{Extended exact solution 1}
Based on both the section 4.1 and the section 6.1, the exact solution 1 under free vibration is extended to the exact solution under damped vibration. The damped vibration of exact solution 1 is derived based on the Eq. (\ref{4.2}) and the Eq. (\ref{6.5}) . The exact solution 1 under the damped vibration is 
\begin{equation}
y(t)=A \cdot \mathrm{exp} \left(-\frac{1}{2} \int_{c}^{t} \beta(u) du \right) \cdot \sqrt{1+(\mathrm{cleaf}_2(\omega t))^2}. \label{6.6}
\end{equation}

The ordinary differential equation is 
\begin{equation}
\frac{\mathrm{d}^2 y(t) }{\mathrm{d}t^2}+ \beta(t) \frac{\mathrm{d}y(t)}{\mathrm{d}t}
+\left( -3 \omega^2 + \frac{1}{4} \beta(t)^2 + \frac{1}{2} \frac{\mathrm{d} \beta(t)}{\mathrm{d}t} \right)y(t)
+ 2 \left( \frac{\omega}{A} \right)^2 \mathrm{exp} \left( \int_{c}^{t} \beta(u) du \right) y(t)^3=0
\label{6.7}.
\end{equation}


\subsection{Extended exact solution 2}
The exact solution 2 under the damped vibration is 
\begin{equation}
y(t)=A \cdot \mathrm{exp} \left(-\frac{1}{2} \int_{c}^{t} \beta(u) du \right) \cdot \sqrt{1-(\mathrm{cleaf}_2(\omega t))^2}. \label{6.8}
\end{equation}

The ordinary differential equation is 
\begin{equation}
\frac{\mathrm{d}^2 y(t) }{\mathrm{d}t^2}+ \beta(t) \frac{\mathrm{d}y(t)}{\mathrm{d}t}
+\left( 3 \omega^2 + \frac{1}{4} \beta(t)^2 + \frac{1}{2} \frac{\mathrm{d} \beta(t)}{\mathrm{d}t} \right)y(t)
- 2 \left( \frac{\omega}{A} \right)^2 \mathrm{exp} \left( \int_{c}^{t} \beta(u) du \right) y(t)^3=0
\label{6.9}.
\end{equation}


\subsection{Extended exact solution 3}
The exact solution 3 under the damped vibration is 
\begin{equation}
y(t)=A \cdot \mathrm{exp} \left(-\frac{1}{2} \int_{c}^{t} \beta(u) du \right) \cdot \sqrt{1+(\mathrm{sleaf}_2(\omega t))^2}. \label{6.10}
\end{equation}

The ordinary differential equation that satisfies the above exact solution is the same as Eq. (\ref{6.7}).


\subsection{Extended exact solution 4}
The exact solution 4 under the damped vibration is 
\begin{equation}
y(t)=A \cdot \mathrm{exp} \left(-\frac{1}{2} \int_{c}^{t} \beta(u) du \right) \cdot \sqrt{1-(\mathrm{sleaf}_2(\omega t))^2}. \label{6.11}
\end{equation}

The ordinary differential equation that satisfies the above exact solution is the same as Eq. (\ref{6.9}).


\subsection{Extended exact solution 5}
The exact solution 5 under the damped vibration is 
\begin{equation}
y(t)=A \cdot \mathrm{exp} \left(-\frac{1}{2} \int_{c}^{t} \beta(u) du \right) \cdot \left(
\sqrt{1 + (\mathrm{sleaf}_2( \omega t))^2} + \sqrt{1 + (\mathrm{cleaf}_2( \omega t))^2} \right)
\label{6.12}
\end{equation}

The ordinary differential equation is 
\begin{equation}
\begin{split}
& \frac{\mathrm{d}^2 y(t) }{\mathrm{d}t^2}+ \beta(t) \frac{\mathrm{d}y(t)}{\mathrm{d}t}
+\left(  -3 \omega^2 ( 2\sqrt{2}+1 ) + \frac{1}{4} \beta(t)^2 + \frac{1}{2} \frac{\mathrm{d} \beta(t)}{\mathrm{d}t} \right)y(t) \\
&+ 2 \left( \frac{\omega}{A} \right)^2 \mathrm{exp} \left( \int_{c}^{t} \beta(u) du \right) y(t)^3=0
\label{6.13}.
\end{split}
\end{equation}


\subsection{Extended exact solution 6}
The exact solution 6 under the damped vibration is 
\begin{equation}
y(t)=A \cdot \mathrm{exp} \left(-\frac{1}{2} \int_{c}^{t} \beta(u) du \right) \cdot \left(
\sqrt{1 + (\mathrm{cleaf}_2( \omega t))^2}-
\sqrt{1 + (\mathrm{sleaf}_2( \omega t))^2}  \right)
\label{6.14}
\end{equation}

The ordinary differential equation is 
\begin{equation}
\begin{split}
&\frac{\mathrm{d}^2 y(t) }{\mathrm{d}t^2}+ \beta(t) \frac{\mathrm{d}y(t)}{\mathrm{d}t}
+\left(  3 \omega^2 ( 2\sqrt{2}-1 ) + \frac{1}{4} \beta(t)^2 + \frac{1}{2} \frac{\mathrm{d} \beta(t)}{\mathrm{d}t} \right)y(t) \\
&+ 2 \left( \frac{\omega}{A} \right)^2 \mathrm{exp} \left( \int_{c}^{t} \beta(u) du \right) y(t)^3=0
\label{6.15}.
\end{split}
\end{equation}


\subsection{Extended exact solution 7}
The extended exact solution 7 is 
\begin{equation}
y(t)=A \cdot \mathrm{exp} \left(-\frac{1}{2} \int_{c}^{t} \beta(u) du \right) \cdot \sqrt{(\mathrm{cleafh}_2(\omega t))^2+1}. \label{6.16}
\end{equation}

The ordinary differential equation that satisfies the above exact solution is the same as Eq. (\ref{6.9}).


\subsection{Extended exact solution 8}
The extended exact solution 8 is 
\begin{equation}
y(t)=A \cdot \mathrm{exp} \left(-\frac{1}{2} \int_{c}^{t} \beta(u) du \right) \cdot \sqrt{(\mathrm{cleafh}_2(\omega t))^2-1}. \label{6.17}
\end{equation}

The ordinary differential equation is 
\begin{equation}
\frac{\mathrm{d}^2 y(t) }{\mathrm{d}t^2}+ \beta(t) \frac{\mathrm{d}y(t)}{\mathrm{d}t}
+\left(- 3 \omega^2 + \frac{1}{4} \beta(t)^2 + \frac{1}{2} \frac{\mathrm{d} \beta(t)}{\mathrm{d}t} \right)y(t)
- 2 \left( \frac{\omega}{A} \right)^2 \mathrm{exp} \left( \int_{c}^{t} \beta(u) du \right) y(t)^3=0
\label{6.18}.
\end{equation}


\section{Extended exact solutions of cubic-quintic Duffing equatoin }
\subsection{Extended exact solution 9}
The extended exact solution 9 is 
\begin{equation}
y(t)=A \cdot \mathrm{exp} \left(-\frac{1}{2} \int_{c}^{t} \beta(u) du \right) \cdot \sqrt{1+\mathrm{sleaf}_2(\omega t)}.
\label{7.1}
\end{equation}

The ordinary differential equation is 
\begin{equation}
\begin{split}
&\frac{\mathrm{d}^2 y(t) }{\mathrm{d}t^2}
+ \beta(t) \frac{\mathrm{d}y(t)}{\mathrm{d}t}
+\left( \frac{3}{2} \omega^2 + \frac{1}{4} \beta(t)^2 
+\frac{1}{2} \frac{\mathrm{d} \beta(t)}{\mathrm{d}t} \right)y(t) \\
&-2 \left( \frac{\omega}{A} \right)^2 
\mathrm{exp} \left( \int_{c}^{t} \beta(u) du \right) y(t)^3 
+ \frac{3}{4}  \frac{\omega^2}{A^4}  
\mathrm{exp} \left( 2 \int_{c}^{t} \beta(u) du \right) y(t)^5=0
\label{7.2}.
\end{split}
\end{equation}

\subsection{Extended exact solution 10}
The extended exact solution 10 is 
\begin{equation}
y(t)=A \cdot \mathrm{exp} \left(-\frac{1}{2} \int_{c}^{t} \beta(u) du \right) \cdot \sqrt{1-\mathrm{sleaf}_2(\omega t)}.
\label{7.3}
\end{equation}

The ordinary differential equation that satisfies the above exact solution is the same as Eq. (\ref{7.2}). 

\subsection{Extended exact solution 11}
The extended exact solution 11 is 
\begin{equation}
y(t)=A \cdot \mathrm{exp} \left(-\frac{1}{2} \int_{c}^{t} \beta(u) du \right) \cdot \sqrt{1+\mathrm{cleaf}_2(\omega t)}. 
\label{7.4}
\end{equation}

The ordinary differential equation that satisfies the above exact solution is the same as Eq. (\ref{7.2}). 

\subsection{Extended exact solution 12}
The extended exact solution 12 is 
\begin{equation}
y(t)=A \cdot \mathrm{exp} \left(-\frac{1}{2} \int_{c}^{t} \beta(u) du \right) \cdot \sqrt{1-\mathrm{cleaf}_2(\omega t)}. 
\label{7.5}
\end{equation}

The ordinary differential equation that satisfies the above exact solution is the same as Eq. (\ref{7.2}). 

\subsection{Extended exact solution 13}
The extended exact solution 13 is 
\begin{equation}
y(t)=A \cdot \mathrm{exp} \left(-\frac{1}{2} \int_{c}^{t} \beta(u) du \right) \cdot \sqrt{\mathrm{cleafh}_2(\omega t)+1}. \label{7.6}
\end{equation}

The ordinary differential equation is 
\begin{equation}
\begin{split}
&\frac{\mathrm{d}^2 y(t) }{\mathrm{d}t^2}
+ \beta(t) \frac{\mathrm{d}y(t)}{\mathrm{d}t}
+\left( -\frac{3}{2} \omega^2 + \frac{1}{4} \beta(t)^2 
+\frac{1}{2} \frac{\mathrm{d} \beta(t)}{\mathrm{d}t} \right)y(t) \\
&+2 \left( \frac{\omega}{A} \right)^2 
\mathrm{exp} \left( \int_{c}^{t} \beta(u) du \right) y(t)^3 
- \frac{3}{4}  \frac{\omega^2}{A^4}  
\mathrm{exp} \left( 2 \int_{c}^{t} \beta(u) du \right) y(t)^5=0
\label{7.7}.
\end{split}
\end{equation}

\subsection{Extended exact solution 14}
The extended exact solution 14 is 
\begin{equation}
y(t)=A \cdot \mathrm{exp} \left(-\frac{1}{2} \int_{c}^{t} \beta(u) du \right) \cdot \sqrt{\mathrm{cleafh}_2(\omega t)-1}. \label{7.8}
\end{equation}

The ordinary differential equation is 
\begin{equation}
\begin{split}
&\frac{\mathrm{d}^2 y(t) }{\mathrm{d}t^2}
+ \beta(t) \frac{\mathrm{d}y(t)}{\mathrm{d}t}
+\left( -\frac{3}{2} \omega^2 + \frac{1}{4} \beta(t)^2 
+\frac{1}{2} \frac{\mathrm{d} \beta(t)}{\mathrm{d}t} \right)y(t) \\
&-2 \left( \frac{\omega}{A} \right)^2 
\mathrm{exp} \left( \int_{c}^{t} \beta(u) du \right) y(t)^3 
- \frac{3}{4}  \frac{\omega^2}{A^4}  
\mathrm{exp} \left( 2 \int_{c}^{t} \beta(u) du \right) y(t)^5=0
\label{7.9}.
\end{split}
\end{equation}


\section{Visualization of exact solutions and extended exact solutions }
The graphs in the section are the waveforms or the divergence obtained by numerical analysis based on differential equations and initial conditions. The waveforms or divergences can also be obtained by using the exact solution. The waveforms or the divergences with respect to the exact solutions are plotted for each graph. The vertical and horizontal axes  represent the variables $x(t)$ and $t$, respectively.

\begin{figure*}[tb]
\begin{center}
\includegraphics[width=0.7 \textwidth]{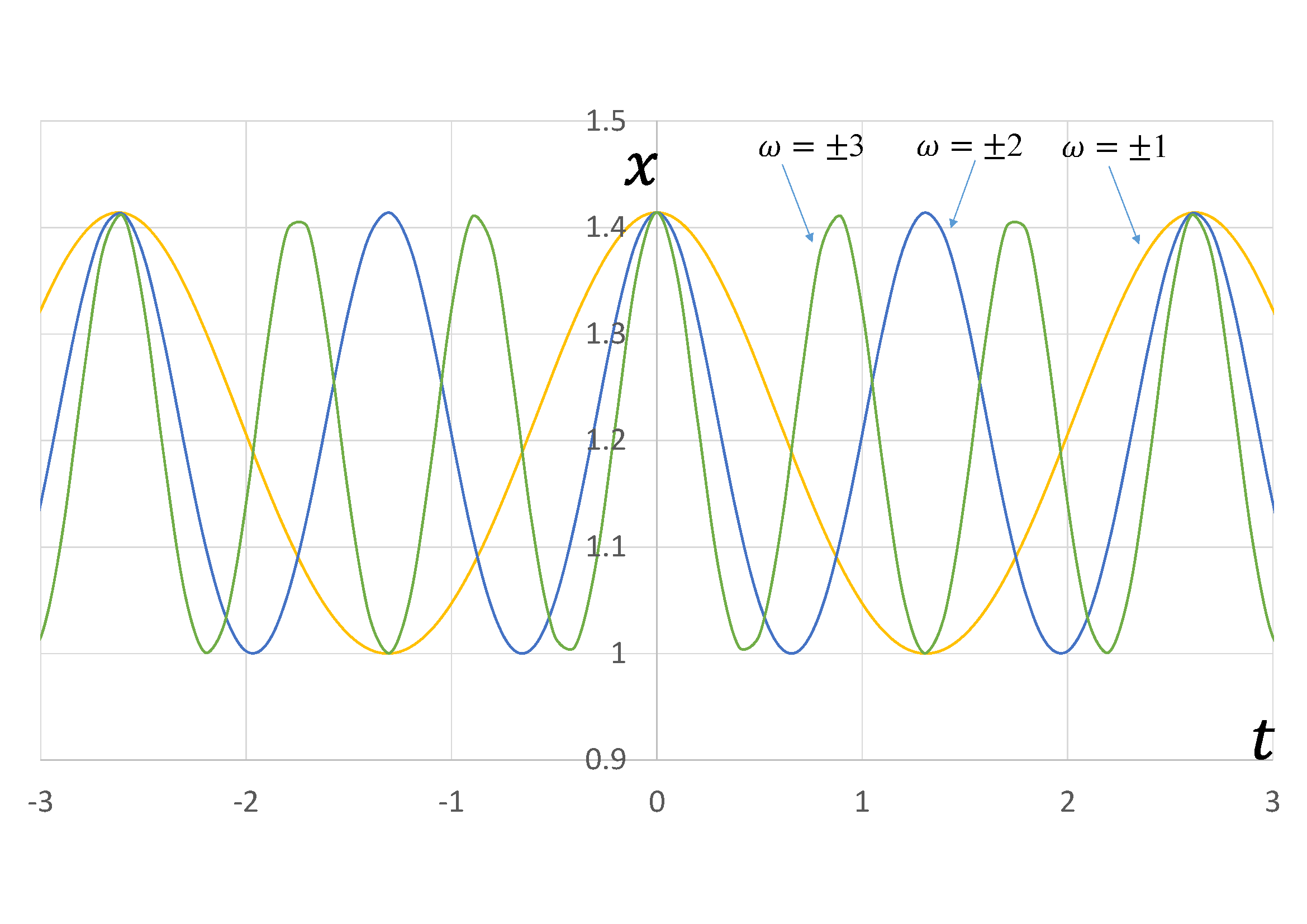}
\caption{Waves of exact solution 1 with respect to variation of the parameter $\omega$ related to the period ($A=1$)}
\label{fig81}      
\end{center}
\end{figure*}

Fig. \ref{fig81} (exact solution 1), Fig. \ref{fig83} (exact solution 2), Fig. \ref{fig85} (exact solution 3), Fig. \ref{fig87} (exact solution 4), Fig. \ref{fig89} (exact solution 5), Fig. \ref{fig811} (exact solution 6), Fig. \ref{fig813} (exact solution 7), Fig. \ref{fig815} (exact solution 8),
Fig. \ref{fig817} (exact solution 9),
Fig. \ref{fig819} (exact solution 10),
Fig. \ref{fig821} (exact solution 11),
Fig. \ref{fig823} (exact solution 12),
Fig. \ref{fig825} (exact solution 13), and
Fig. \ref{fig827} (exact solution 14),
 show the graphs when the period parameter $\omega$ changes, with the amplitude parameter set to $A=1$. 
From the graphs, one can observe that the amplitude of the wave, $x(t)$, does not change even if the parameter $\omega$ changes. The range of $x(t)$ with respect to time $t$ satisfies the following inequalities: 

(i) For Figs. \ref{fig81} and \ref{fig85} \ (exact solutions 1 and 3),
\begin{equation}
1 \leqq  x(t) \leqq \sqrt{2}.
\label{8.1}
\end{equation}

(ii)  For Figs. \ref{fig83} and \ref{fig87} \ (exact solutions 2 and 4),
\begin{equation}
0 \leqq  x(t) \leqq 1.
\label{8.2}
\end{equation}

(iii)  For Fig. \ref{fig89} \ (exact solution 5),
\begin{equation}
2^{\frac{5}{4}} (\fallingdotseq 2.378 \cdots) \leqq  x(t) \leqq 1+\sqrt{2} (\fallingdotseq 2.414 \cdots).
\label{8.3}
\end{equation}

(iv)  For Fig. \ref{fig811} \ (exact solution 6),
\begin{equation}
-(\sqrt{2}-1) (\fallingdotseq -0.414 \cdots) \leqq  x(t) \leqq \sqrt{2}-1 (\fallingdotseq 0.414 \cdots).
\label{8.4}
\end{equation}

(v)  For Fig. \ref{fig813} \ (exact solution 7),
\begin{equation}
1.0 \leqq  x(t).
\label{8.5}
\end{equation}

(vi) For Figs. \ref{fig815} and \ref{fig827} \ (exact solution 8 and 14),
\begin{equation}
0.0 \leqq  x(t).
\label{8.6}
\end{equation}

(vii) For Figs. \ref{fig817}, \ref{fig819}, \ref{fig821}, and \ref{fig823} (exact solutions 9, 10, 11, and 12),
\begin{equation}
0.0 \leqq  x(t) \leqq \sqrt{2}.
\label{8.7}
\end{equation}

(viii) For Fig. \ref{fig825} \ (exact solution 13),
\begin{equation}
\sqrt{2} \leqq  x(t).
\label{8.8}
\end{equation}

The period of wave $x(t)$ is defined as $T$. The relationships with the parameter $\omega$ are as follows: 

(i) For Figs. \ref{fig81}, \ref{fig83}, \ref{fig85}, \ref{fig87}, and \ref{fig811} \ (exact solutions 1, 2, 3, 4, and 6)
\begin{equation}
T=\frac{ \pi_2 }{ | \omega | }.
\label{8.9}
\end{equation}

(ii)  For Fig. \ref{fig89} \ (exact solution 5),
\begin{equation}
T=\frac{ \pi_2 }{ 2| \omega | }.
\label{8.10}
\end{equation}

(iii)  For Figs. \ref{fig813}, \ref{fig815}, \ref{fig825}, and \ref{fig827}, \ (exact solutions 7, 8, 13, and 14)
\begin{equation}
T=\frac{ 2 \eta_2 }{ | \omega | }.
\label{8.11}
\end{equation}

(iv)  For Figs. \ref{fig817}, \ref{fig819}, \ref{fig821}, and \ref{fig823}, \ (exact solution 9, 10, 11, and 12)
\begin{equation}
T=\frac{ 2 \pi_2 }{ | \omega | }.
\label{8.12}
\end{equation}

\begin{figure*}[tb]
\begin{center}
\includegraphics[width=0.7 \textwidth]{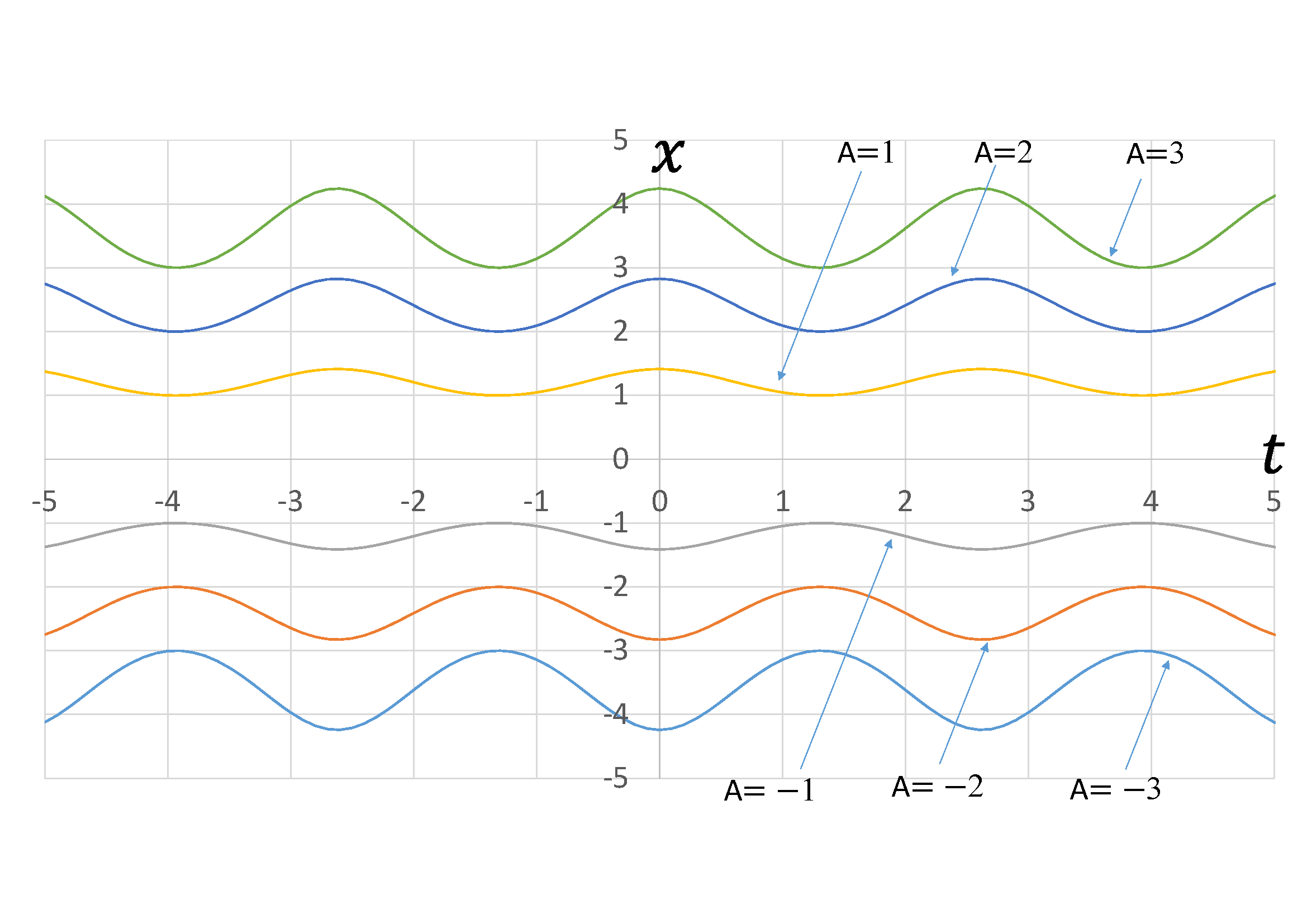}
\caption {Waves of exact solution 1 with respect to variation of the parameter $A$ related to the amplitude ($\omega=1$)}
\label{fig82}      
\end{center}
\end{figure*}

\begin{figure*}[tb]
\begin{center}
\includegraphics[width=0.7 \textwidth]{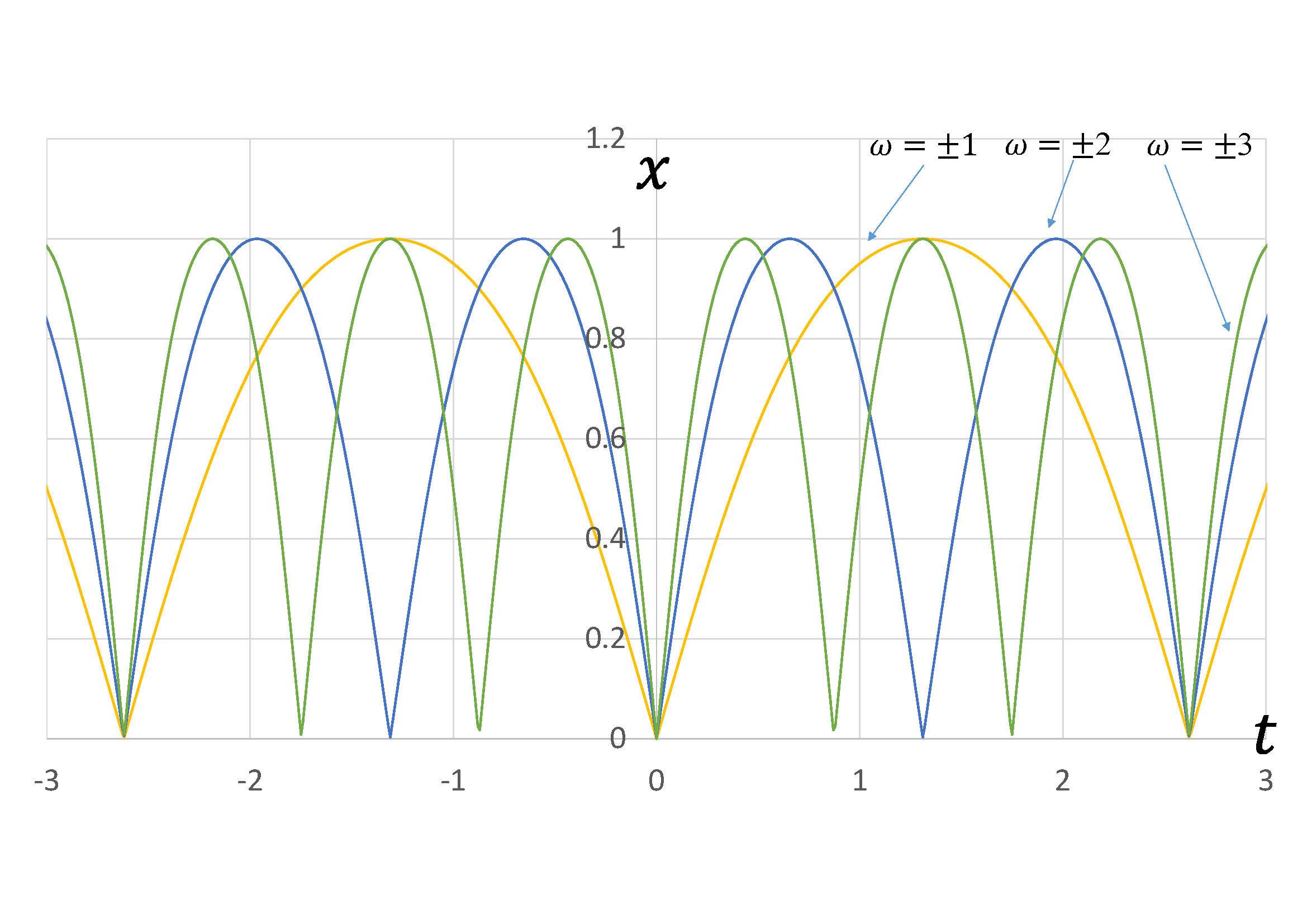}
\caption{Waves of exact solution 2 with respect to variation of the parameter $\omega$ related to the period ($A=1$)}
\label{fig83}      
\end{center}
\end{figure*}

\begin{figure*}[tb]
\begin{center}
\includegraphics[width=0.7 \textwidth]{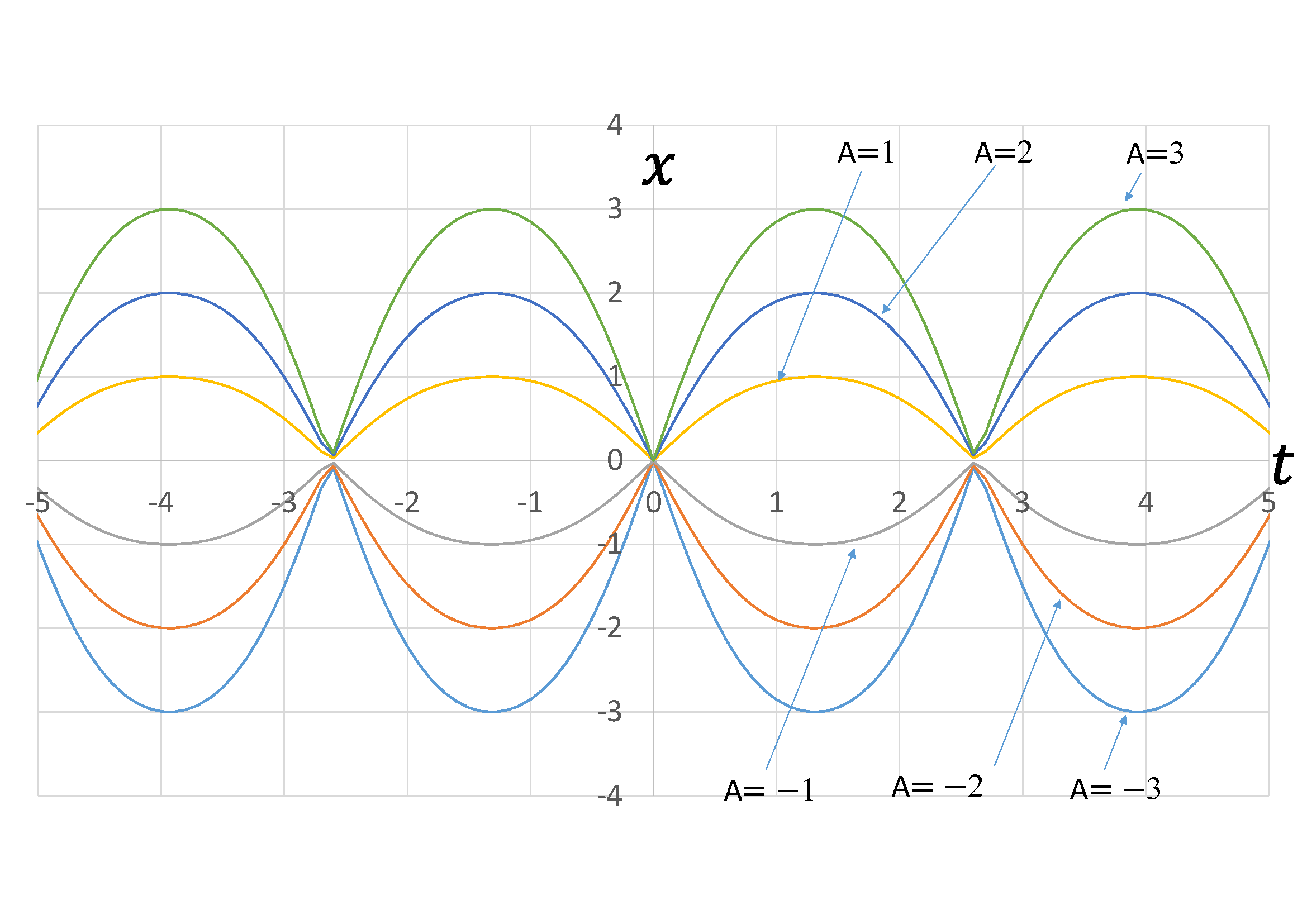}
\caption{Waves of exact solution 2 with respect to variation of the parameter $A$ related to the amplitude ($\omega=1$)}
\label{fig84}      
\end{center}
\end{figure*}

\begin{figure*}[tb]
\begin{center}
\includegraphics[width=0.7 \textwidth]{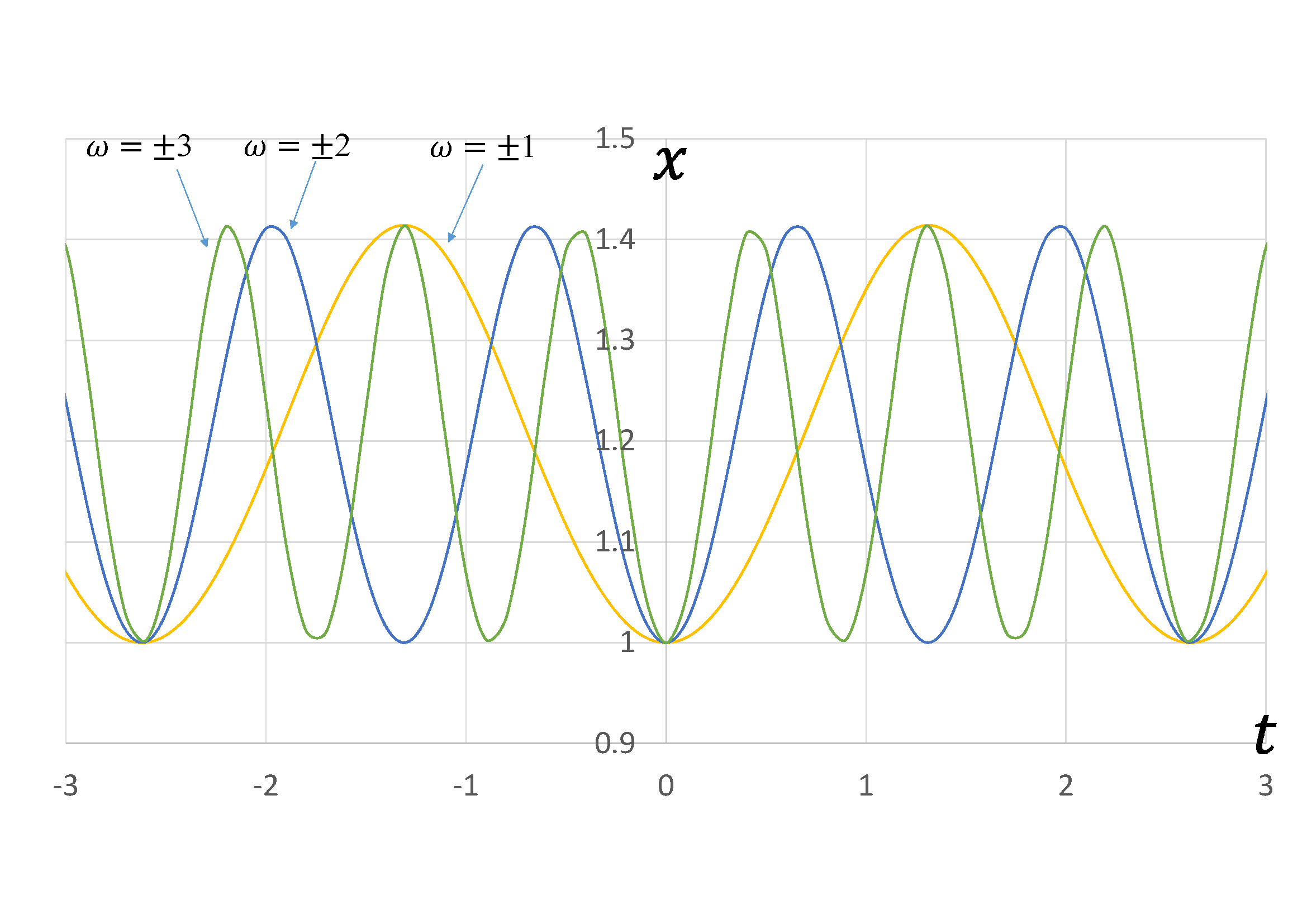}
\caption{Waves of exact solution 3 with respect to variation of the parameter $\omega$ related to the period ($A=1$)}
\label{fig85}      
\end{center}
\end{figure*}

\begin{figure*}[tb]
\begin{center}
\includegraphics[width=0.7 \textwidth]{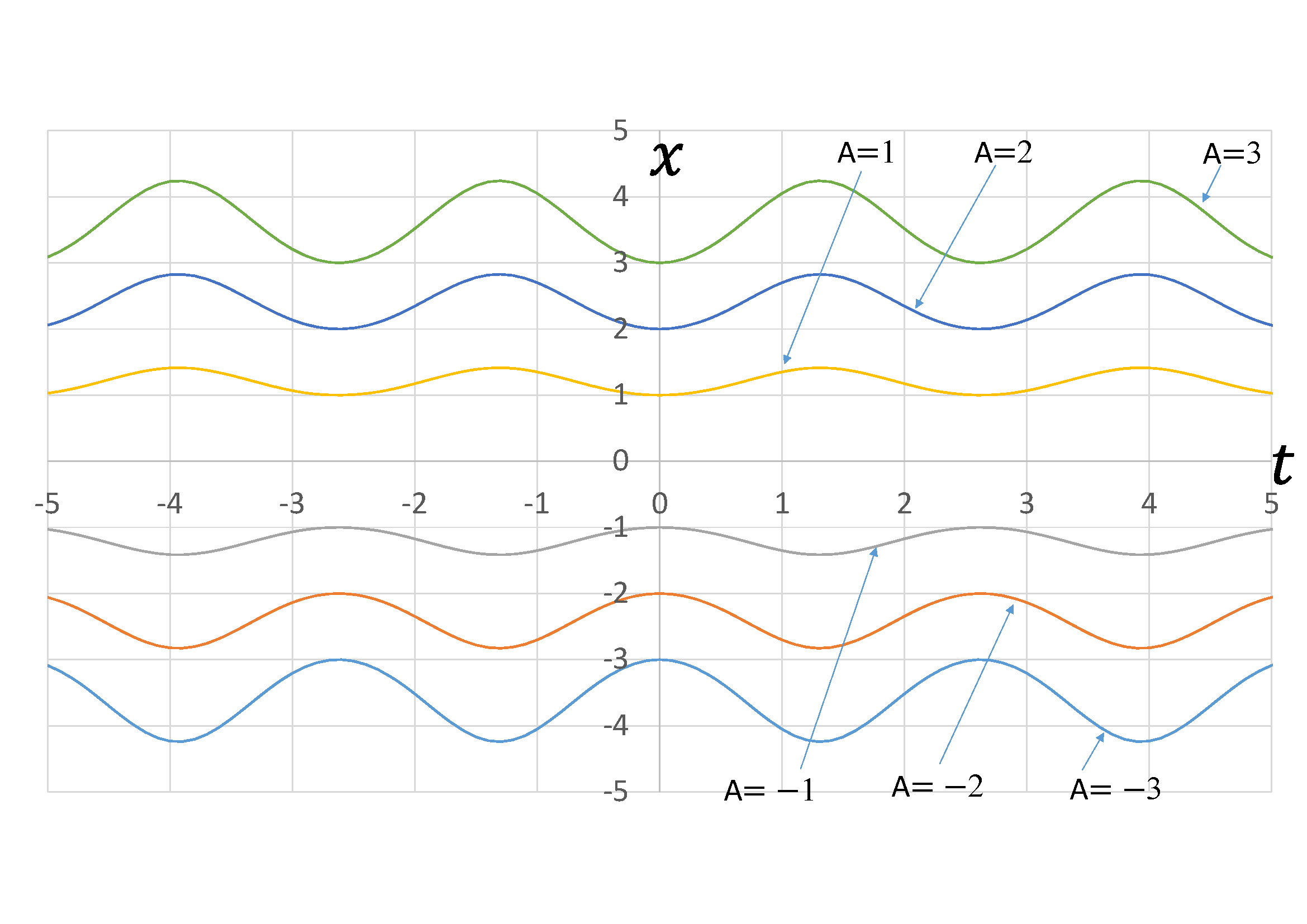}
\caption{Waves of exact solution 3 with respect to variation of the parameter $A$ related to the amplitude ($\omega$=1)}
\label{fig86}      
\end{center}
\end{figure*}

\begin{figure*}[tb]
\begin{center}
\includegraphics[width=0.7 \textwidth]{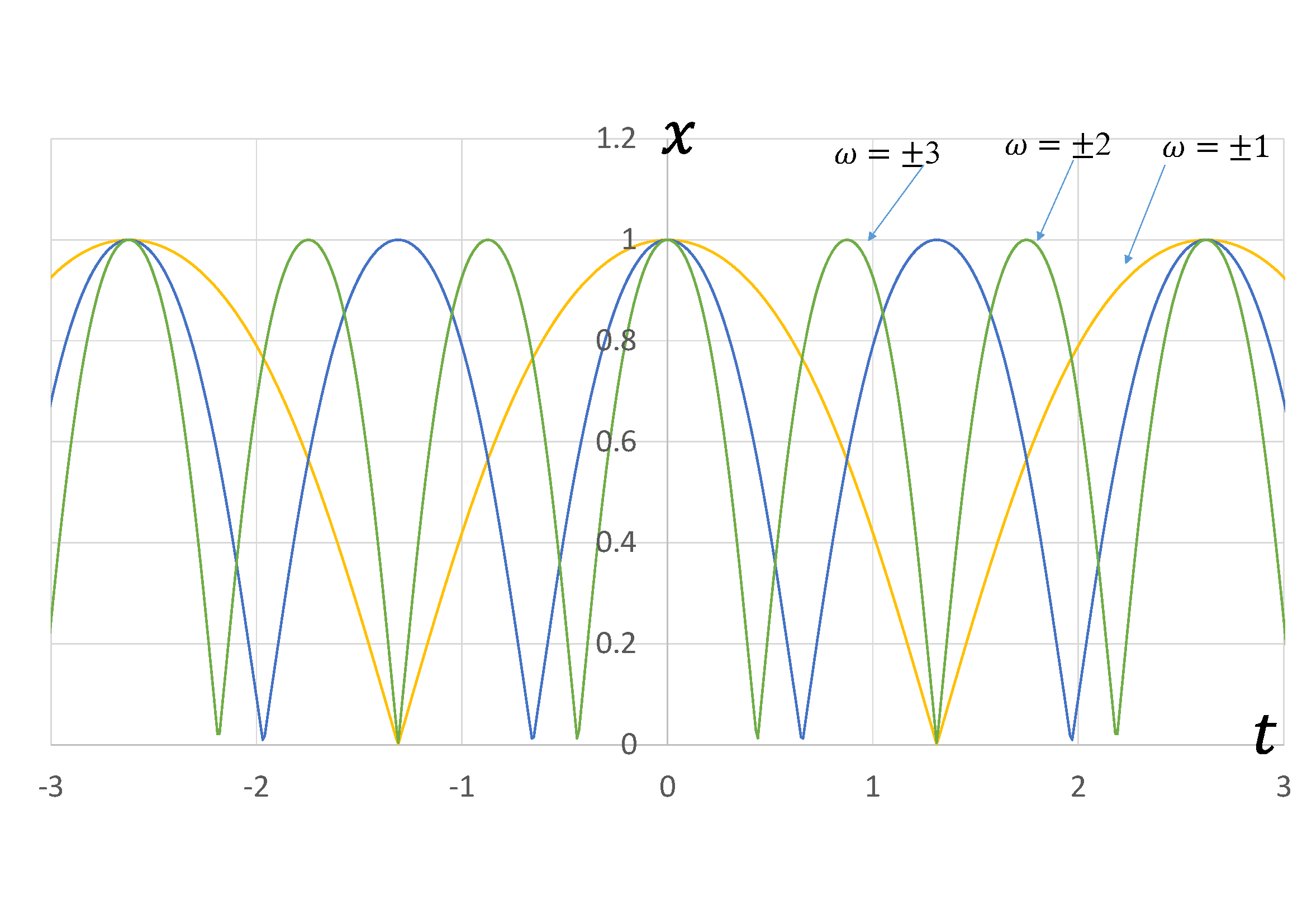}
\caption{Waves of exact solution 4 with respect to variation of the parameter $\omega$ related to the period ($A=1$)}
\label{fig87}      
\end{center}
\end{figure*}

\begin{figure*}[tb]
\begin{center}
\includegraphics[width=0.7 \textwidth]{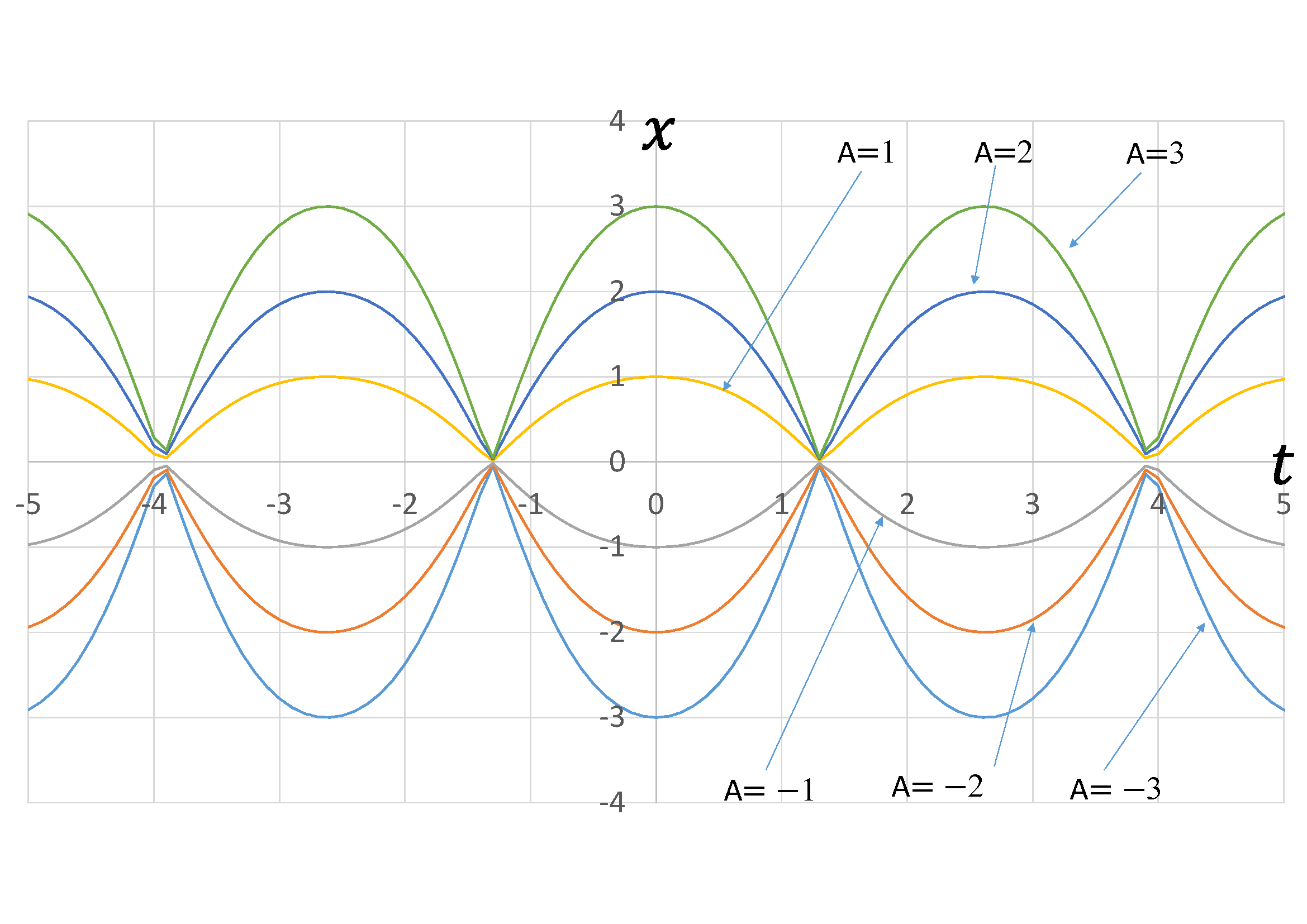}
\caption{Waves of exact solution 4 with respect to variation of the parameter $A$ related to the amplitude ($\omega=1$)}
\label{fig88}      
\end{center}
\end{figure*}

\begin{figure*}[tb]
\begin{center}
\includegraphics[width=0.7 \textwidth]{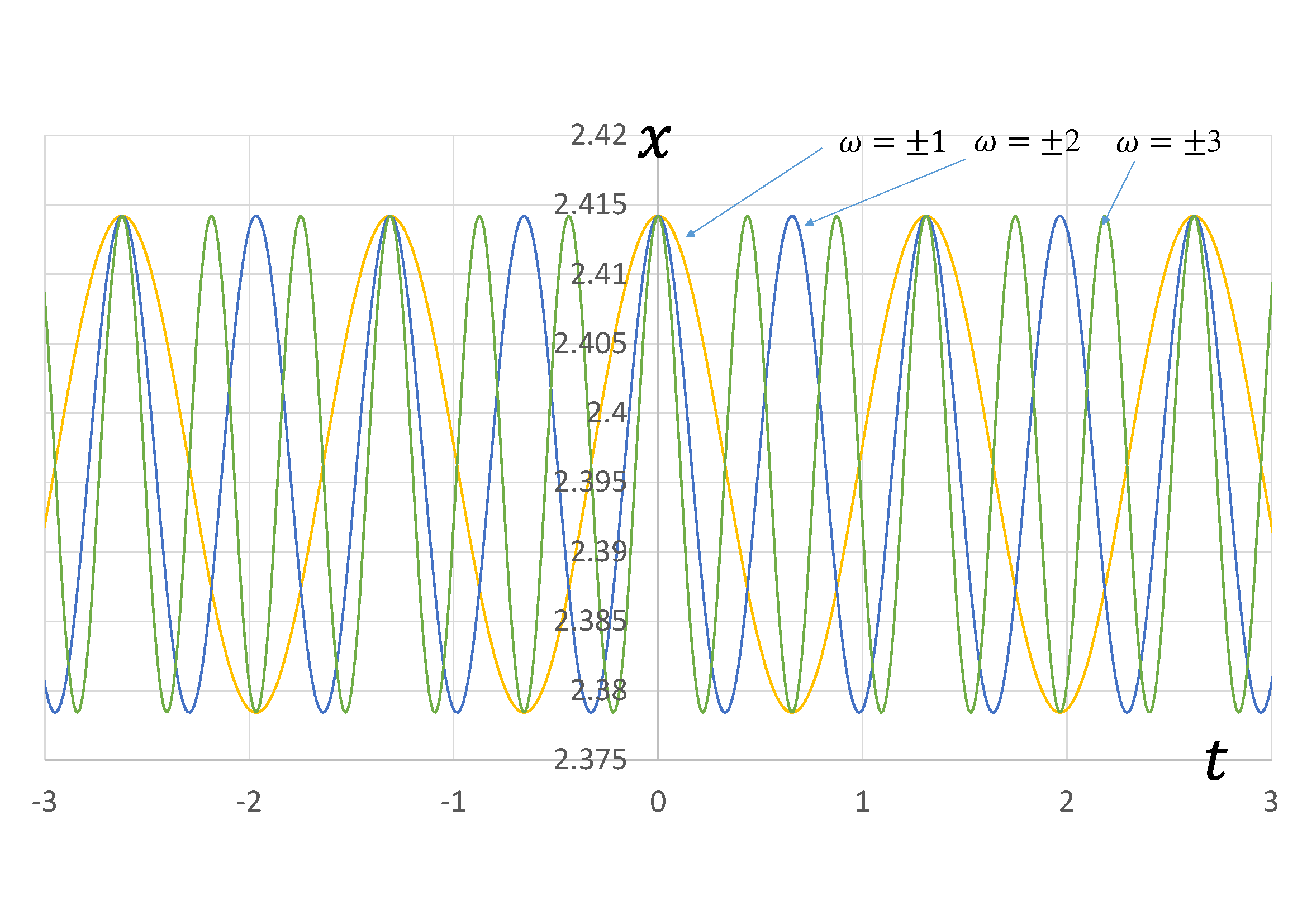}
\caption{Waves of  exact solution 5 with respect to variation of the parameter $\omega$ related to the period ($A=1$)}
\label{fig89}      
\end{center}
\end{figure*}

\begin{figure*}[tb]
\begin{center}
\includegraphics[width=0.7 \textwidth]{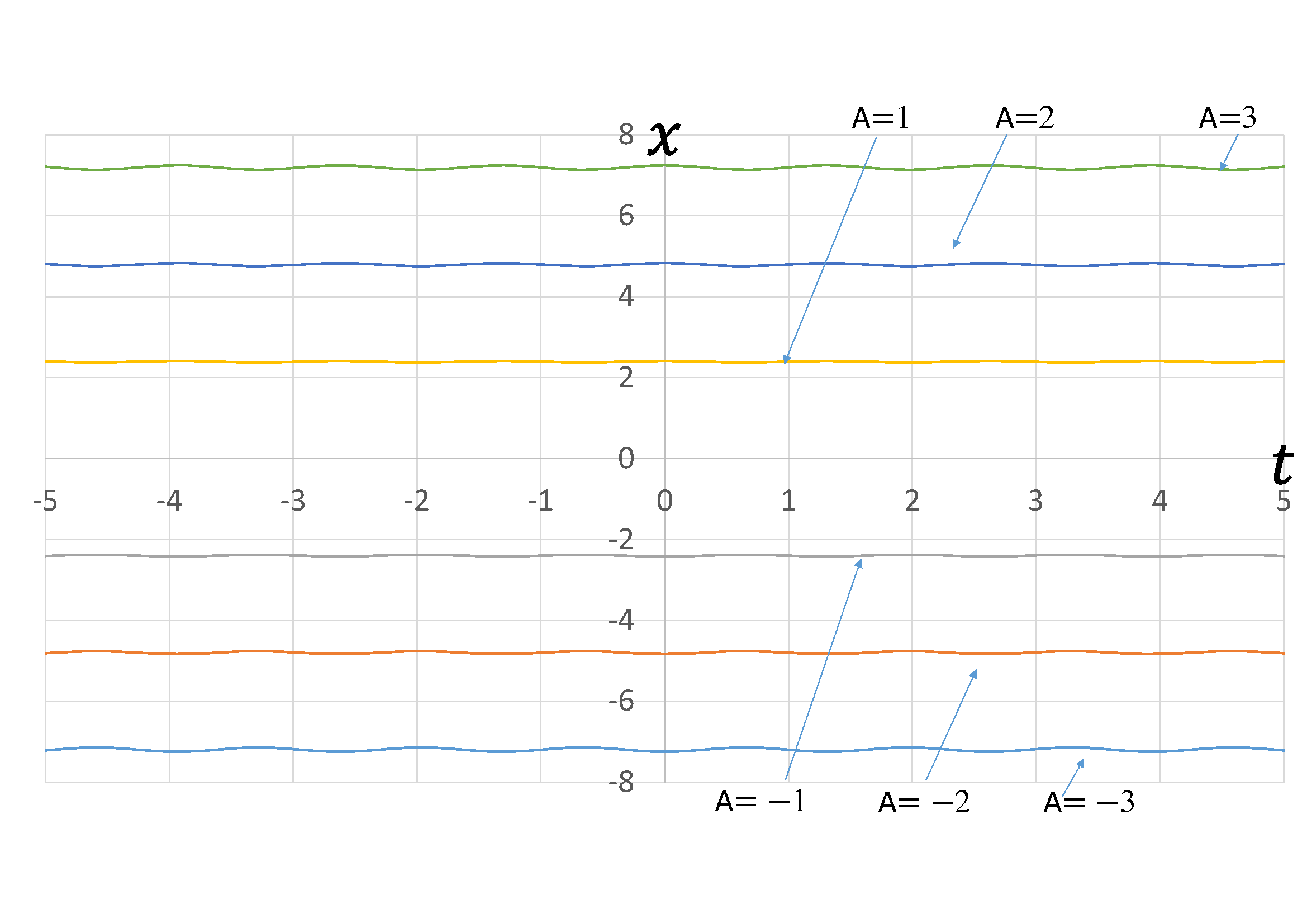}
\caption{Waves of exact solution 5 with respect to variation of the parameter $A$ related to the amplitude ($\omega=1$)}
\label{fig810}      
\end{center}
\end{figure*}

\begin{figure*}[tb]
\begin{center}
\includegraphics[width=0.7 \textwidth]{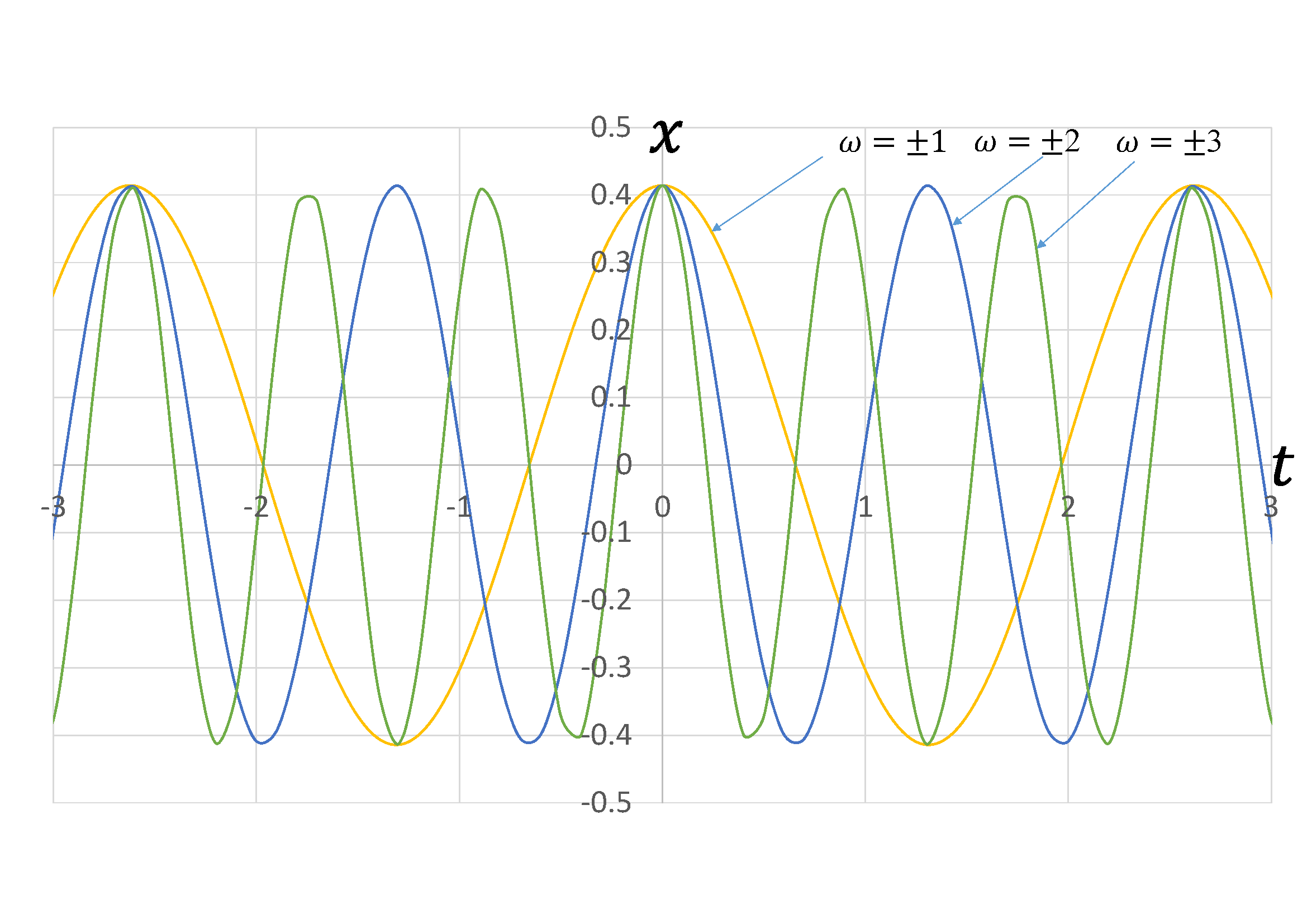}
\caption{Waves of exact solution 6 with respect to variation of the parameter $\omega$ related to the period ($A=1$)}
\label{fig811}      
\end{center}
\end{figure*}

\begin{figure*}[tb]
\begin{center}
\includegraphics[width=0.7 \textwidth]{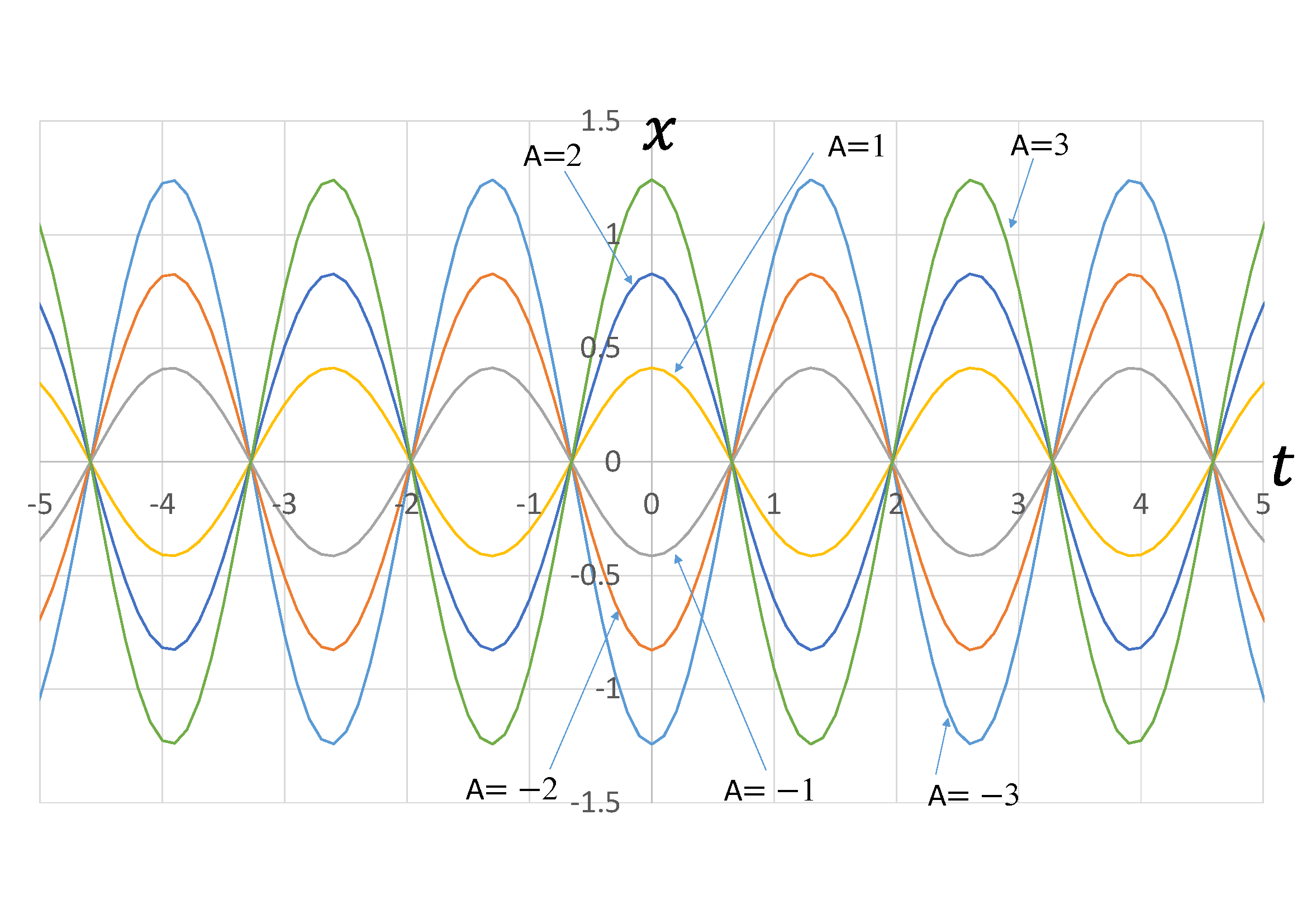}
\caption{Waves of exact solution 6 with respect to variation of the parameter $A$ related to the amplitude ($\omega=1$)}
\label{fig812}      
\end{center}
\end{figure*}

\begin{figure*}[tb]
\begin{center}
\includegraphics[width=0.7 \textwidth]{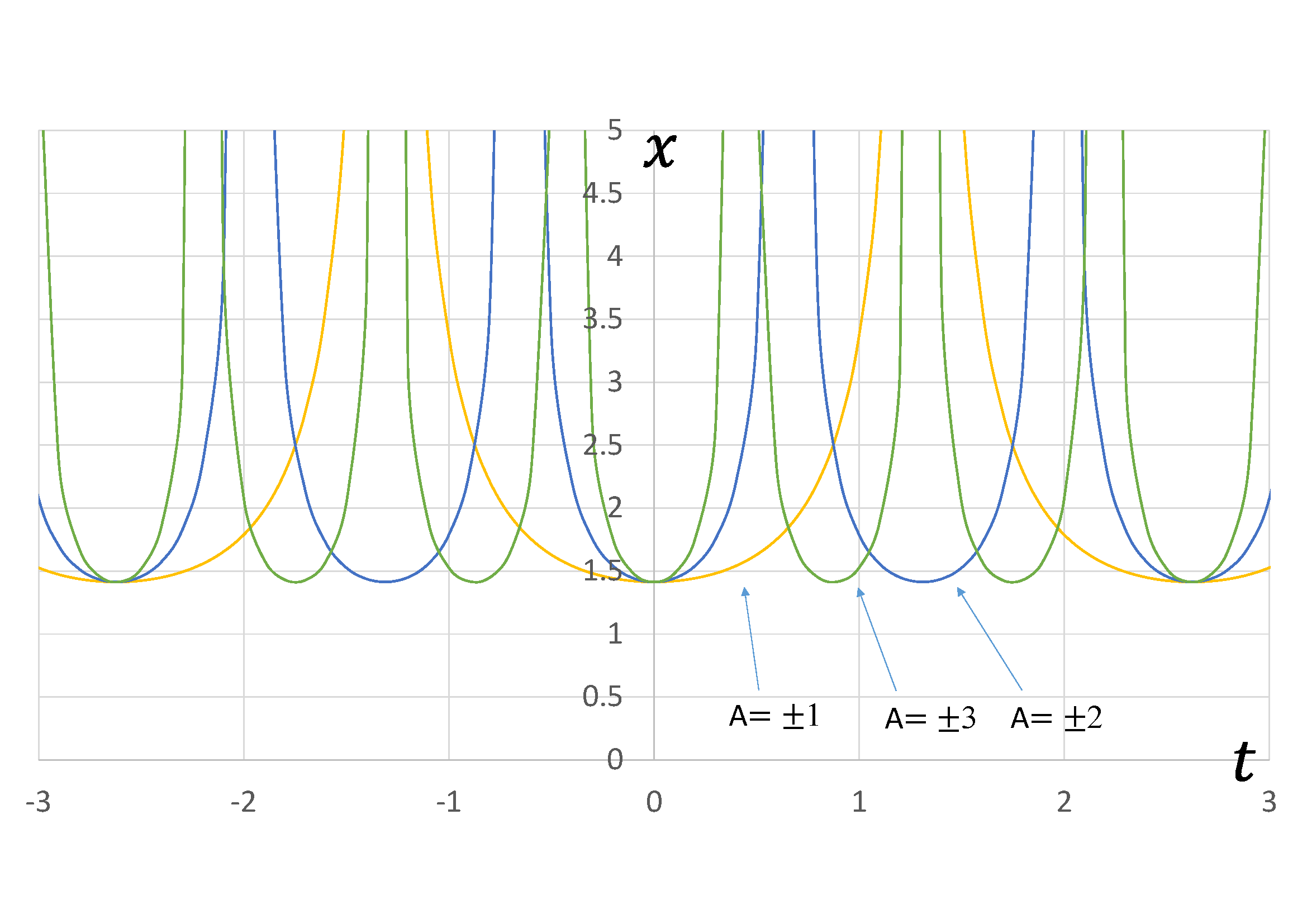}
\caption{Divergences of exact solution 7 with respect to variation of the parameter $\omega$ related to the period ($A=1$)}
\label{fig813}      
\end{center}
\end{figure*}

\begin{figure*}[tb]
\begin{center}
\includegraphics[width=0.7 \textwidth]{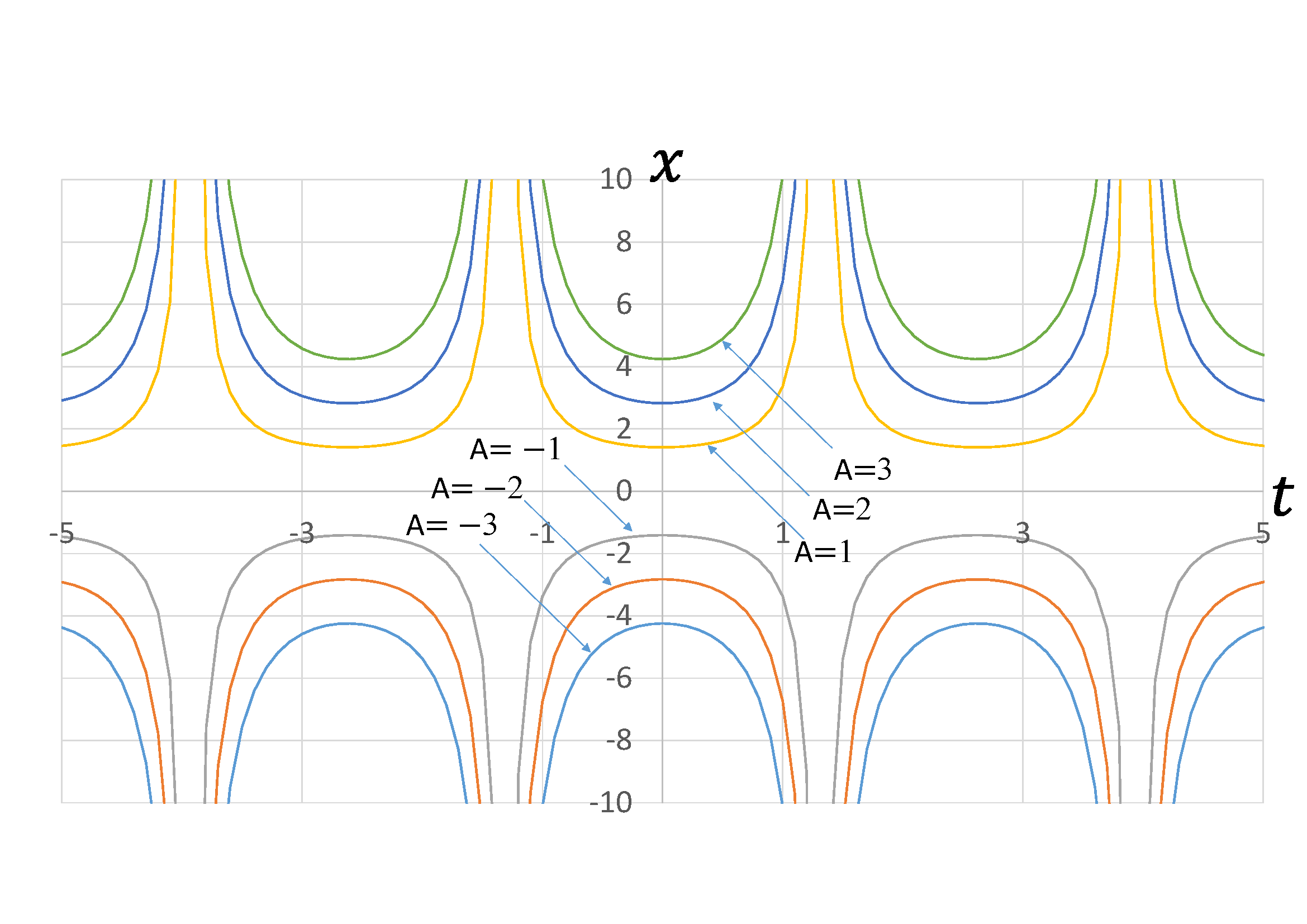}
\caption{Divergences of exact solution 7 with respect to variation of the parameter $A$ related to the amplitude ($\omega=1$)}
\label{fig814}      
\end{center}
\end{figure*}

\begin{figure*}[tb]
\begin{center}
\includegraphics[width=0.7 \textwidth]{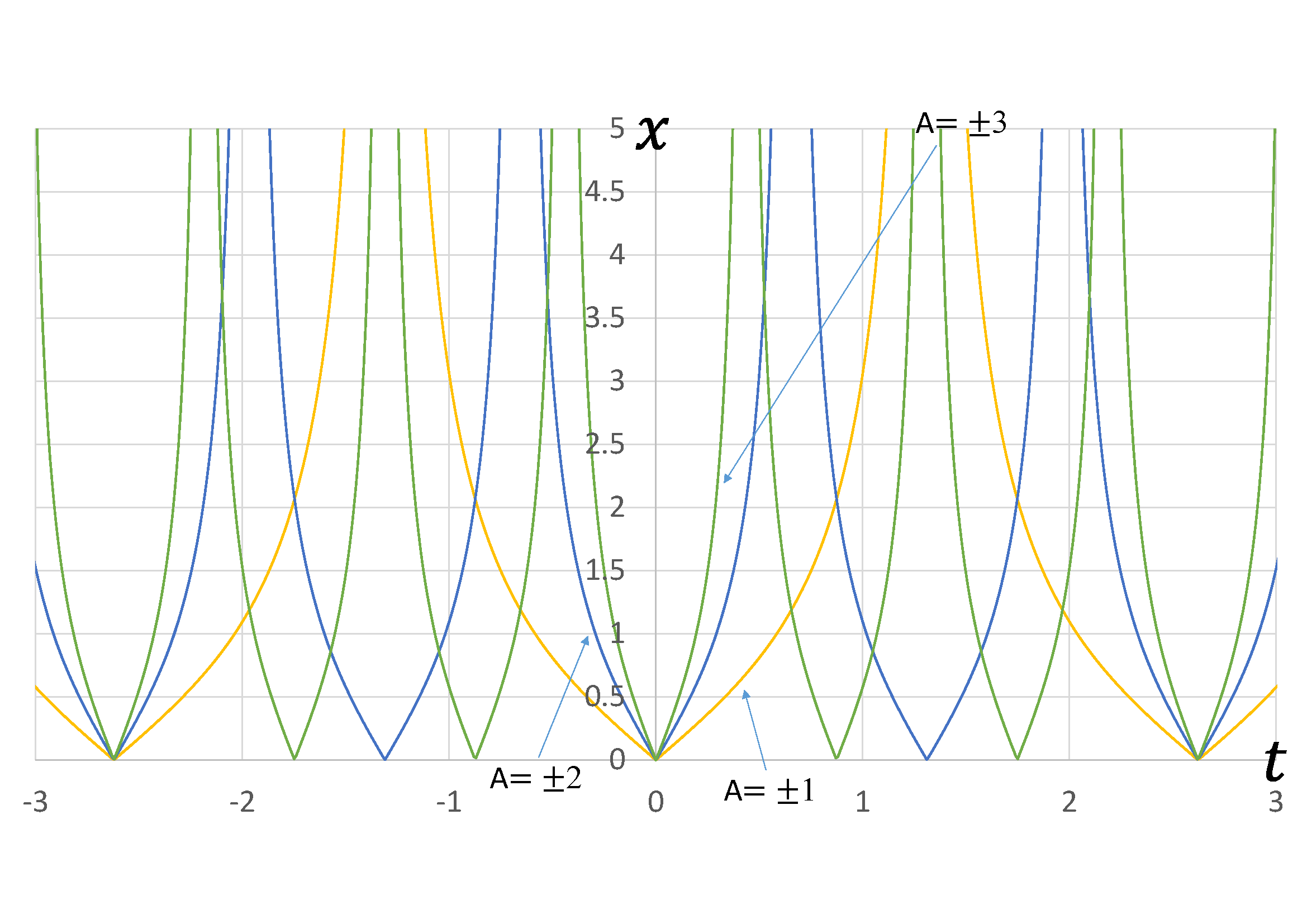}
\caption{Divergences of exact solution 8 with respect to variation of the parameter $\omega$ related to the period ($A=1$)}
\label{fig815}      
\end{center}
\end{figure*}

\begin{figure*}[tb]
\begin{center}
\includegraphics[width=0.7 \textwidth]{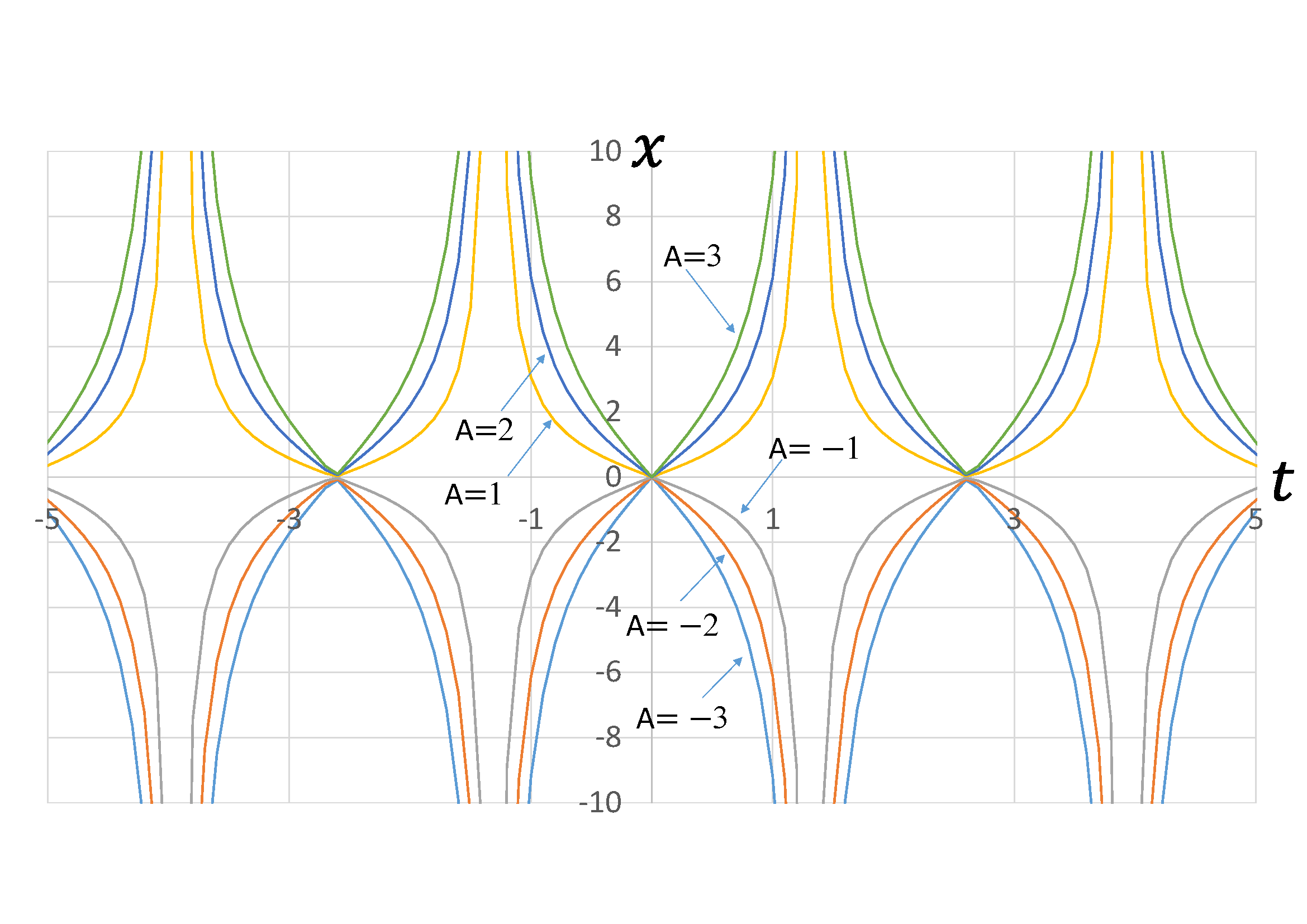}
\caption{Divergences of exact solution 8 with respect to variation of the parameter $A$ related to the amplitude ($\omega=1$)}
\label{fig816}      
\end{center}
\end{figure*}

\begin{figure*}[tb]
\begin{center}
\includegraphics[width=0.7 \textwidth]{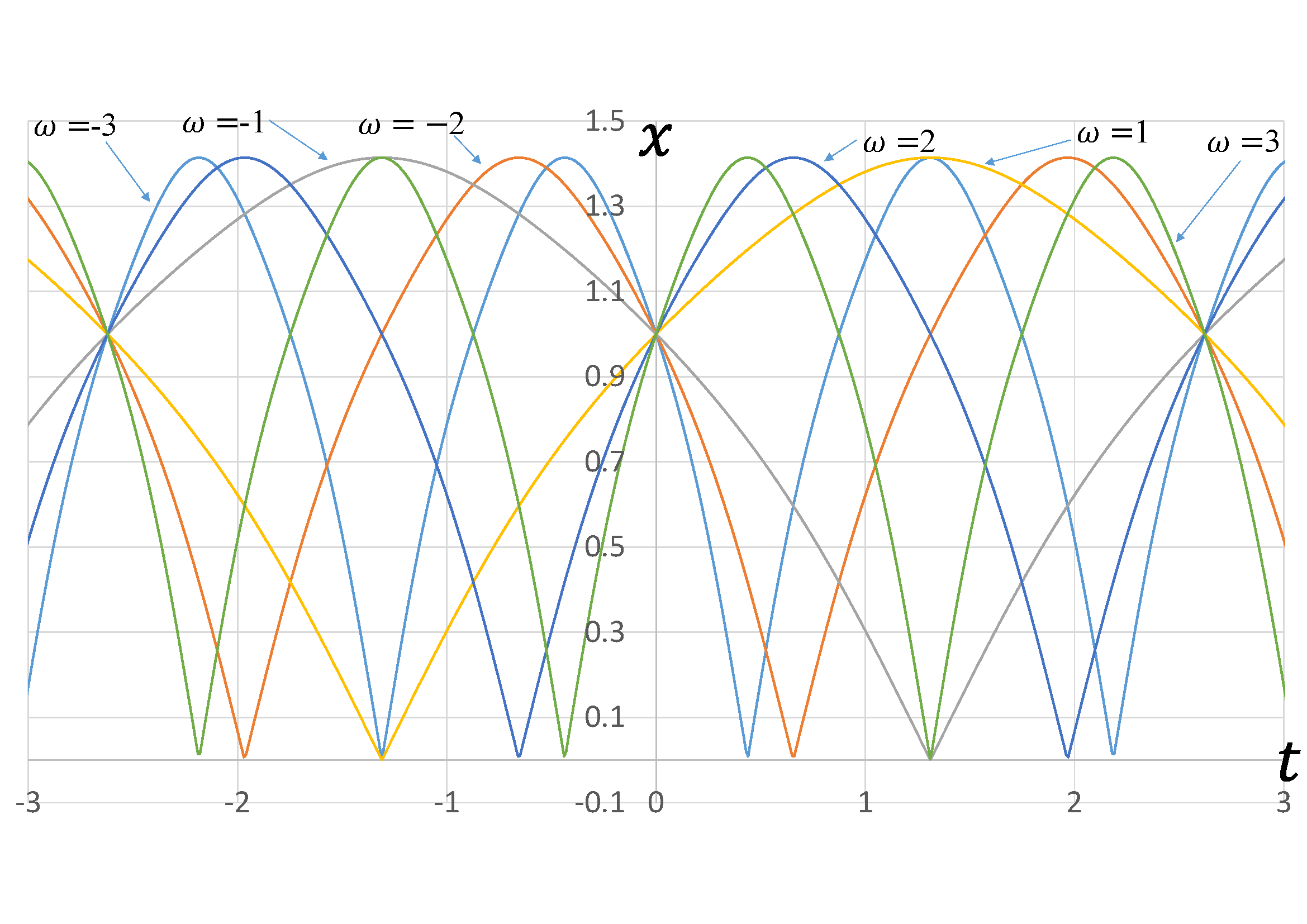}
\caption{Waves of exact solution 9 with respect to variation of the parameter $\omega$ related to the period ($A=1$)}
\label{fig817}      
\end{center}
\end{figure*}

\begin{figure*}[tb]
\begin{center}
\includegraphics[width=0.7 \textwidth]{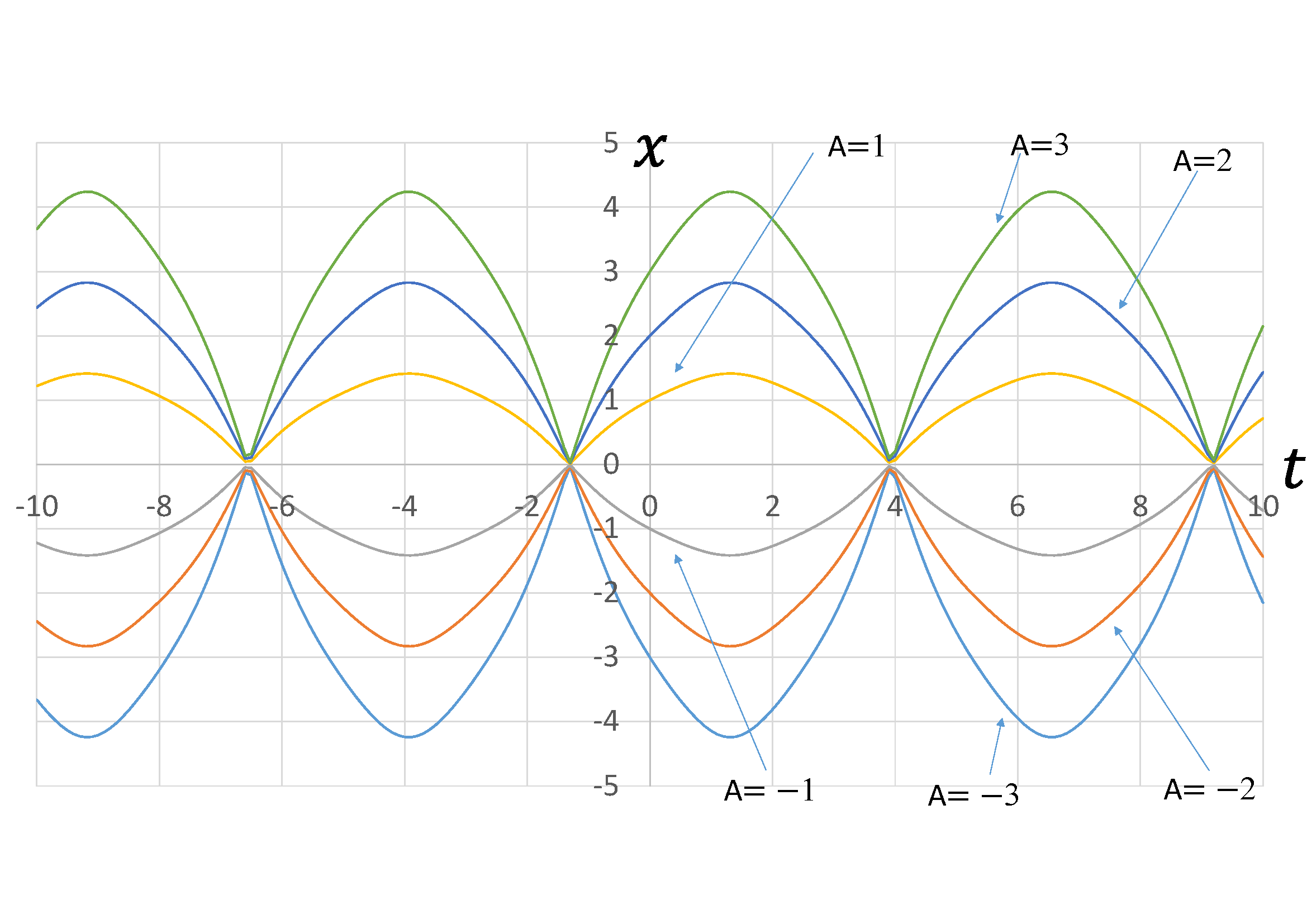}
\caption{Waves of exact solution 9 with respect to variation of the parameter $A$ related to the amplitude ($\omega=1$)}
\label{fig818}      
\end{center}
\end{figure*}

\begin{figure*}[tb]
\begin{center}
\includegraphics[width=0.7 \textwidth]{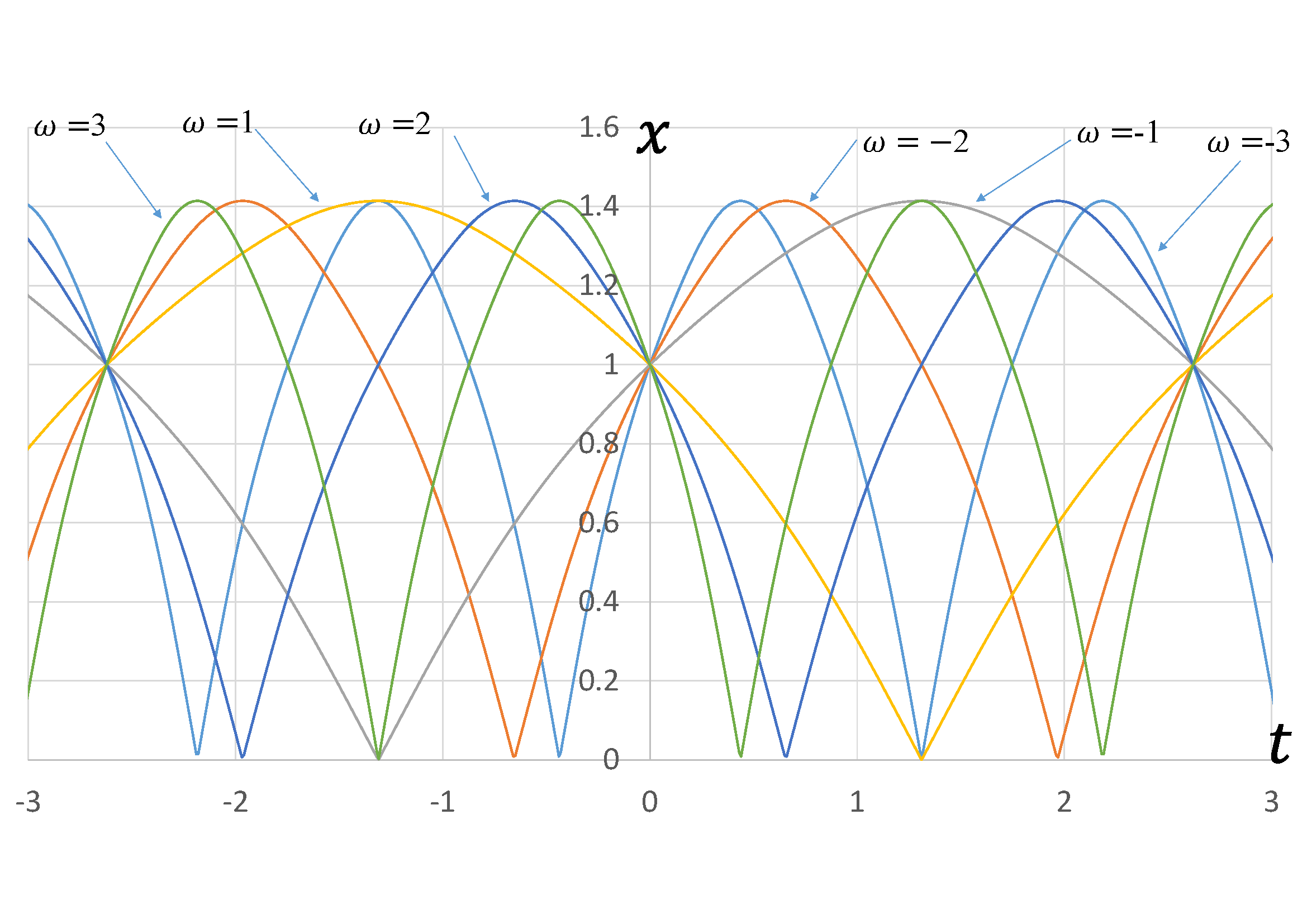}
\caption{Waves of exact solution 10 with respect to variation of the parameter $\omega$ related to the period ($A=1$)}
\label{fig819}      
\end{center}
\end{figure*}

\begin{figure*}[tb]
\begin{center}
\includegraphics[width=0.7 \textwidth]{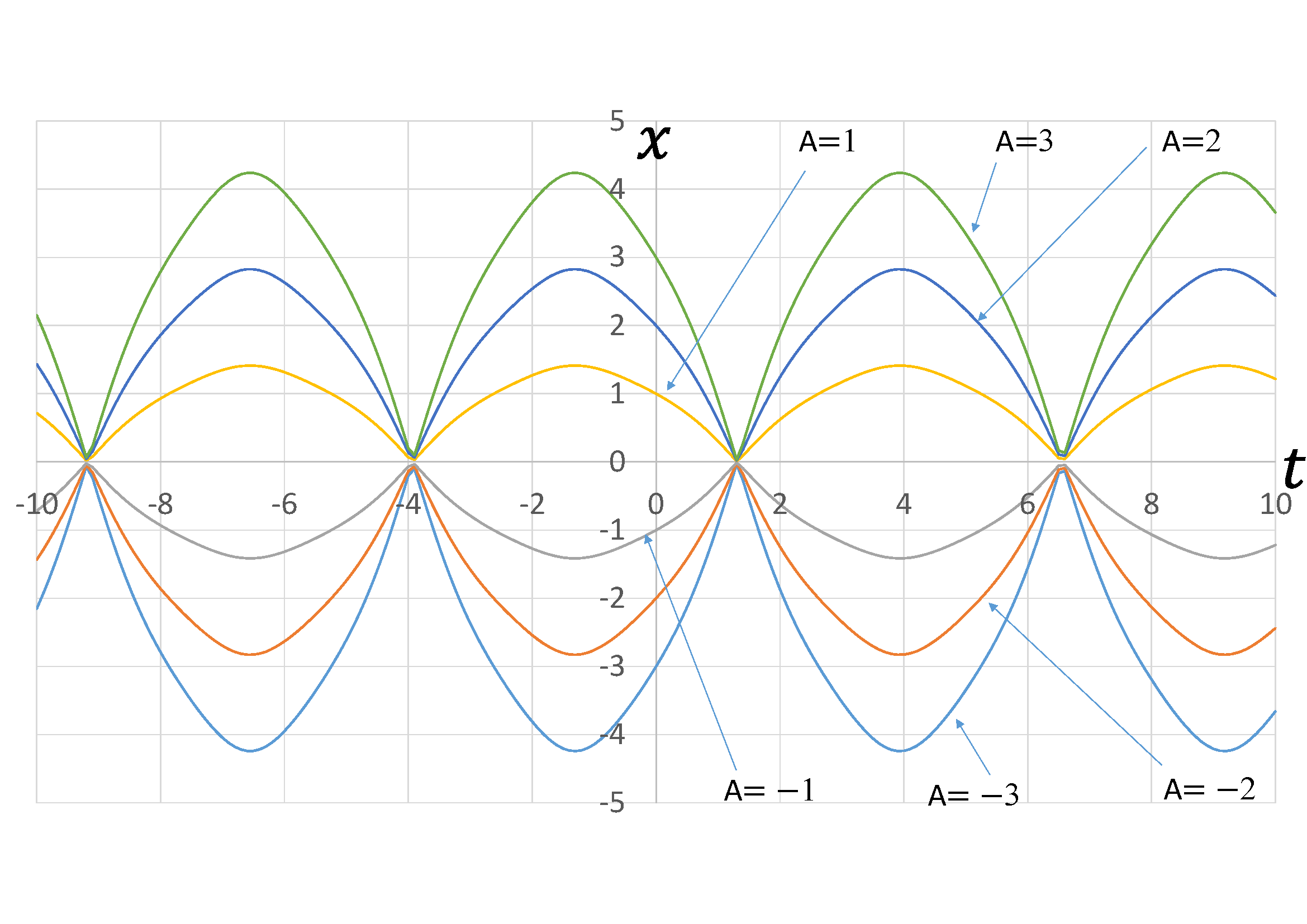}
\caption{Waves of exact solution 10 with respect to variation of the parameter $A$ related to the amplitude ($\omega=1$)}
\label{fig820}      
\end{center}
\end{figure*}

\begin{figure*}[tb]
\begin{center}
\includegraphics[width=0.7 \textwidth]{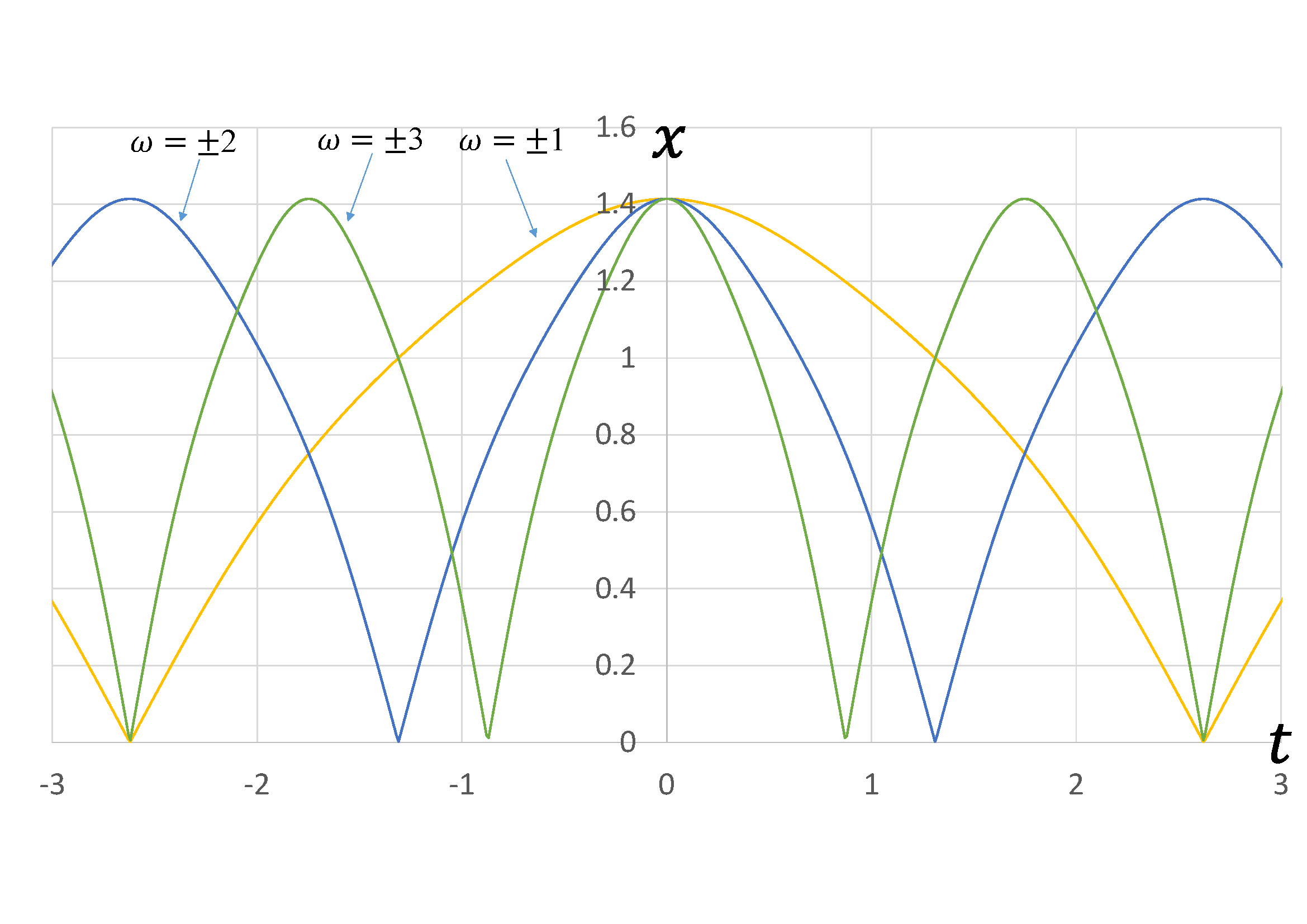}
\caption{Waves of exact solution 11 with respect to variation of the parameter $\omega$ related to the period ($A=1$)}
\label{fig821}      
\end{center}
\end{figure*}

\begin{figure*}[tb]
\begin{center}
\includegraphics[width=0.7 \textwidth]{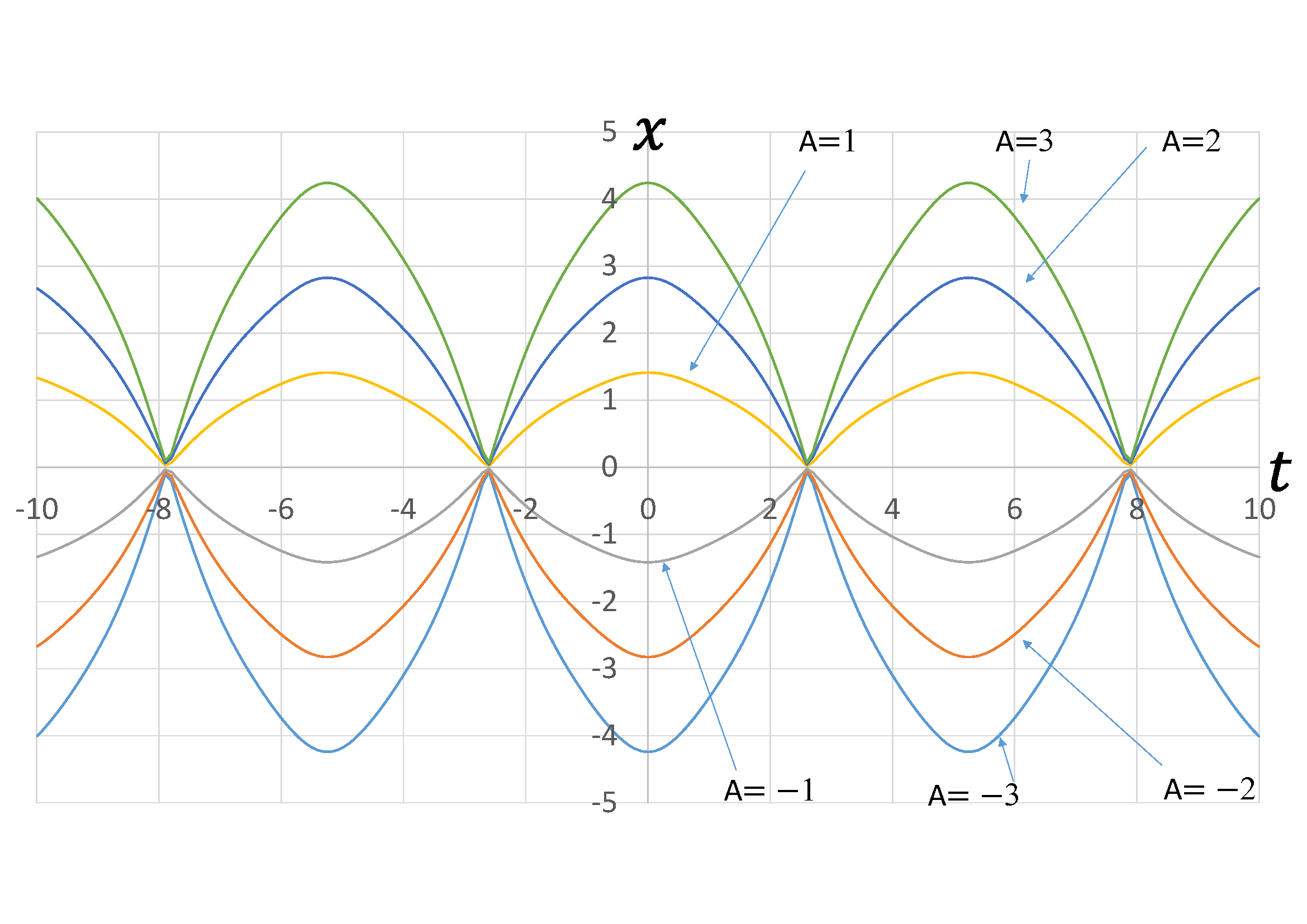}
\caption{Waves of exact solution 11 with respect to variation of the parameter $A$ related to the amplitude ($\omega=1$)}
\label{fig822}      
\end{center}
\end{figure*}

\begin{figure*}[tb]
\begin{center}
\includegraphics[width=0.7 \textwidth]{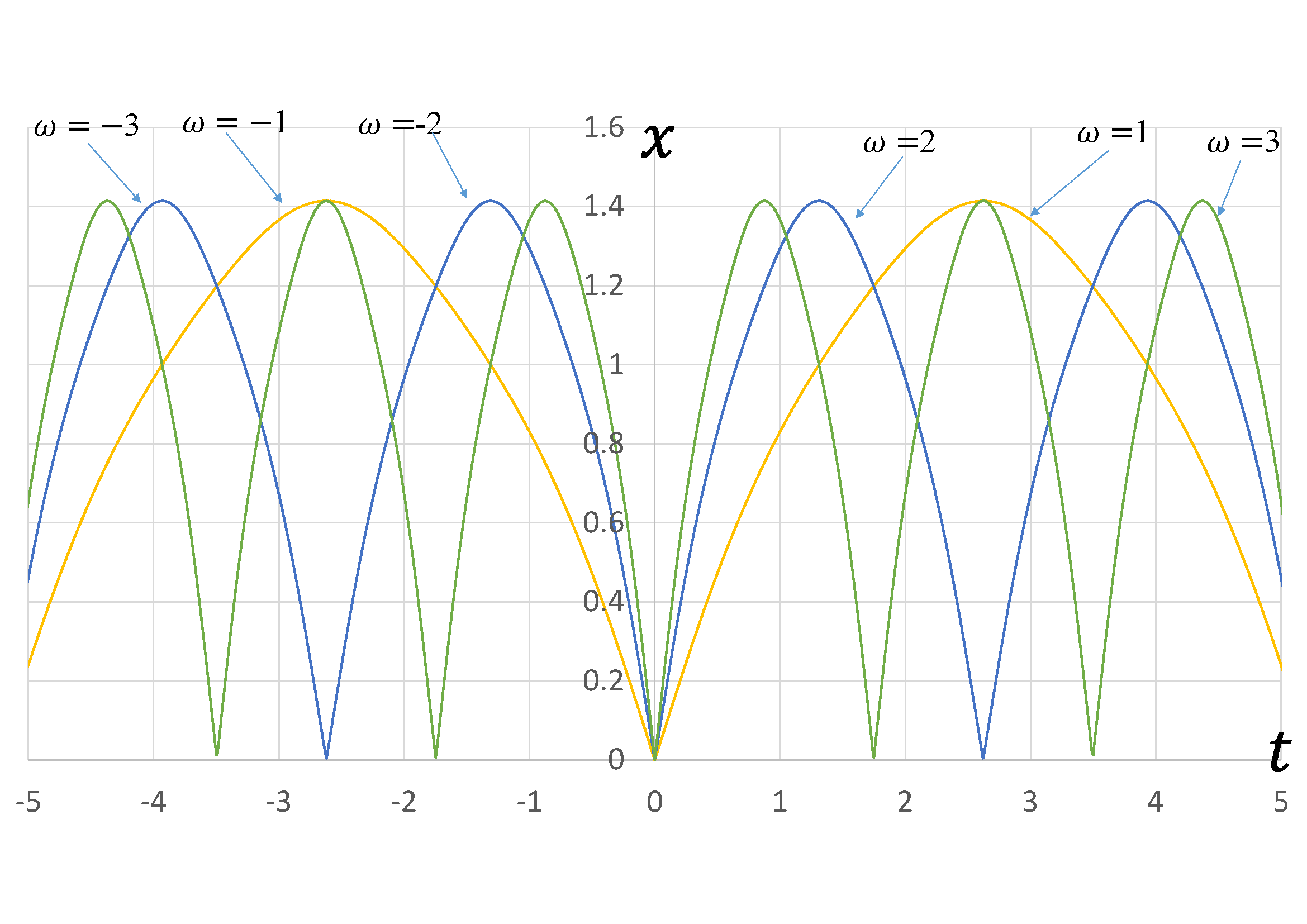}
\caption{Waves of exact solution 12 with respect to variation of the parameter $\omega$ related to the period ($A=1$)}
\label{fig823}      
\end{center}
\end{figure*}

\begin{figure*}[tb]
\begin{center}
\includegraphics[width=0.7 \textwidth]{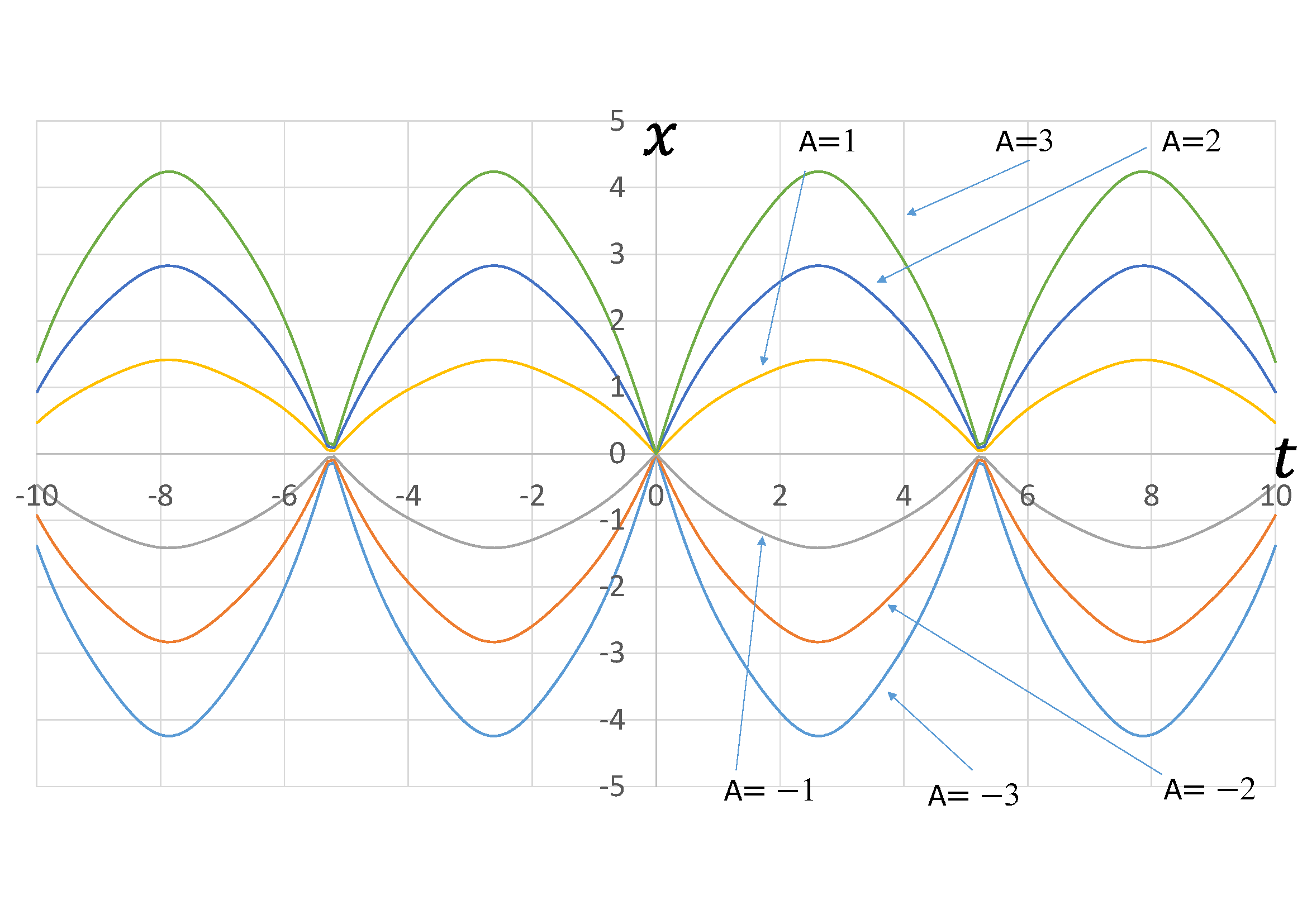}
\caption{Waves of exact solution 12 with respect to variation of the parameter $A$ related to the amplitude ($\omega=1$)}
\label{fig824}      
\end{center}
\end{figure*}

\begin{figure*}[tb]
\begin{center}
\includegraphics[width=0.7 \textwidth]{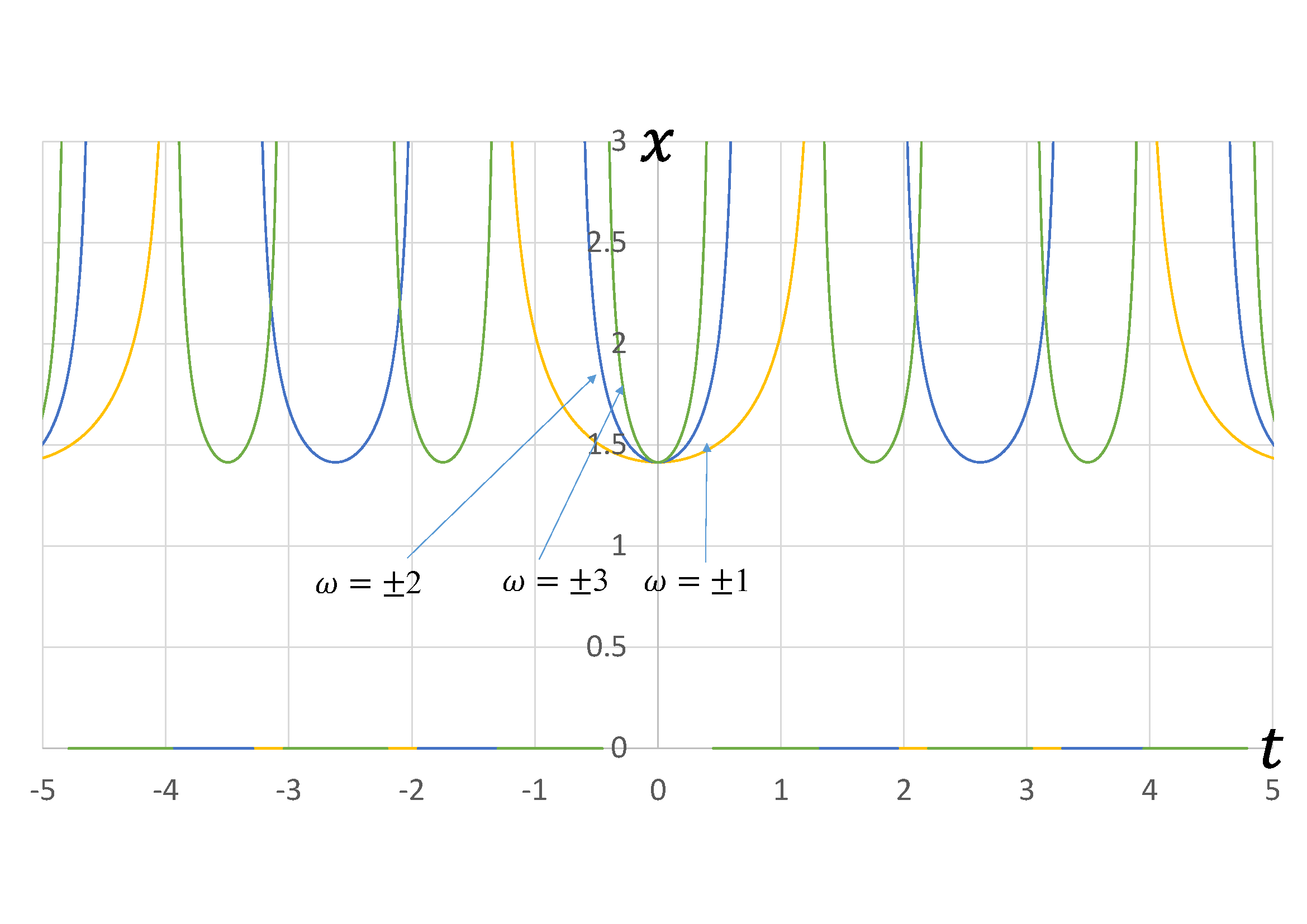}
\caption{Divergences of exact solution 13 with respect to variation of the parameter $\omega$ related to the period ($A=1$)}
\label{fig825}      
\end{center}
\end{figure*}

\begin{figure*}[tb]
\begin{center}
\includegraphics[width=0.7 \textwidth]{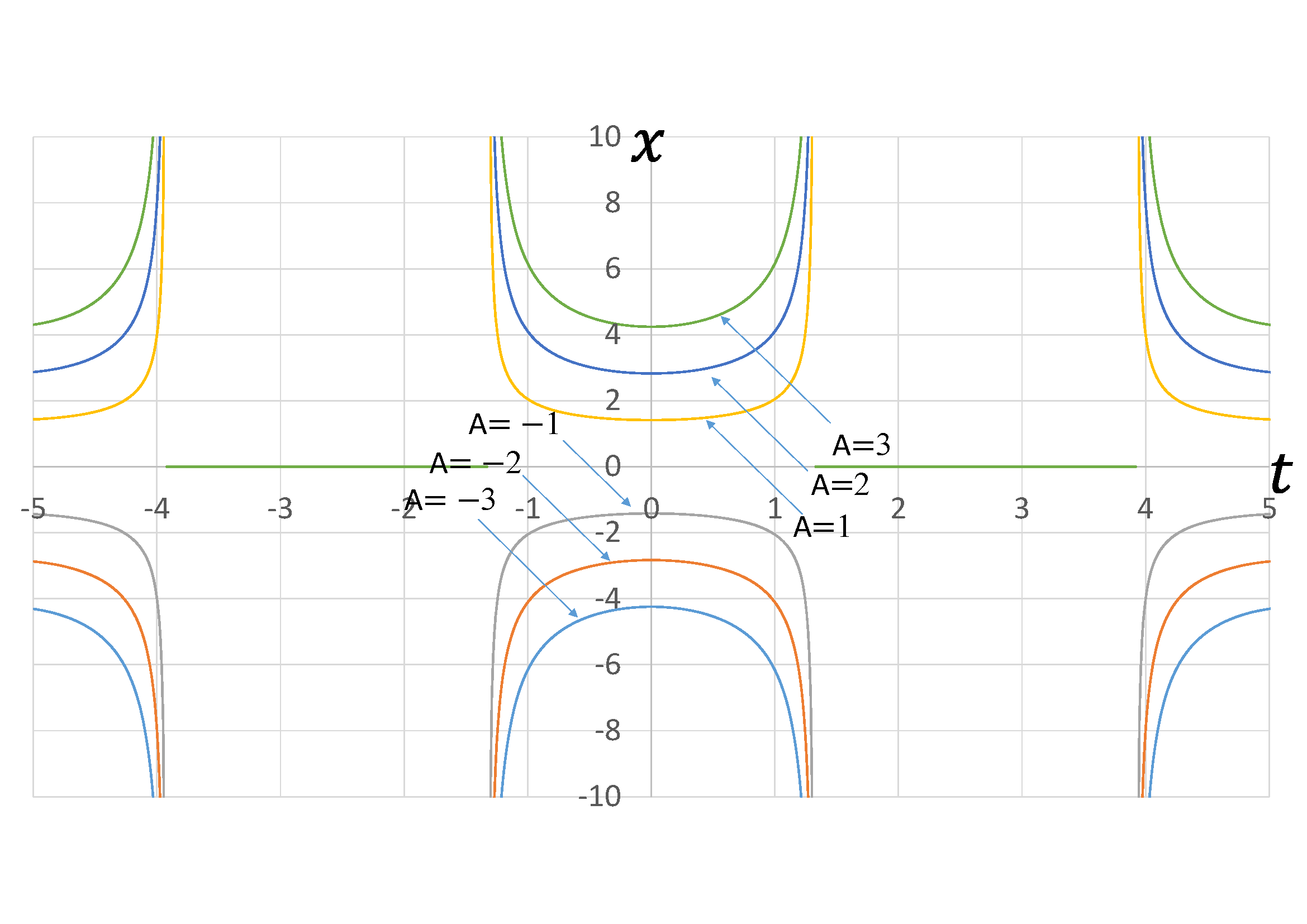}
\caption{Divergences of exact solution 13 with respect to variation of the parameter $A$ related to the amplitude ($\omega=1$)}
\label{fig826}      
\end{center}
\end{figure*}

\begin{figure*}[tb]
\begin{center}
\includegraphics[width=0.7 \textwidth]{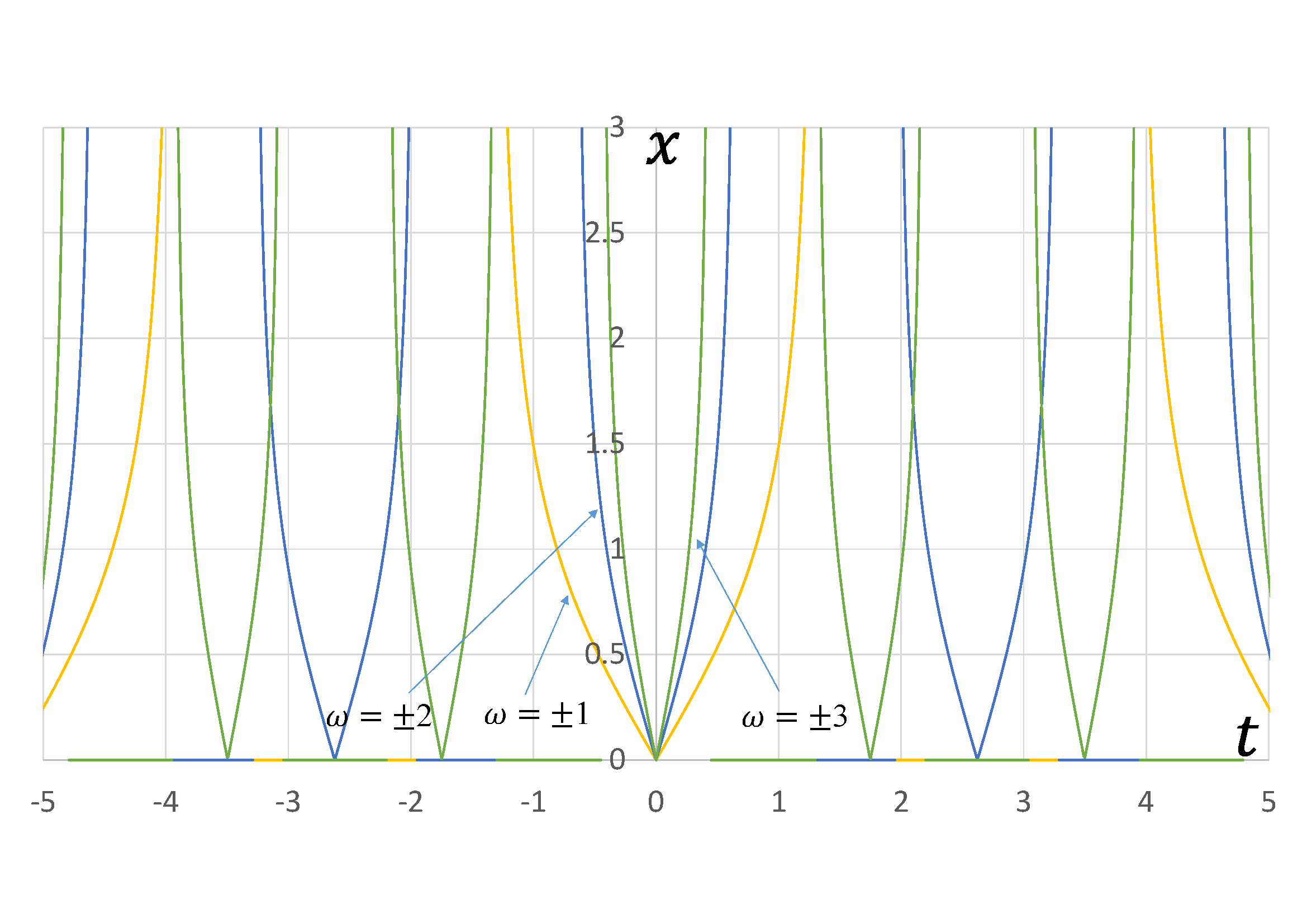}
\caption{Divergences of exact solution 14 with respect to variation of the parameter $\omega$ related to the period ($A=1$)}
\label{fig827}      
\end{center}
\end{figure*}

\begin{figure*}[tb]
\begin{center}
\includegraphics[width=0.7 \textwidth]{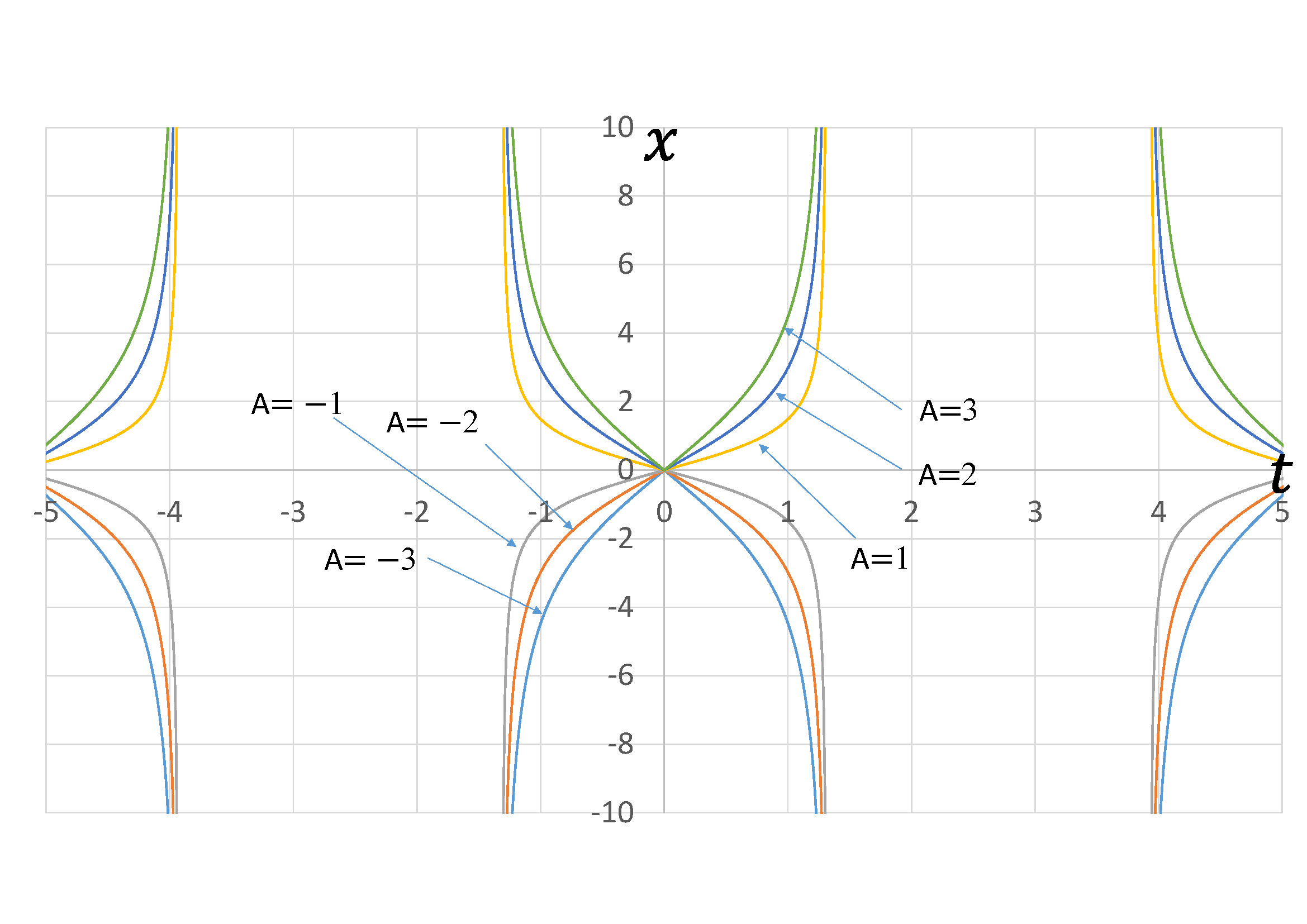}
\caption{Divergences of exact solution 14 with respect to variation of the parameter $A$ related to the amplitude($\omega=1$)}
\label{fig828}      
\end{center}
\end{figure*}

Figs. \ref{fig82}, \ref{fig84}, \ref{fig86}, \ref{fig88}, \ref{fig810}, \ref{fig812}, \ref{fig814}, \ref{fig816}, \ref{fig818}, \ref{fig820},\ref{fig822}, \ref{fig824}, \ref{fig826}, and \ref{fig828} show the variations of the amplitude parameter $A$ with $\omega=1$. Because the parameter $\omega$ is fixed at $\omega=1$, the wave periods are $T=\pi_2$ (Figs. \ref{fig82}, \ref{fig84}, \ref{fig86}, \ref{fig88}, and \ref{fig812}), $T=\frac{\pi_2}{2}$ (Fig. \ref{fig810}), $T=2 \eta_2$ (Figs. \ref{fig814}, \ref{fig816}, \ref{fig826}, and \ref{fig828}), and $T=2 \pi_2$ (Figs. \ref{fig818}, \ref{fig820}, \ref{fig822}, and \ref{fig824}). As the absolute value $|A|$ increases, the wave  amplitude increases with the period kept constant.  As the absolute value $|A|$ decreases, the wave amplitude decreases with the period kept constant. The range of variable $x(t)$ satisfies the following inequalities: 

(i)  For Figs. \ref{fig82} and \ref{fig86} (exact solutions 1 and 3),
\begin{equation}
|A| \leqq  x(t) \leqq \sqrt{2} |A|.
\label{8.13}
\end{equation}

(ii)  For Figs. \ref{fig84} and \ref{fig88} (exact solutions 2 and 4)s
\begin{equation}
0 \leqq  x(t) \leqq |A|.
\label{8.14}
\end{equation}

(iii)  For Fig. \ref{fig810} (exact solution 5),
\begin{equation}
2^{\frac{5}{4}} |A| \leqq  x(t) \leqq (1+\sqrt{2}) |A|.
\label{8.15}
\end{equation}

(iv)  For Fig. \ref{fig812} (exact solution 6),
\begin{equation}
-(\sqrt{2}-1) |A|\leqq  x(t) \leqq (\sqrt{2}-1) |A|.
\label{8.16}
\end{equation}

(v)  For Figs. \ref{fig814} and \ref{fig826} (exact solution 7 and 13),
\begin{equation}
\sqrt{2} A \leqq  x(t)  \ (A \geqq 0),
\label{8.17}
\end{equation}
\begin{equation}
 x(t) \leqq \sqrt{2} A \ (A \leqq 0).
\label{8.18}
\end{equation}

(vi)  For Figs. \ref{fig818}, \ref{fig820}, \ref{fig822}, and \ref{fig824} (exact solutions 9, 10, 11, and 12),
\begin{equation}
0 \leqq  x(t)  \leqq \sqrt{2}|A|.
\label{8.19}
\end{equation}

(vii)   For Figs. \ref{fig816} and \ref{fig828} (exact solutions 8 and 14),
\begin{equation}
0 \leqq  x(t)  \ (A \geqq 0),
\label{8.20}
\end{equation}
\begin{equation}
 x(t) \leqq 0 \ (A \leqq 0).
\label{8.21}
\end{equation}

\begin{figure*}[tb]
\begin{center}
\includegraphics[width=0.70 \textwidth]{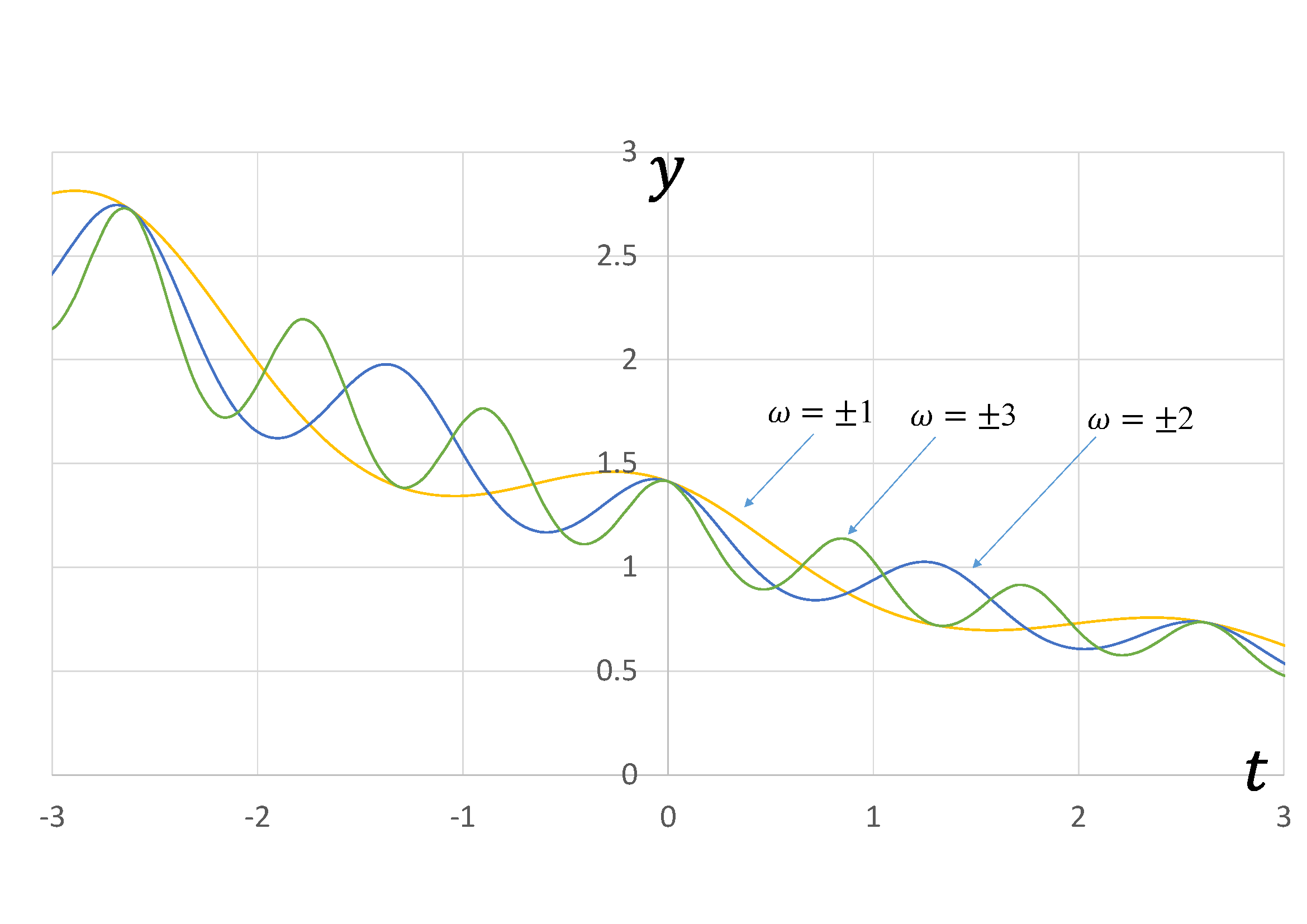}
\caption{Extended exact solution 1with respect to variation of the parameter $\omega$ related to the period ($A=1$)}
\label{fig829}      
\end{center}
\end{figure*}

\begin{figure*}[tb]
\begin{center}
\includegraphics[width=0.7 \textwidth]{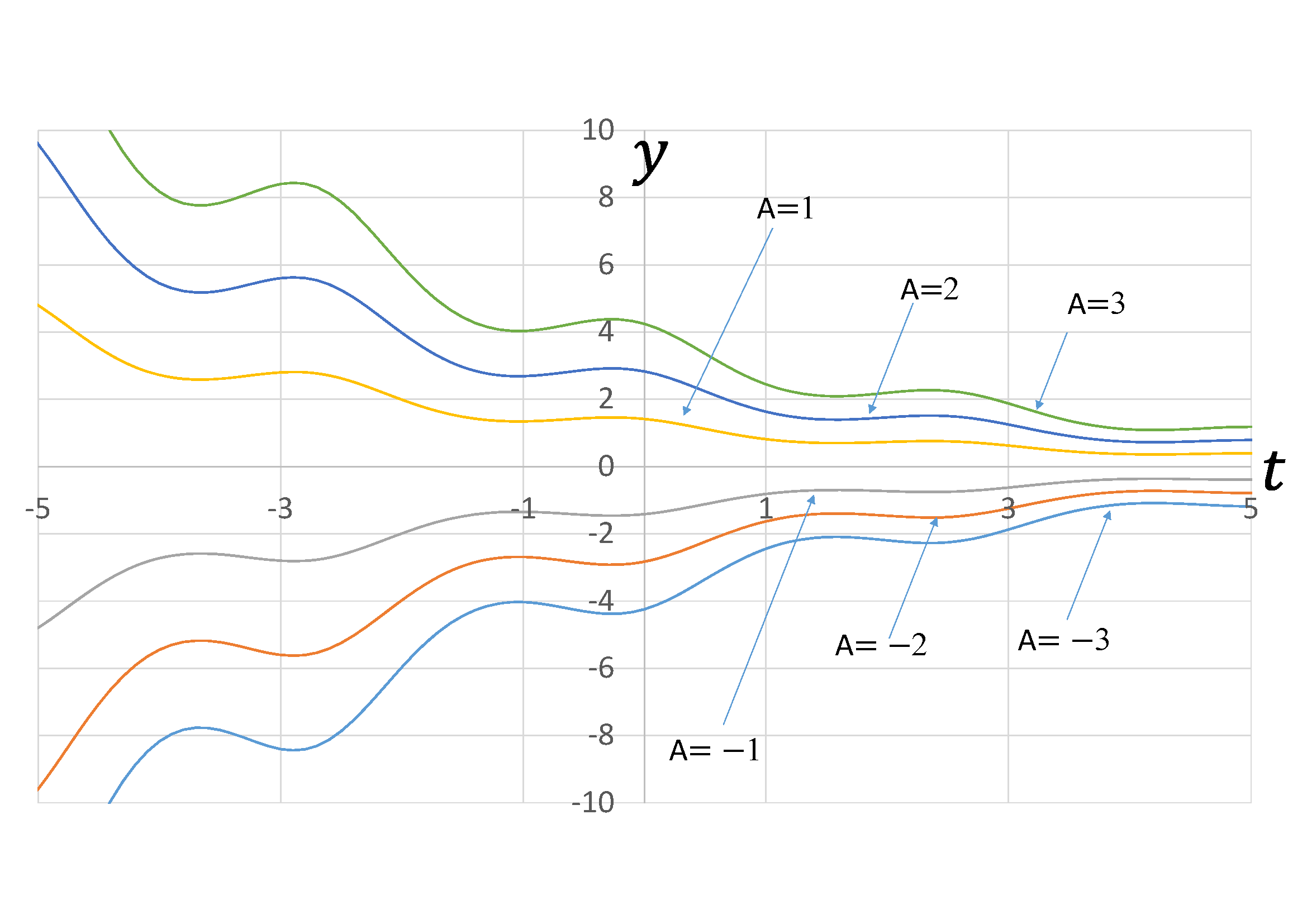}
\caption{Extended exact solution 1 with respect to variation of the parameter $A$ related to the amplitude($\omega=1$)}
\label{fig830}      
\end{center}
\end{figure*}

\begin{figure*}[tb]
\begin{center}
\includegraphics[width=0.70 \textwidth]{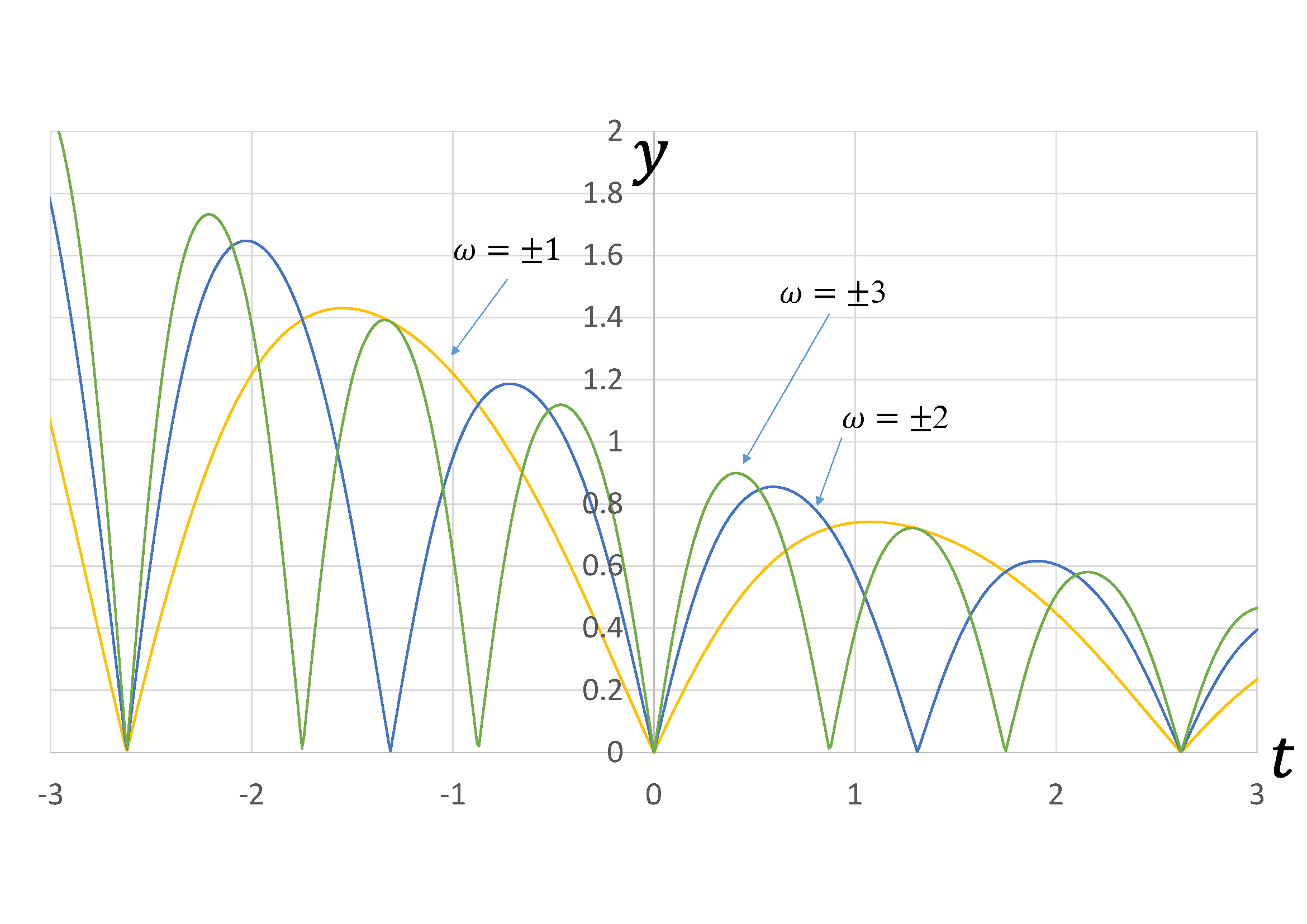}
\caption{Extended exact solution 2 with respect to variation of the parameter $\omega$ related to the period ($A=1$)}
\label{fig831}      
\end{center}
\end{figure*}

\begin{figure*}[tb]
\begin{center}
\includegraphics[width=0.70 \textwidth]{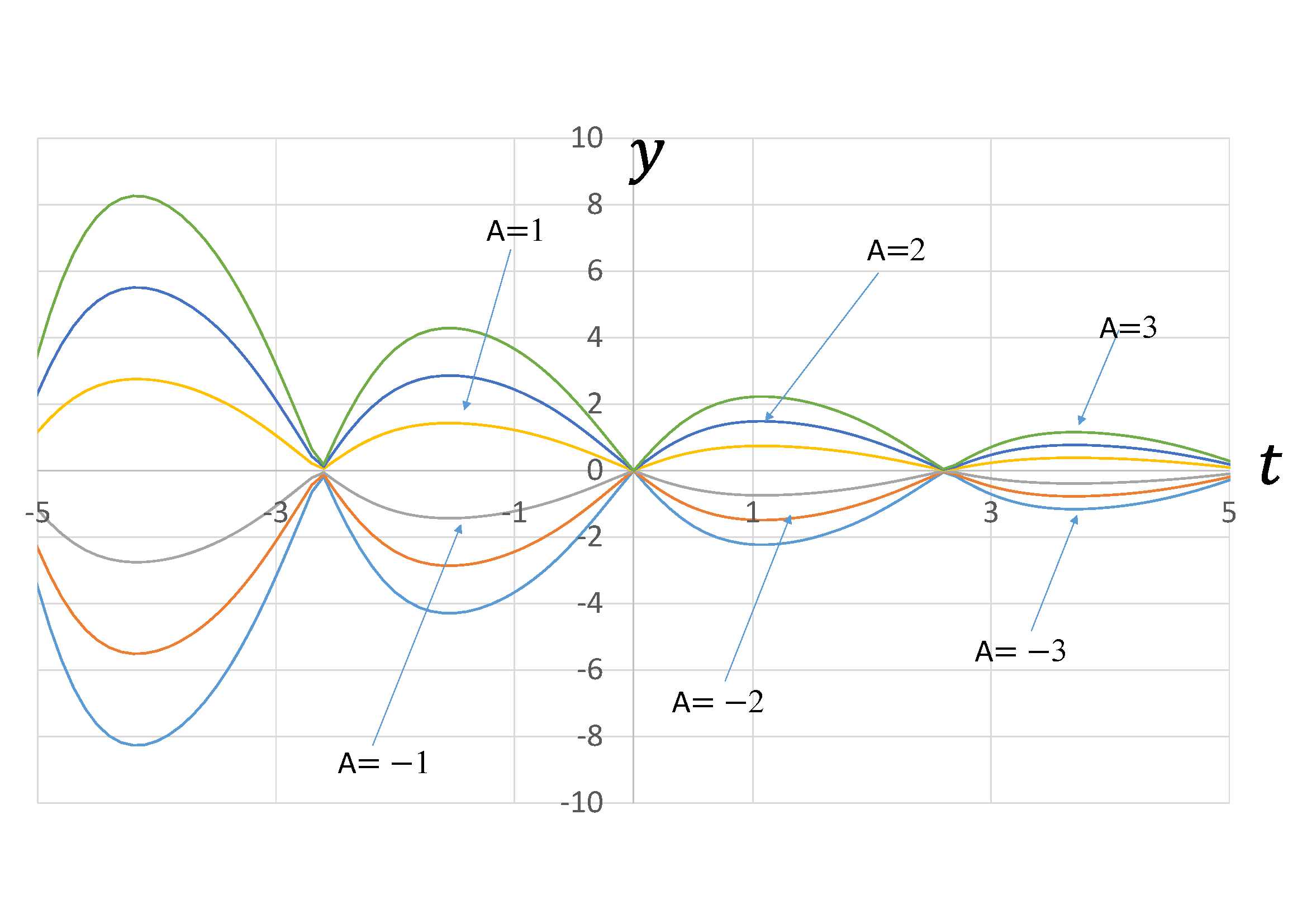}
\caption{Extended exact solution 2 with respect to variation of the parameter $A$ related to the amplitude($\omega=1$)}
\label{fig832}      
\end{center}
\end{figure*}

\begin{figure*}[tb]
\begin{center}
\includegraphics[width=0.70 \textwidth]{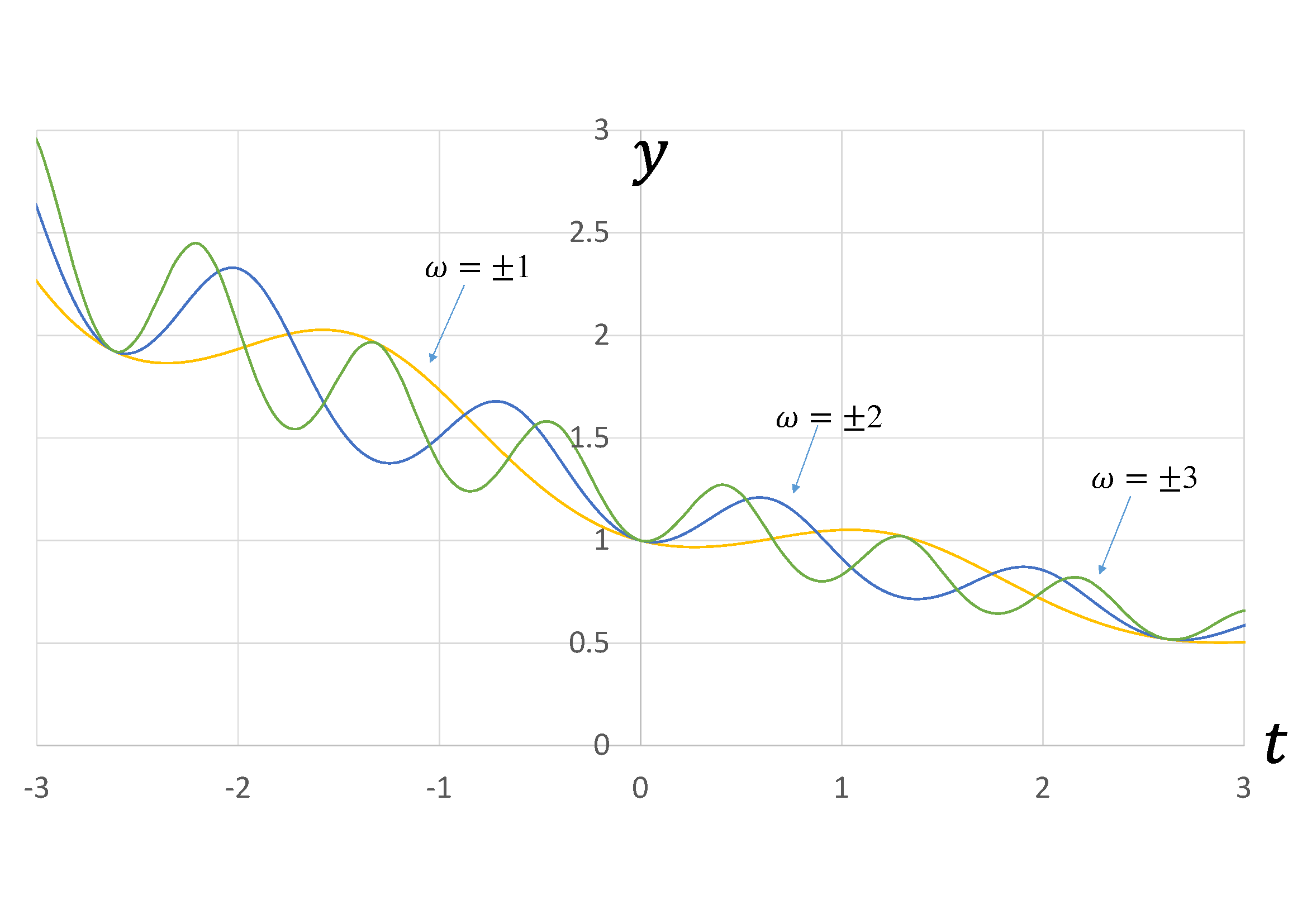}
\caption{Extended exact solution 3 with respect to variation of the parameter $\omega$ related to the period ($A=1$)}
\label{fig833}      
\end{center}
\end{figure*}

\begin{figure*}[tb]
\begin{center}
\includegraphics[width=0.70 \textwidth]{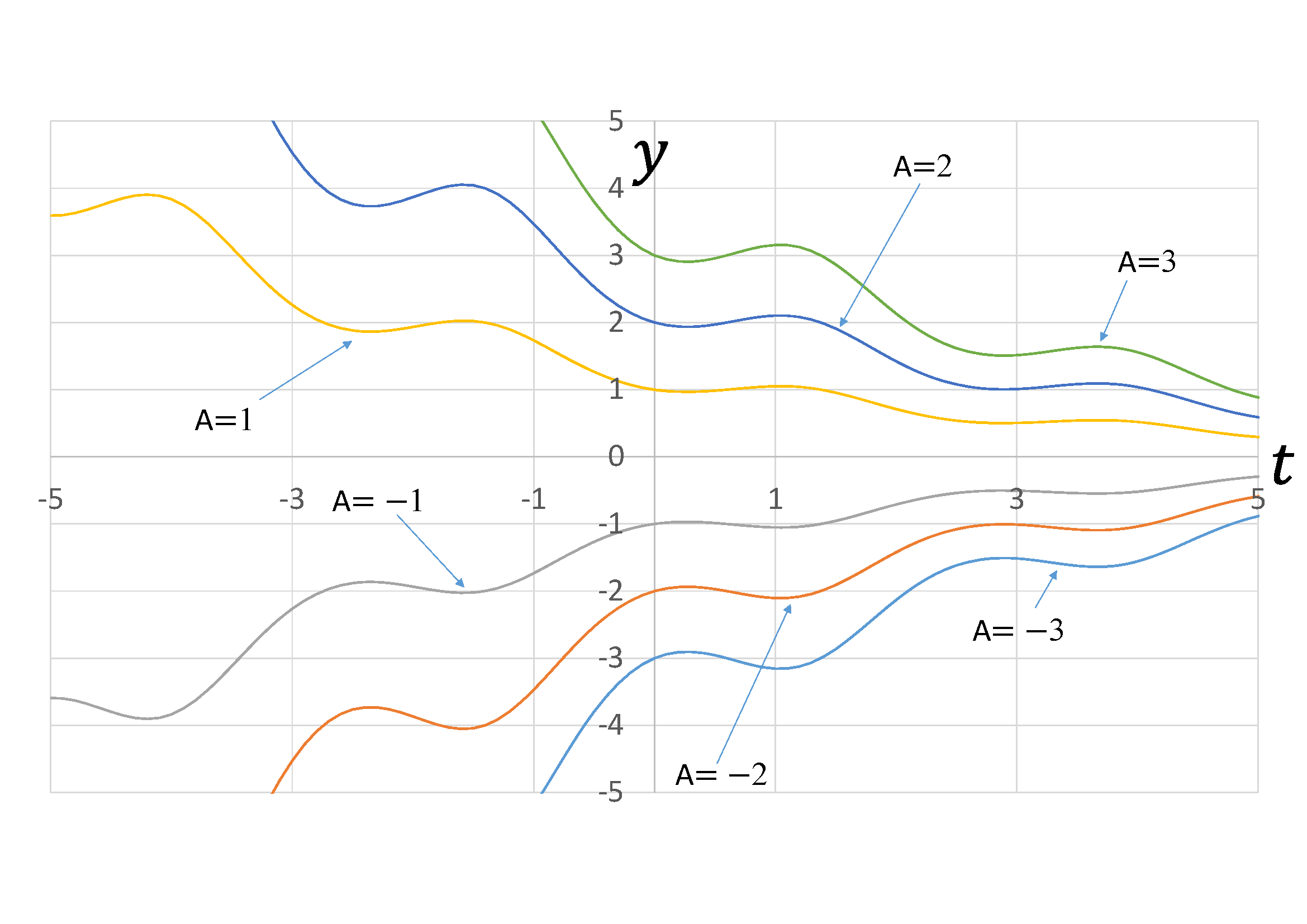}
\caption{Extended exact solution 3 with respect to variation of the parameter $A$ related to the amplitude($\omega=1$)}
\label{fig834}      
\end{center}
\end{figure*}

\begin{figure*}[tb]
\begin{center}
\includegraphics[width=0.70 \textwidth]{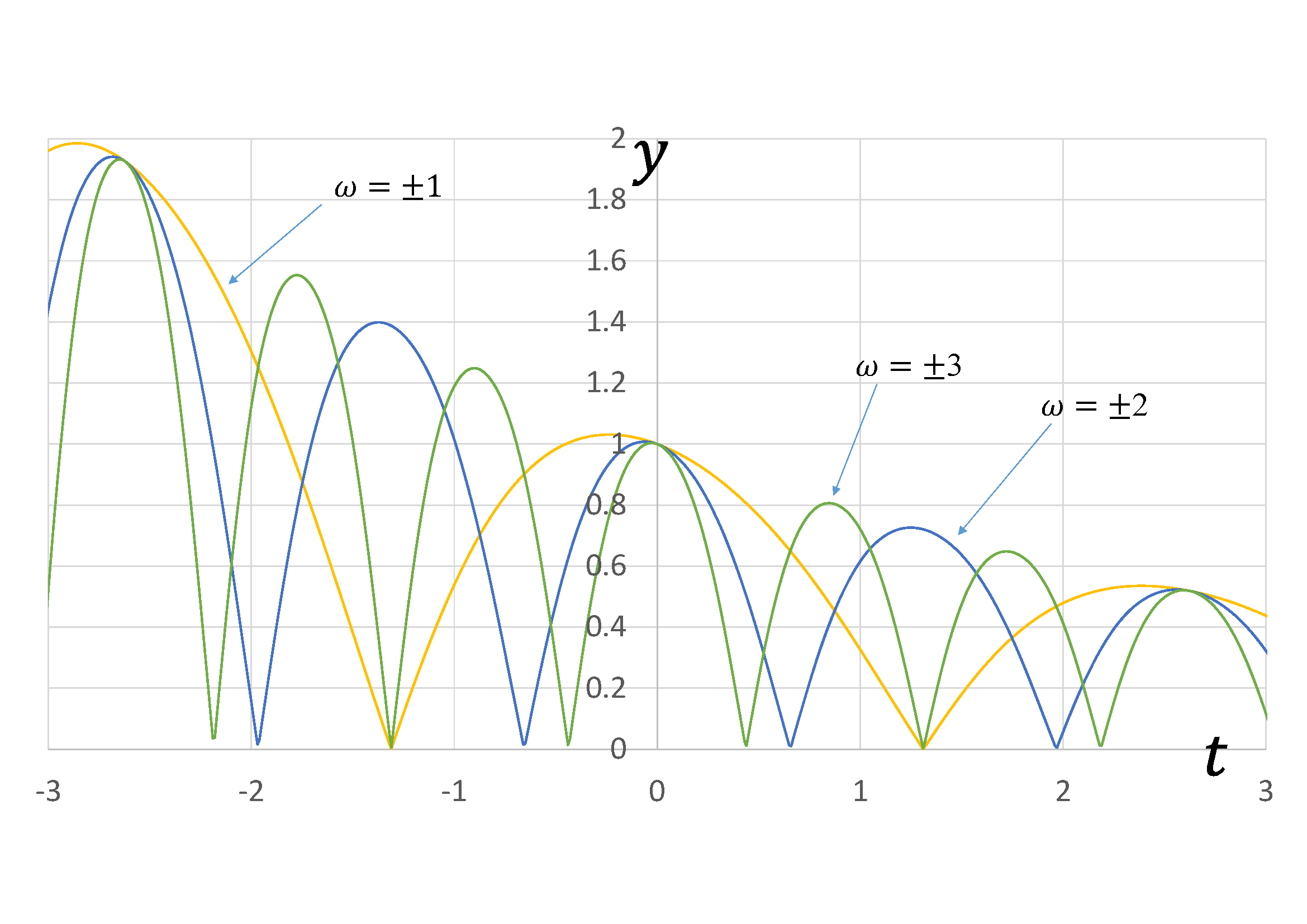}
\caption{Extended exact solution 4 with respect to variation of the parameter $\omega$ related to the period ($A=1$)}
\label{fig835}      
\end{center}
\end{figure*}

\begin{figure*}[tb]
\begin{center}
\includegraphics[width=0.70 \textwidth]{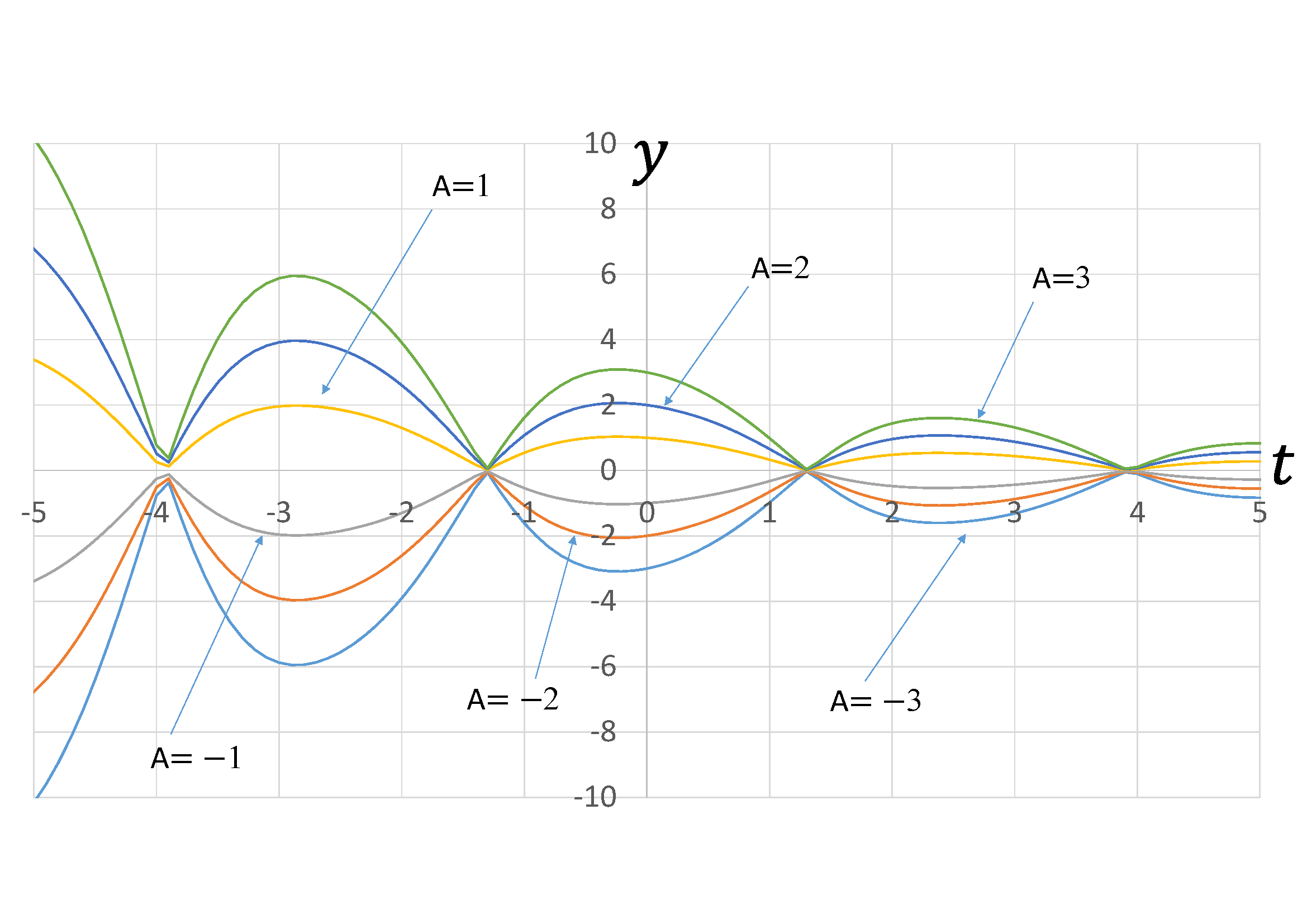}
\caption{Extended exact solution 4 with respect to variation of the parameter $A$ related to the amplitude($\omega=1$)}
\label{fig836}      
\end{center}
\end{figure*}

\begin{figure*}[tb]
\begin{center}
\includegraphics[width=0.70 \textwidth]{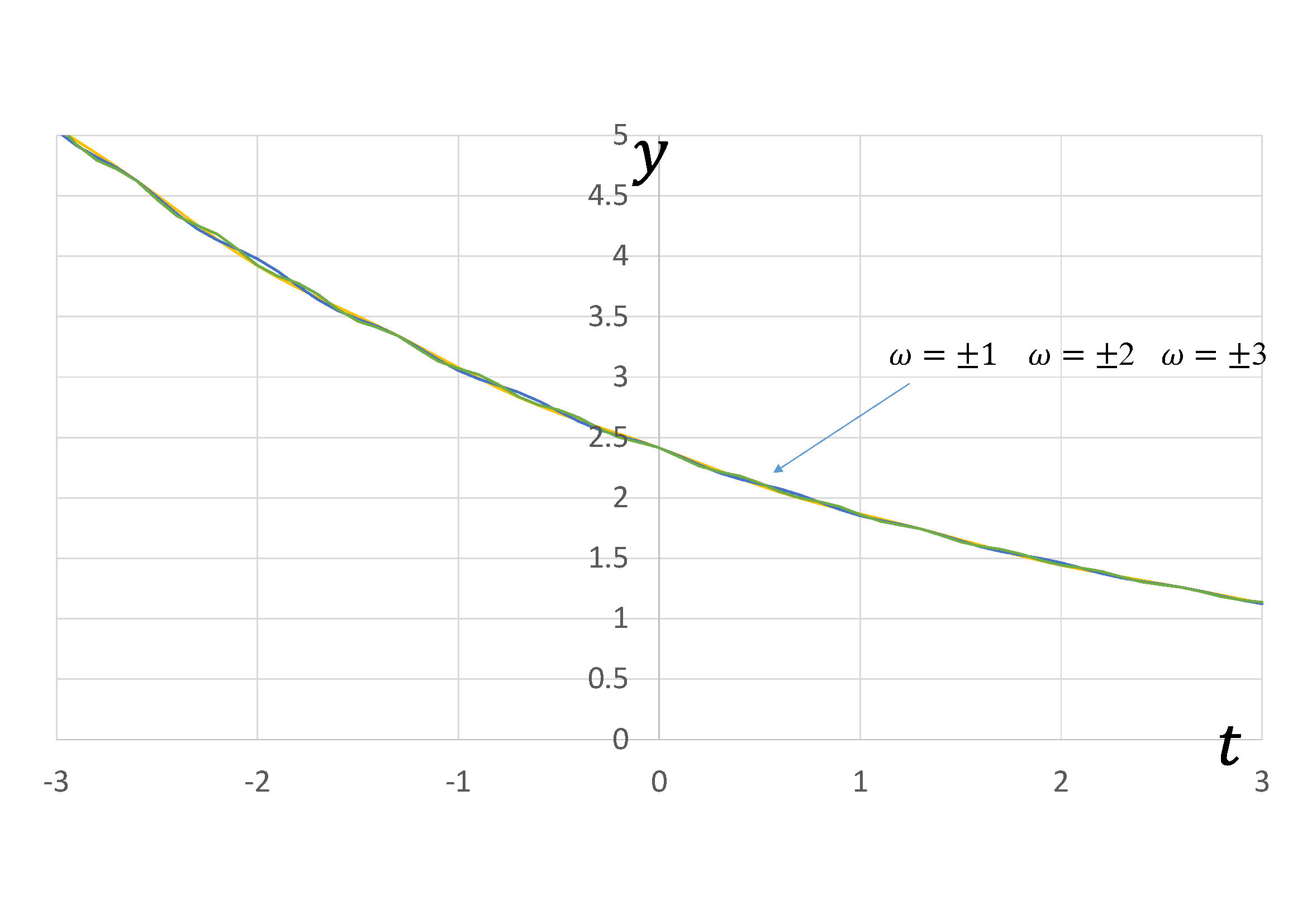}
\caption{Extended exact solution 5 with respect to variation of the parameter $\omega$ related to the period ($A=1$)}
\label{fig837}      
\end{center}
\end{figure*}

\begin{figure*}[tb]
\begin{center}
\includegraphics[width=0.70 \textwidth]{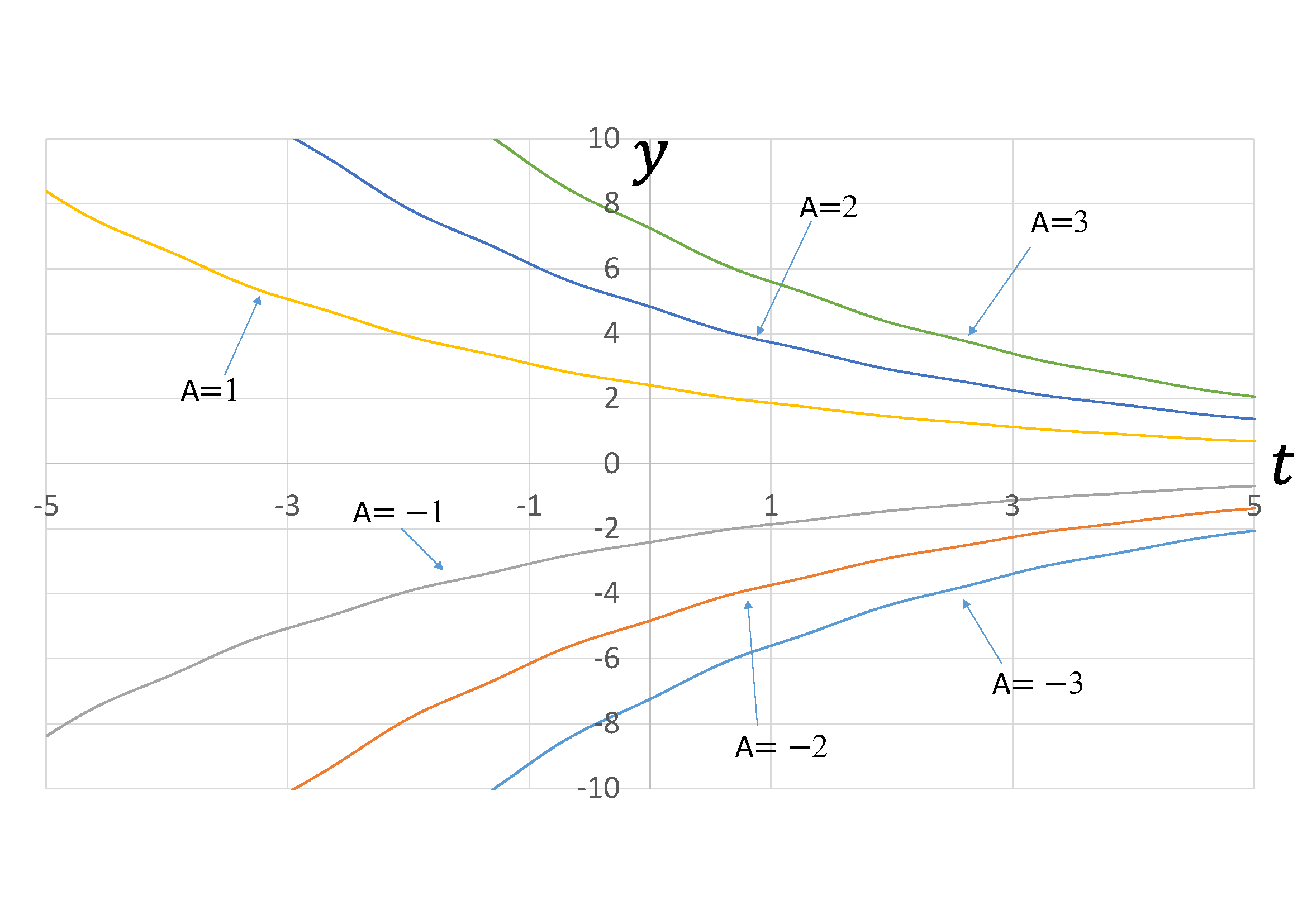}
\caption{Extended exact solution 5 with respect to variation of the parameter $A$ related to the amplitude($\omega=1$)}
\label{fig838}      
\end{center}
\end{figure*}

\begin{figure*}[tb]
\begin{center}
\includegraphics[width=0.70 \textwidth]{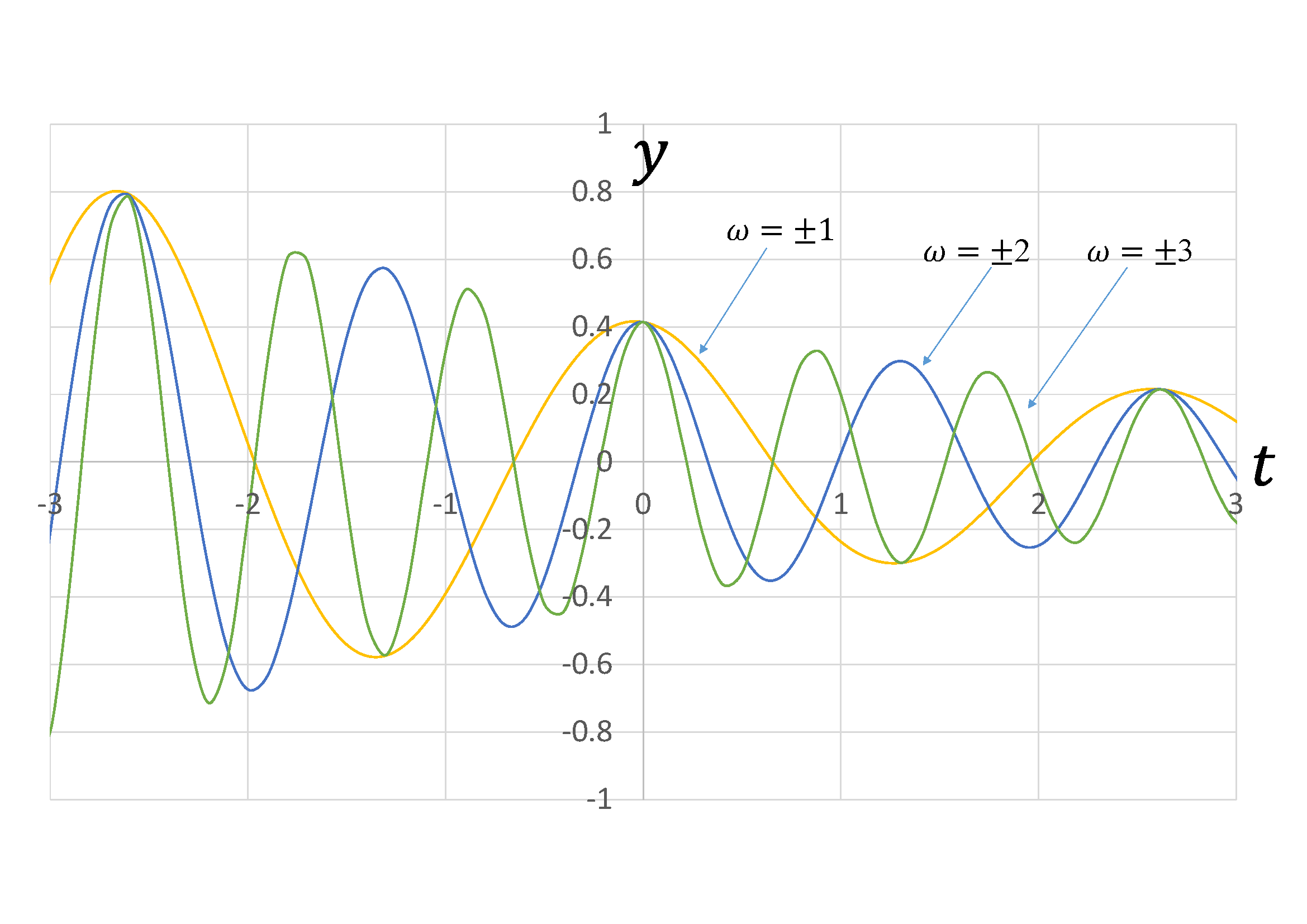}
\caption{Extended exact solution 6 with respect to variation of the parameter $\omega$ related to the period ($A=1$)}
\label{fig839}      
\end{center}
\end{figure*}

\begin{figure*}[tb]
\begin{center}
\includegraphics[width=0.70 \textwidth]{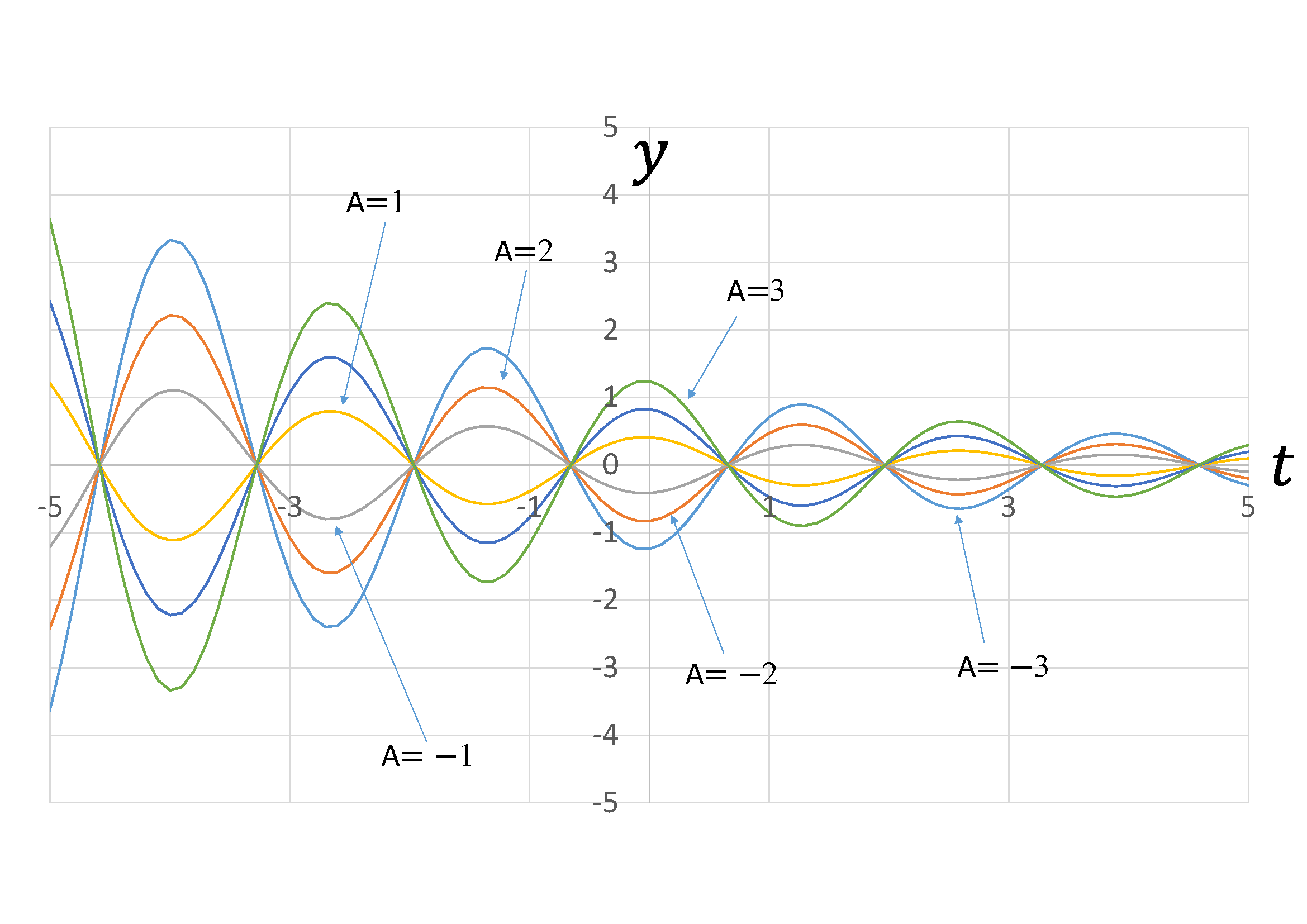}
\caption{Extended exact solution 6 with respect to variation of the parameter $A$ related to the amplitude($\omega=1$)}
\label{fig840}      
\end{center}
\end{figure*}

\begin{figure*}[tb]
\begin{center}
\includegraphics[width=0.70 \textwidth]{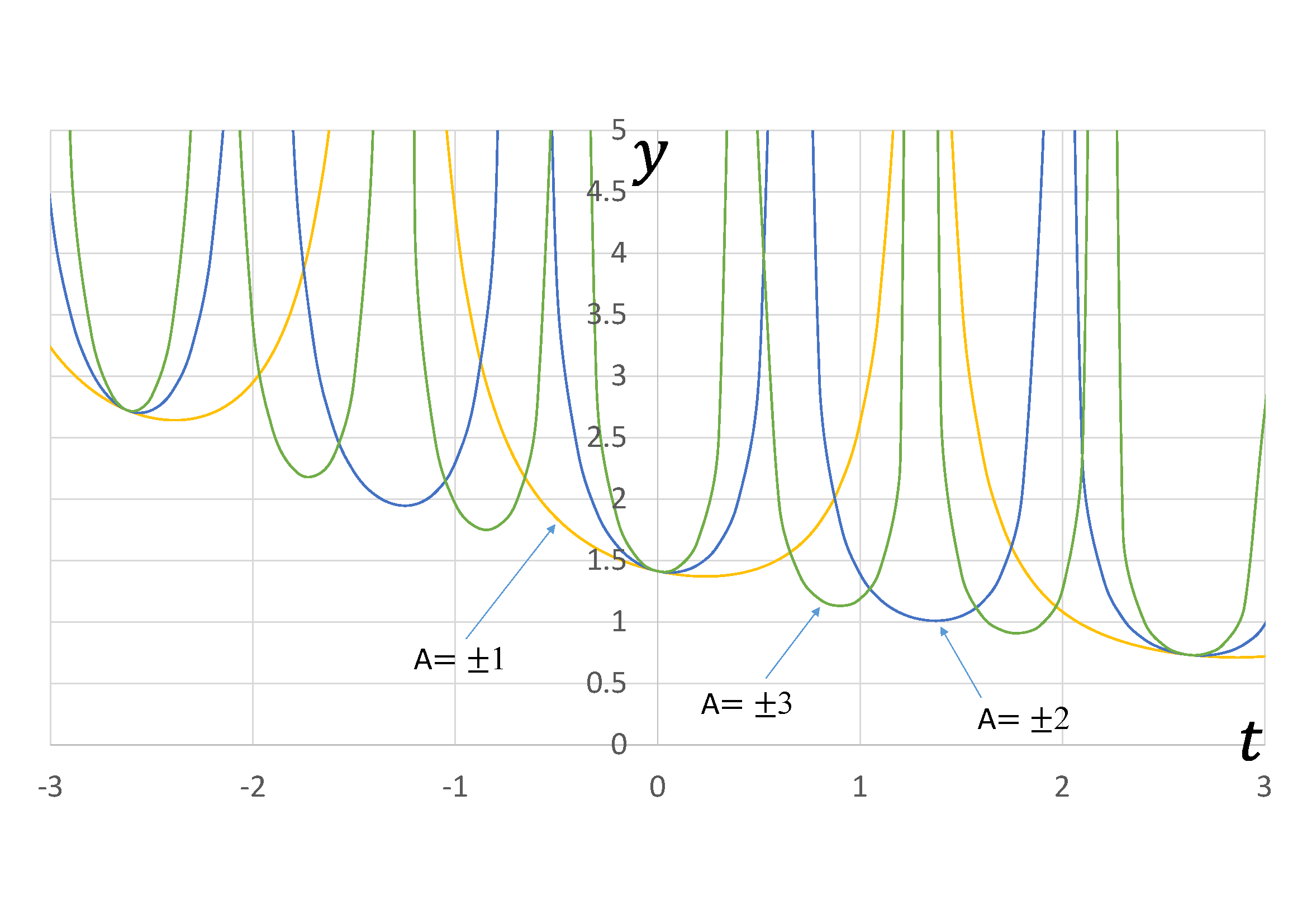}
\caption{Extended exact solution 7 with respect to variation of the parameter $\omega$ related to the period ($A=1$)}
\label{fig841}      
\end{center}
\end{figure*}

\begin{figure*}[tb]
\begin{center}
\includegraphics[width=0.70 \textwidth]{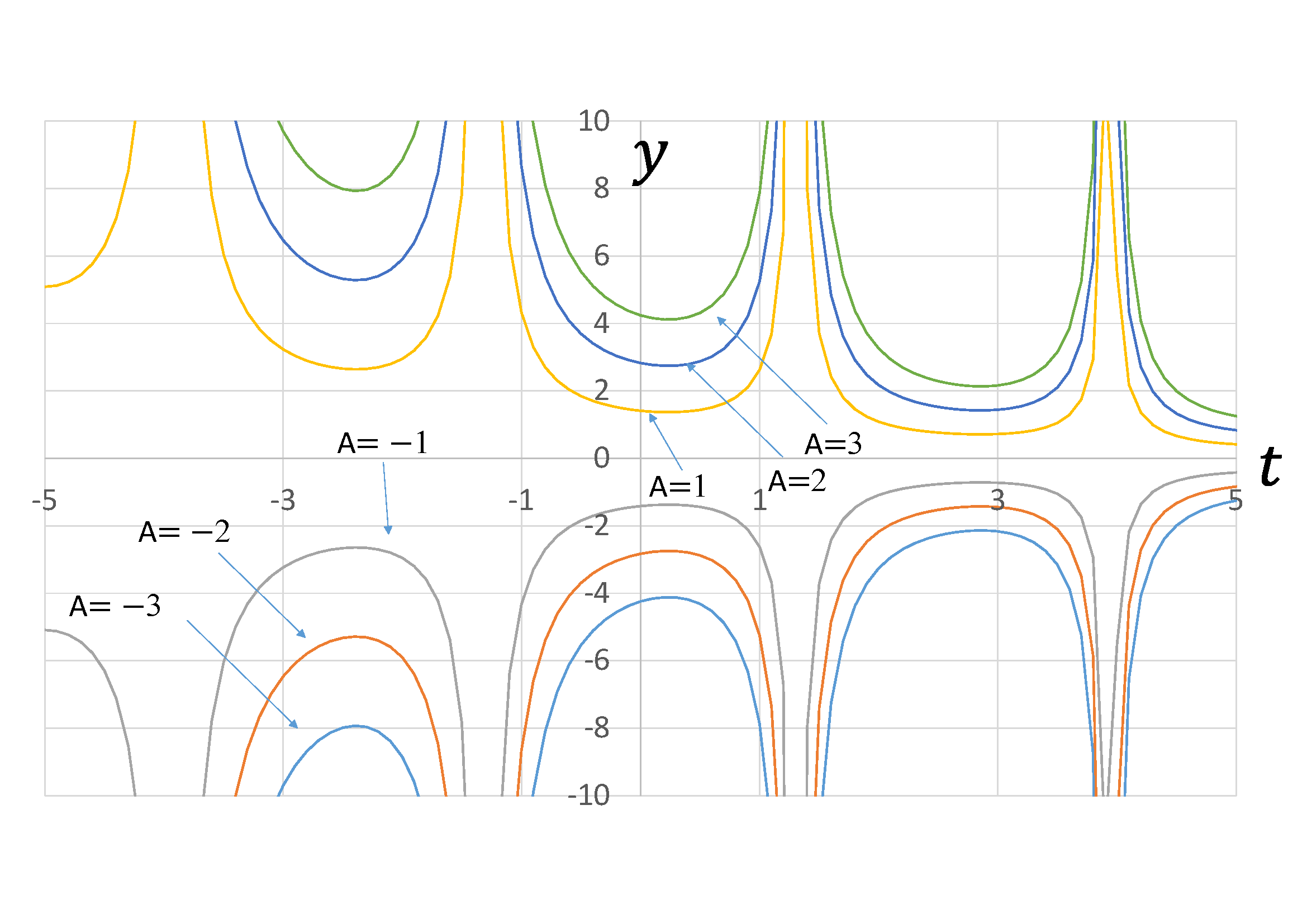}
\caption{Extended exact solution 7 with respect to variation of the parameter $A$ related to the amplitude($\omega=1$)}
\label{fig842}      
\end{center}
\end{figure*}

\begin{figure*}[tb]
\begin{center}
\includegraphics[width=0.70 \textwidth]{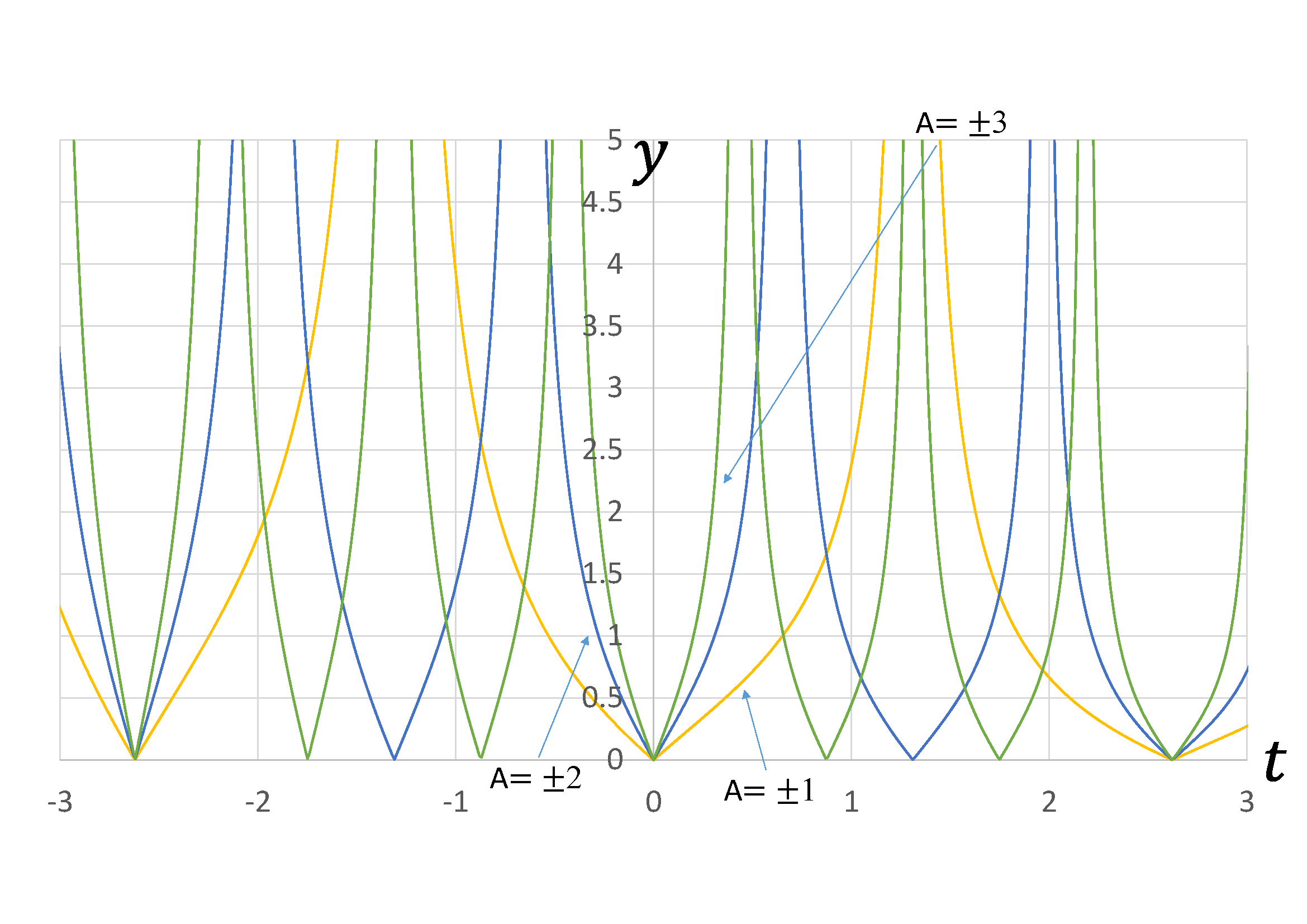}
\caption{Extended exact solution 8 with respect to variation of the parameter $\omega$ related to the period ($A=1$)}
\label{fig843}      
\end{center}
\end{figure*}

\begin{figure*}[tb]
\begin{center}
\includegraphics[width=0.70 \textwidth]{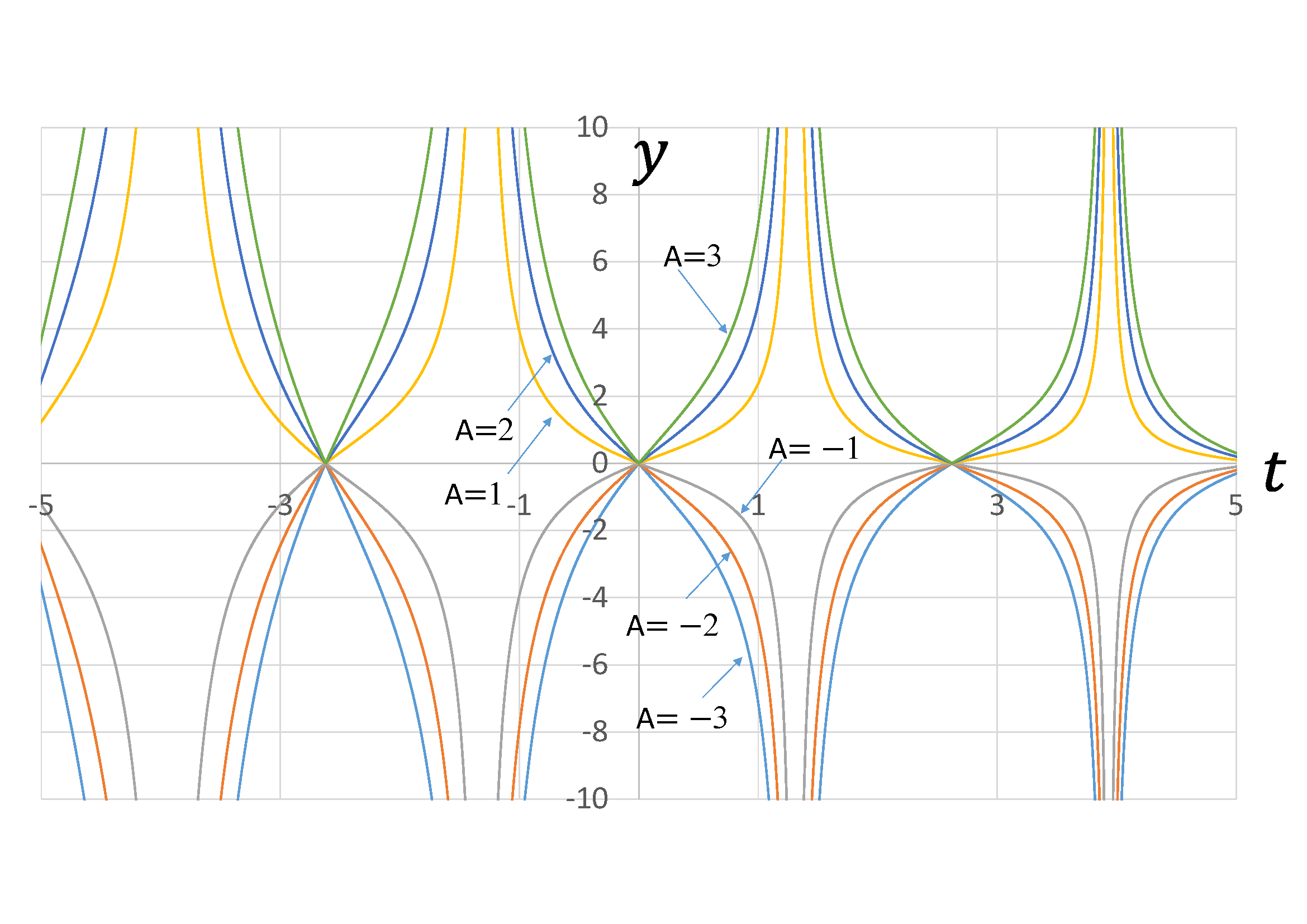}
\caption{Extended exact solution 8 with respect to variation of the parameter $A$ related to the amplitude($\omega=1$)}
\label{fig844}      
\end{center}
\end{figure*}

\begin{figure*}[tb]
\begin{center}
\includegraphics[width=0.70 \textwidth]{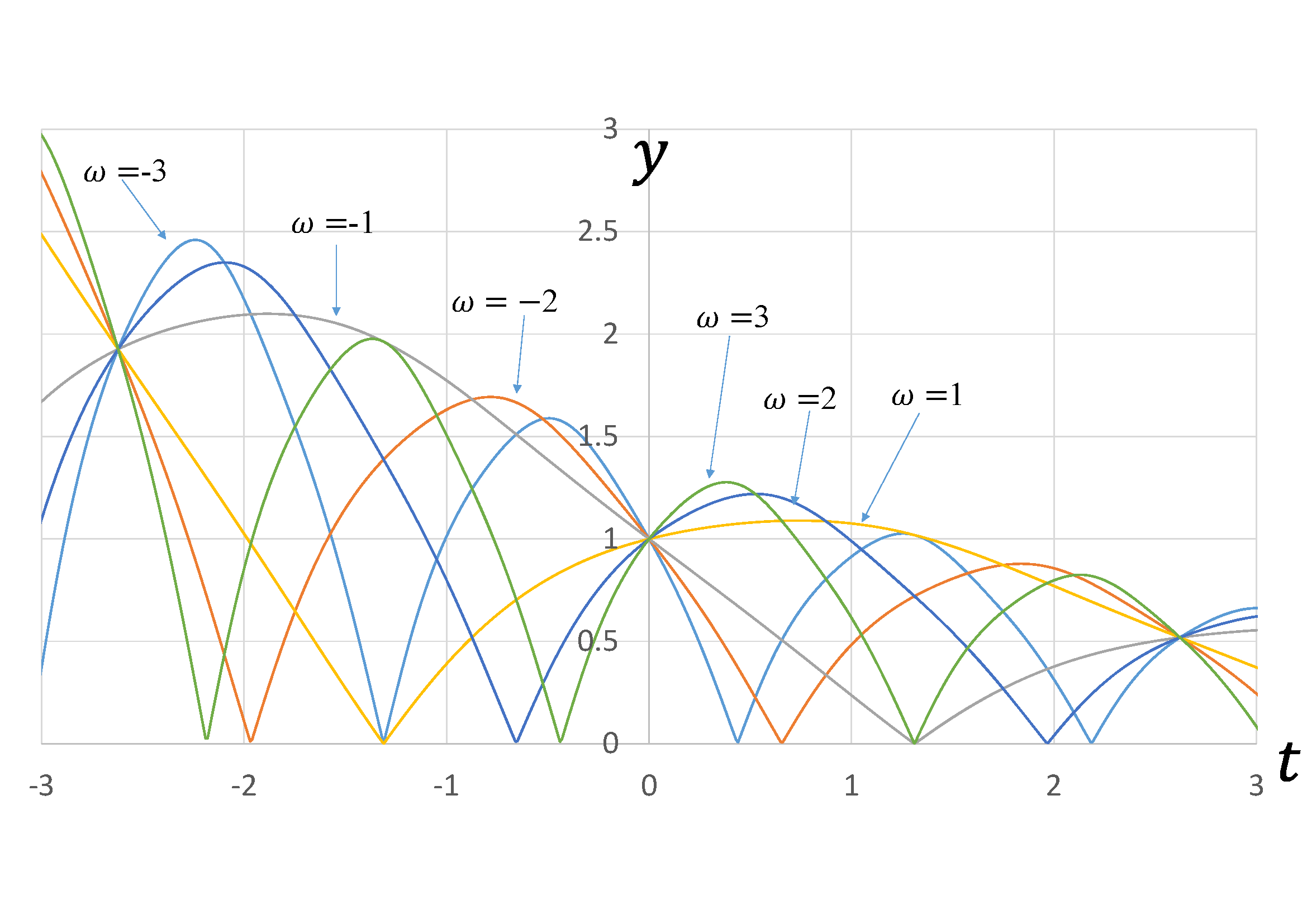}
\caption{Extended exact solution 9 with respect to variation of the parameter $\omega$ related to the period ($A=1$)}
\label{fig845}      
\end{center}
\end{figure*}

\begin{figure*}[tb]
\begin{center}
\includegraphics[width=0.70 \textwidth]{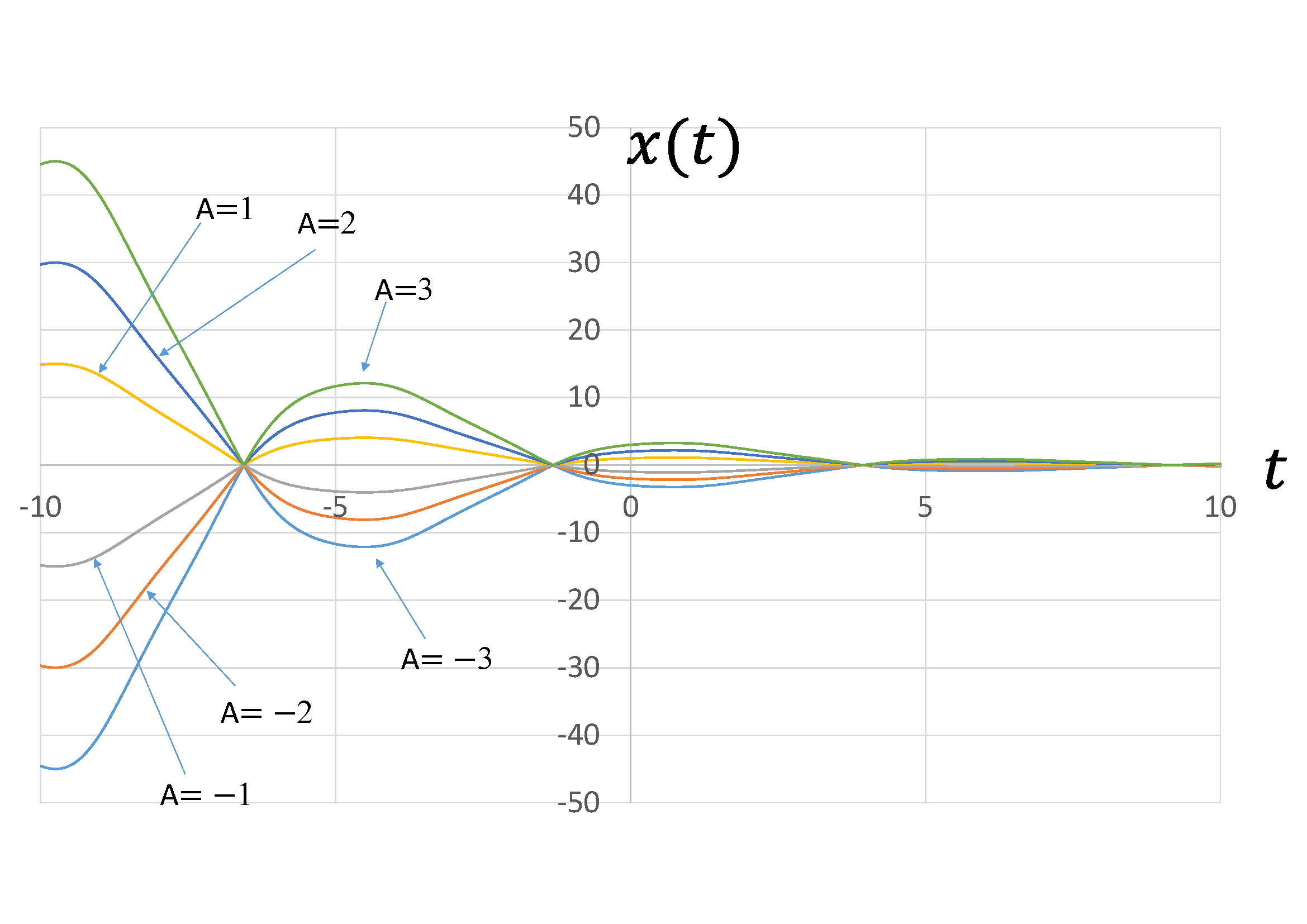}
\caption{Extended exact solution 9 with respect to variation of the parameter $A$ related to the amplitude($\omega=1$)}
\label{fig846}      
\end{center}
\end{figure*}

\begin{figure*}[tb]
\begin{center}
\includegraphics[width=0.70 \textwidth]{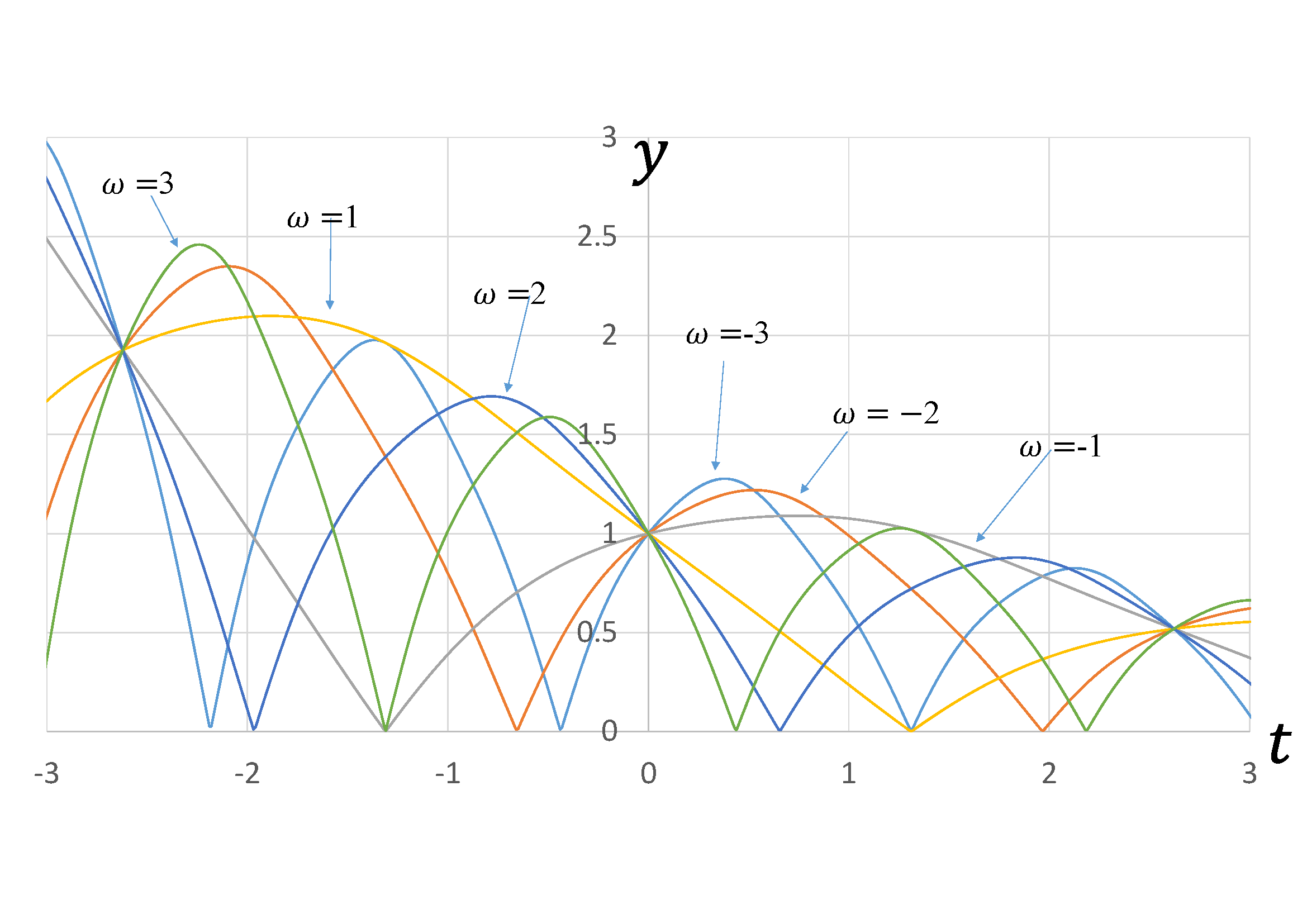}
\caption{Extended exact solution 10 with respect to variation of the parameter $\omega$ related to the period ($A=1$)}
\label{fig847}      
\end{center}
\end{figure*}

\begin{figure*}[tb]
\begin{center}
\includegraphics[width=0.70 \textwidth]{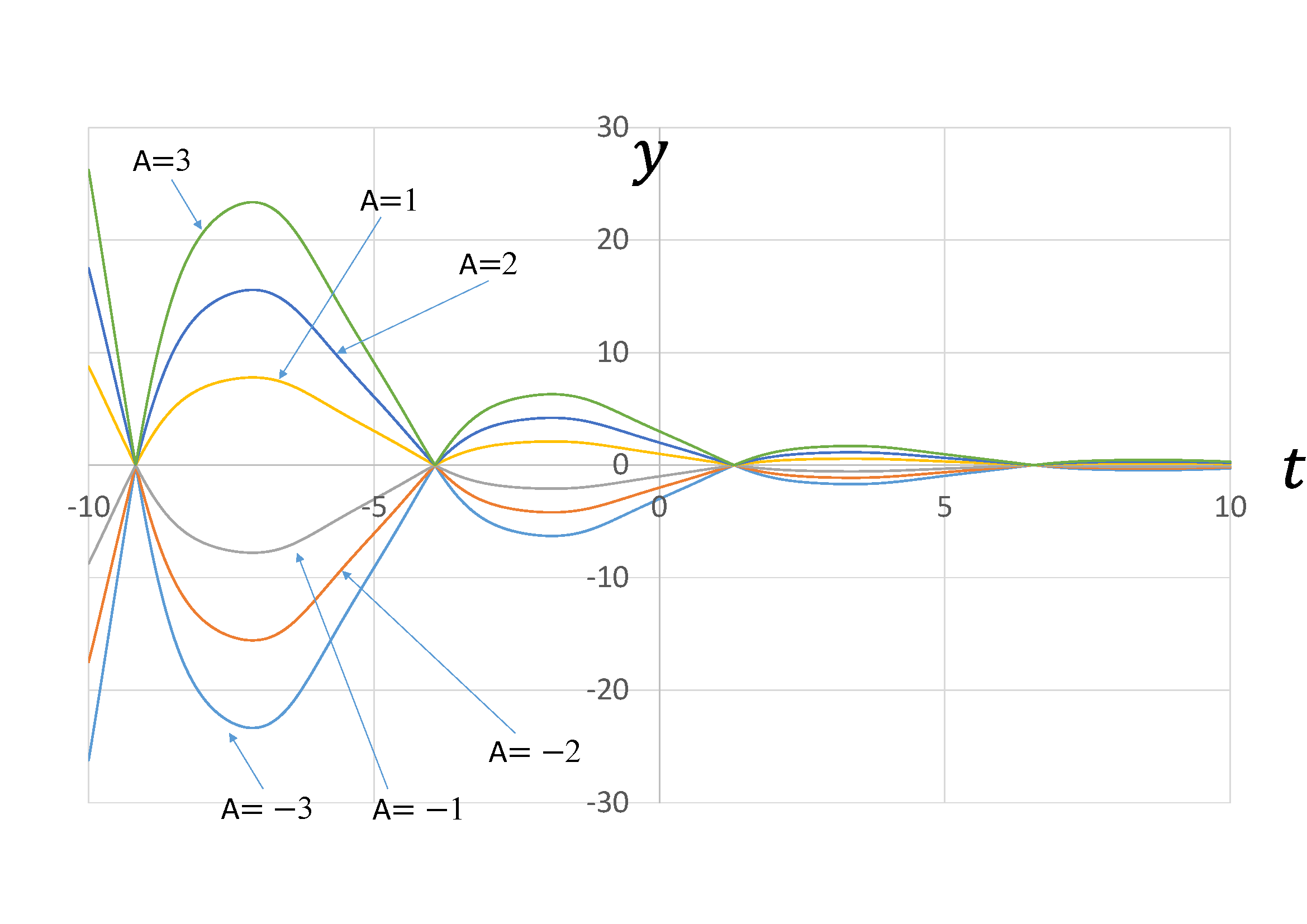}
\caption{Extended exact solution 10 with respect to variation of the parameter $A$ related to the amplitude($\omega=1$)}
\label{fig848}      
\end{center}
\end{figure*}

\begin{figure*}[tb]
\begin{center}
\includegraphics[width=0.70 \textwidth]{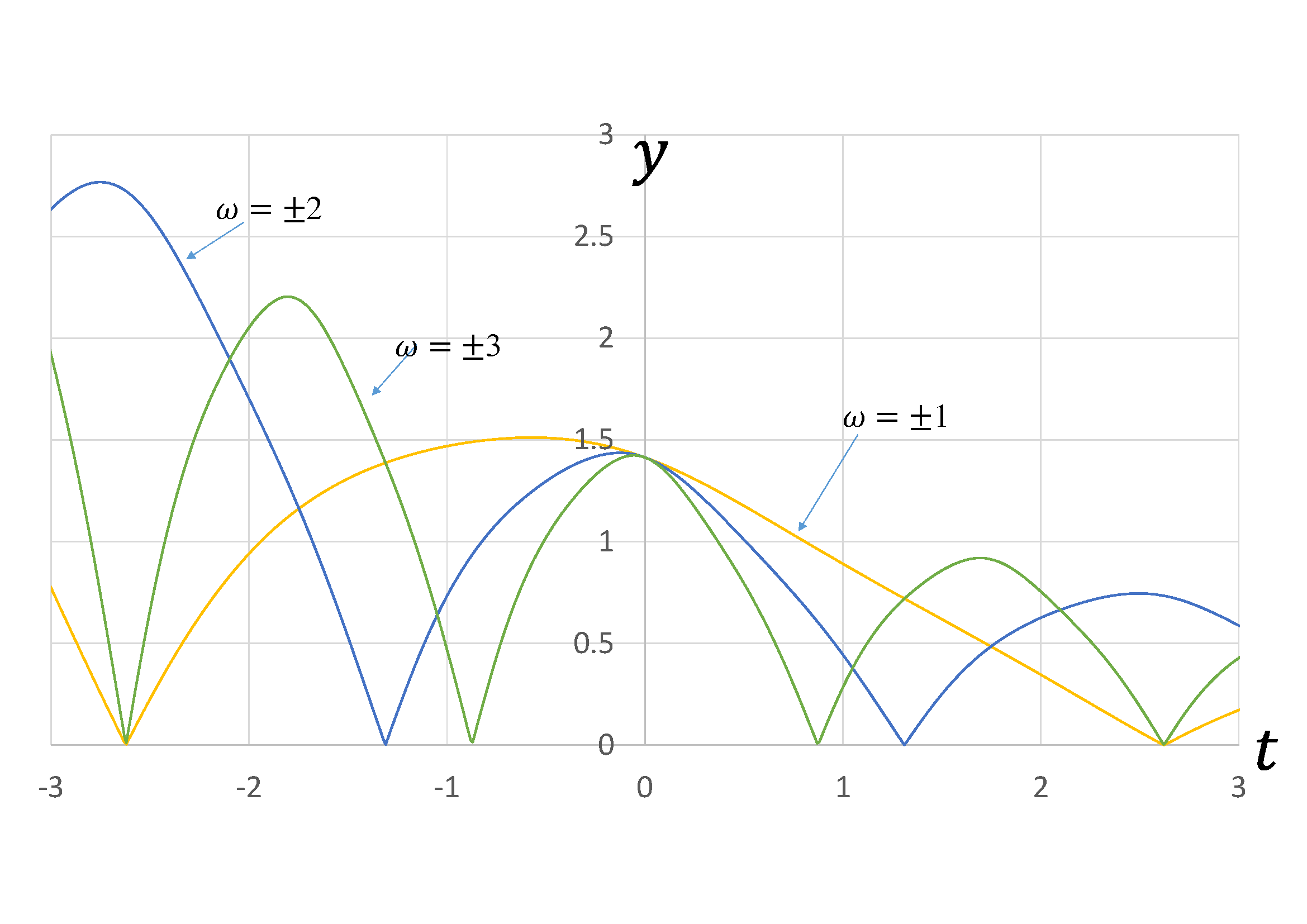}
\caption{Extended exact solution 11 with respect to variation of the parameter $\omega$ related to the period ($A=1$)}
\label{fig849}      
\end{center}
\end{figure*}

\begin{figure*}[tb]
\begin{center}
\includegraphics[width=0.70 \textwidth]{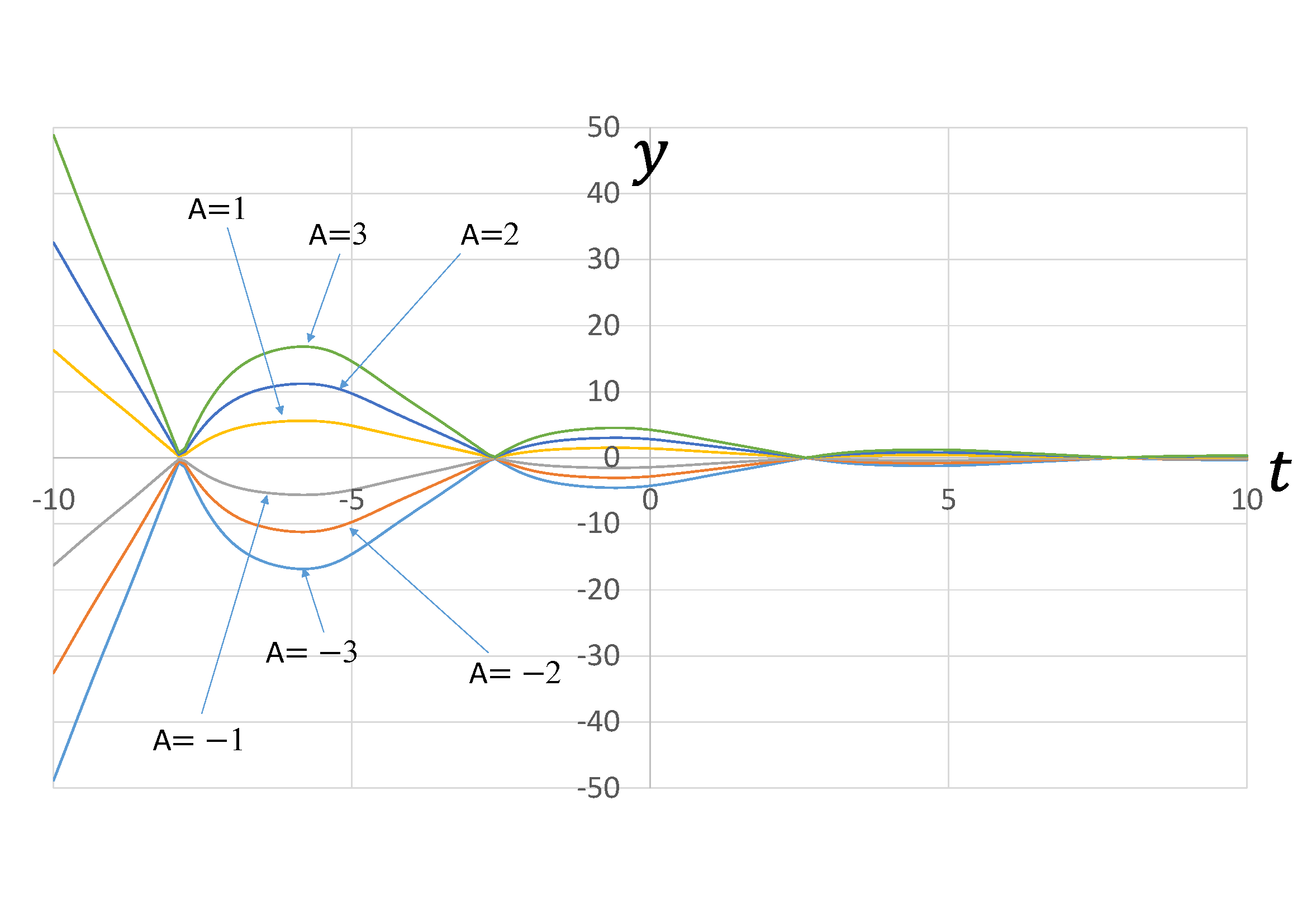}
\caption{Extended exact solution 11 with respect to variation of the parameter $A$ related to the amplitude($\omega=1$)}
\label{fig850}      
\end{center}
\end{figure*}

\begin{figure*}[tb]
\begin{center}
\includegraphics[width=0.70 \textwidth]{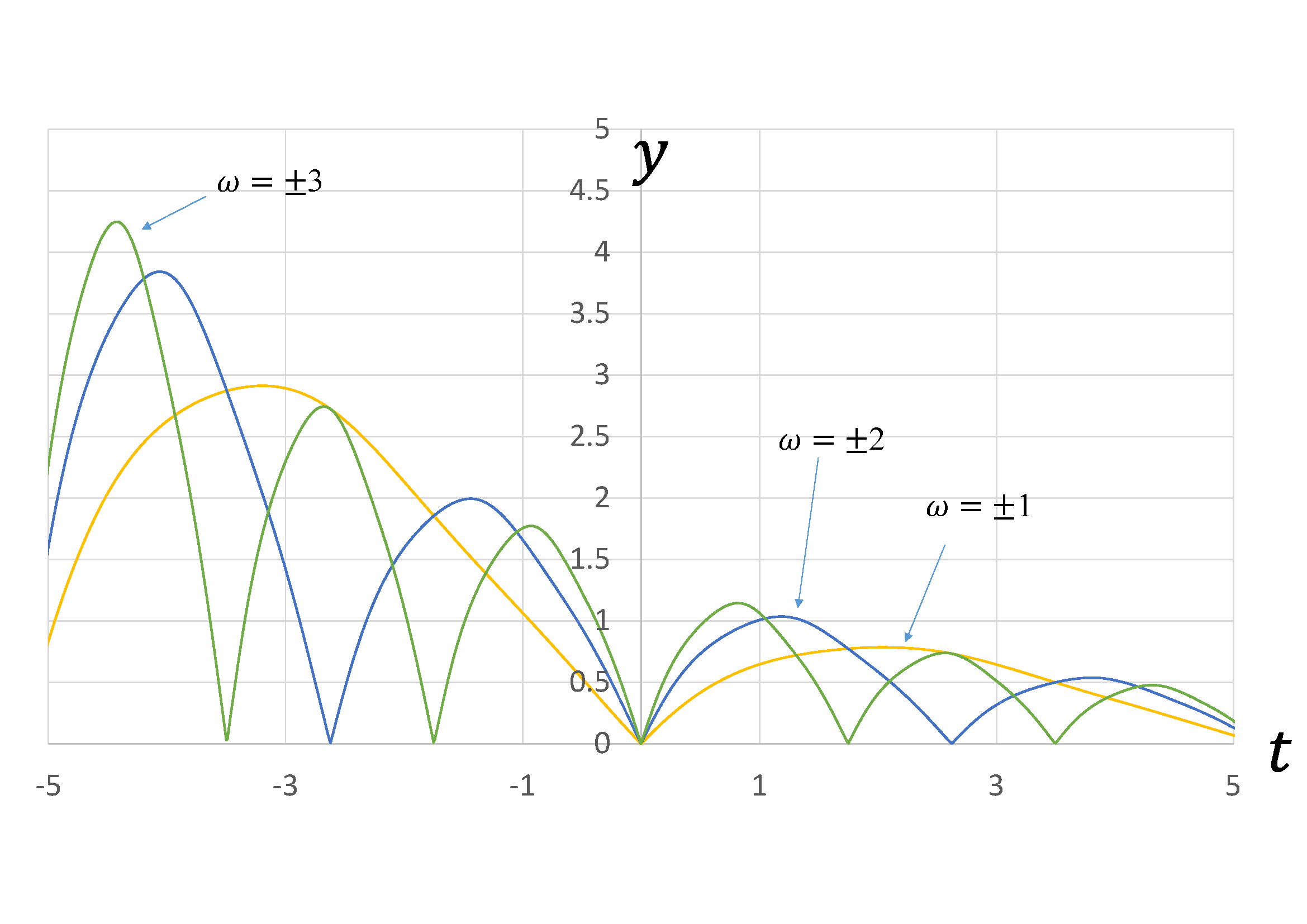}
\caption{Extended exact solution 12 with respect to variation of the parameter $\omega$ related to the period ($A=1$)}
\label{fig851}      
\end{center}
\end{figure*}

\begin{figure*}[tb]
\begin{center}
\includegraphics[width=0.70 \textwidth]{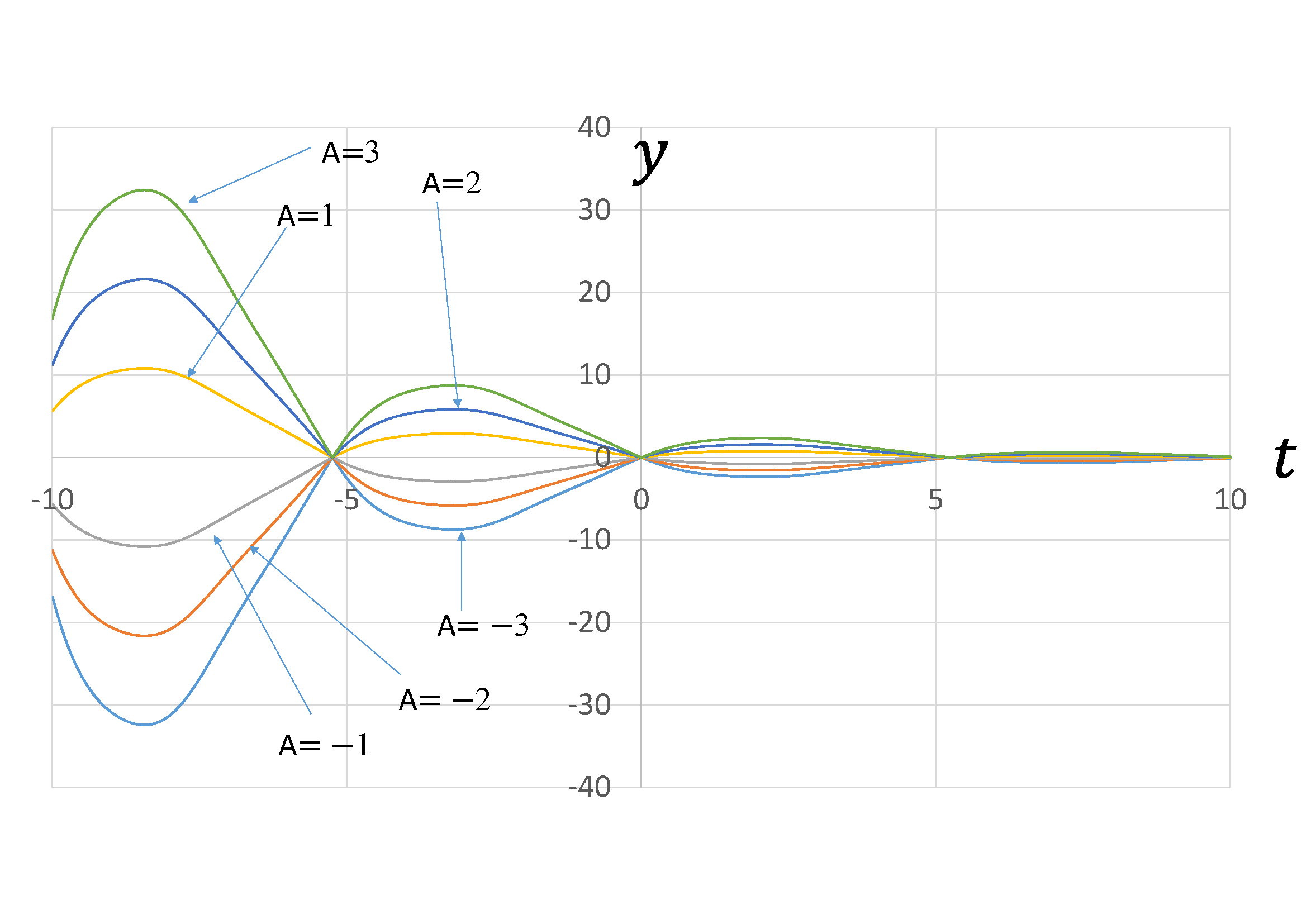}
\caption{Extended exact solution 12 with respect to variation of the parameter $A$ related to the amplitude($\omega=1$)}
\label{fig852}      
\end{center}
\end{figure*}

\begin{figure*}[tb]
\begin{center}
\includegraphics[width=0.70 \textwidth]{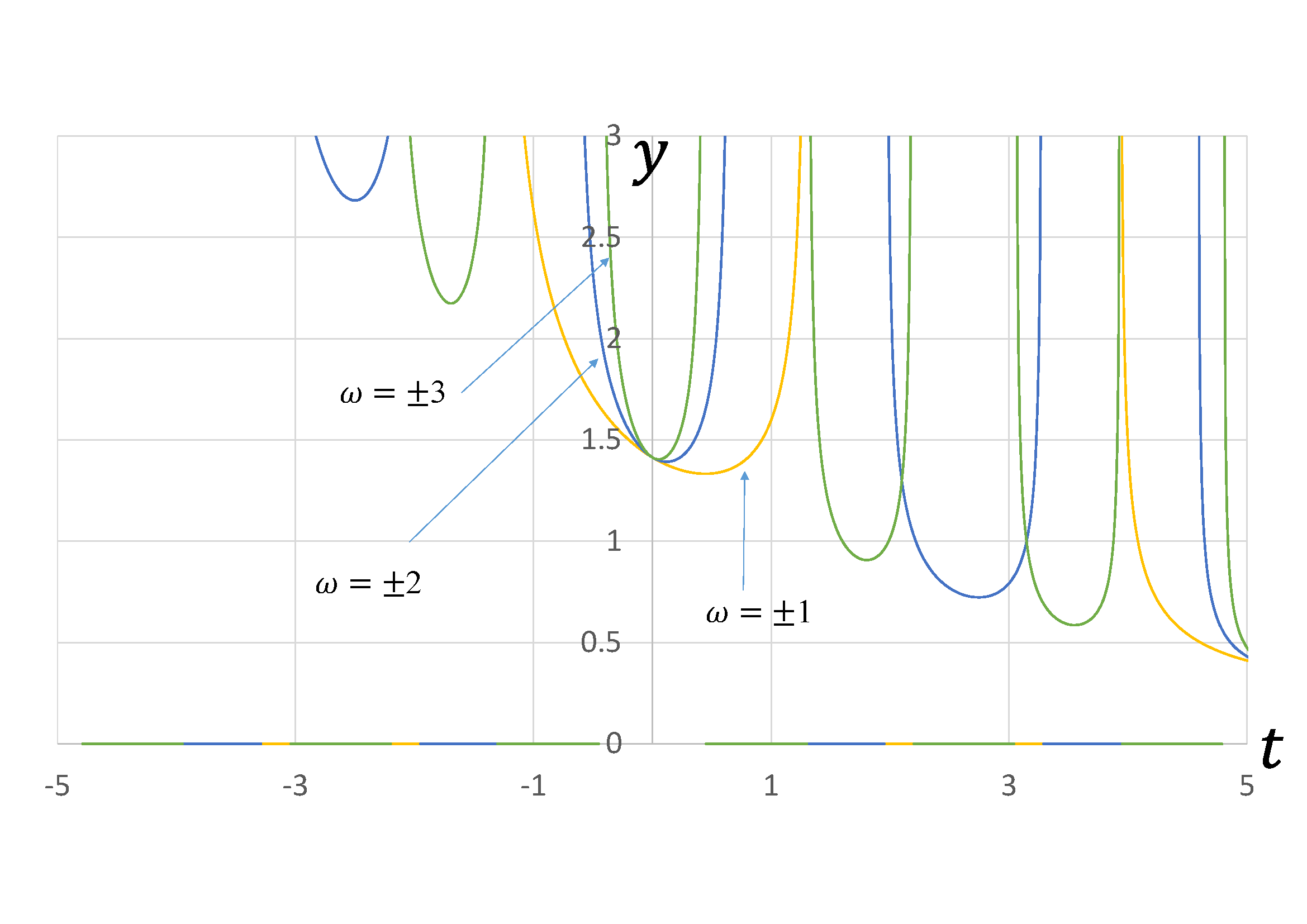}
\caption{Extended exact solution 13 with respect to variation of the parameter $\omega$ related to the period ($A=1$)}
\label{fig853}      
\end{center}
\end{figure*}

\begin{figure*}[tb]
\begin{center}
\includegraphics[width=0.70 \textwidth]{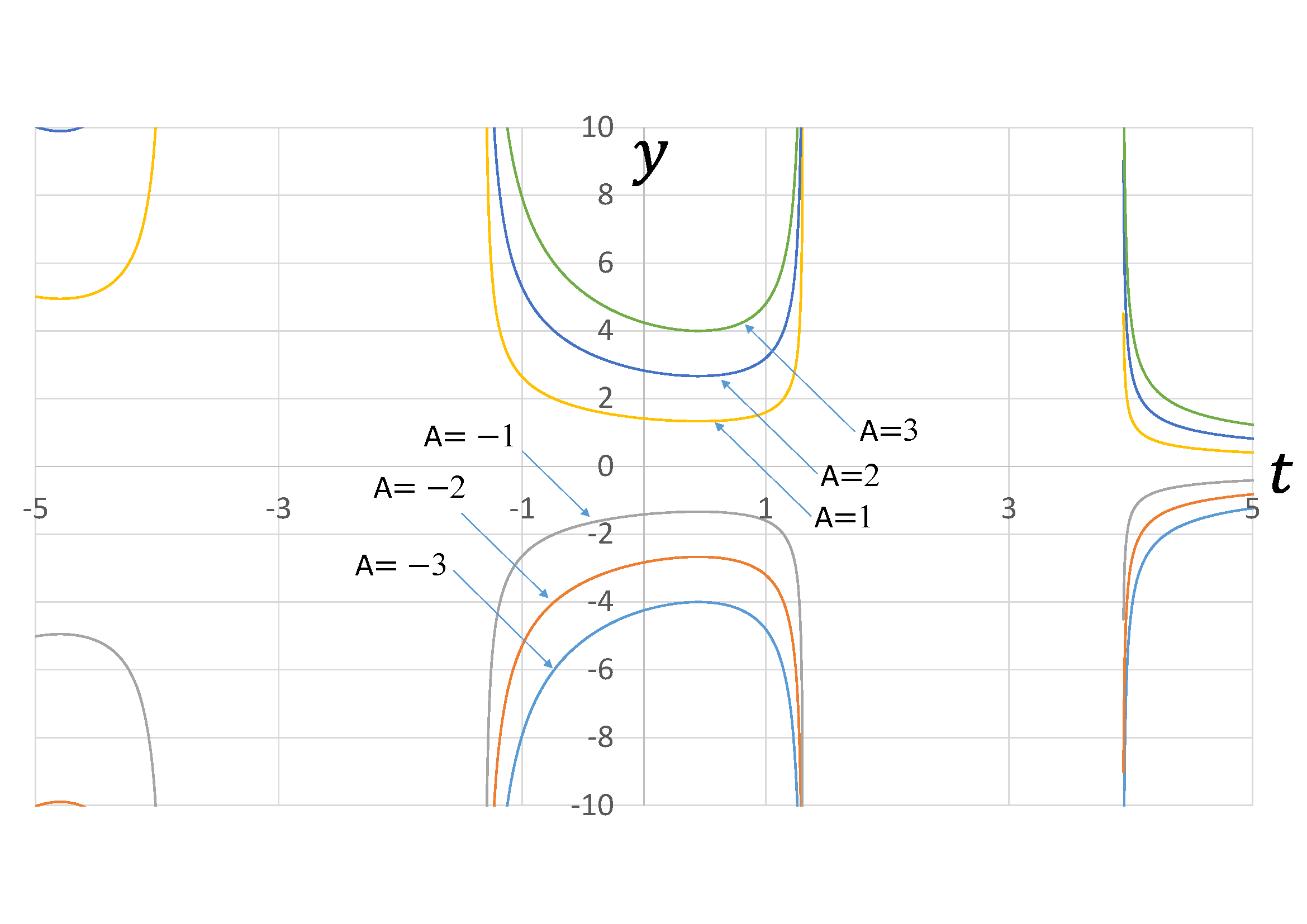}
\caption{Extended exact solution 13 with respect to variation of the parameter $A$ related to the amplitude($\omega=1$)}
\label{fig854}      
\end{center}
\end{figure*}

\begin{figure*}[tb]
\begin{center}
\includegraphics[width=0.70 \textwidth]{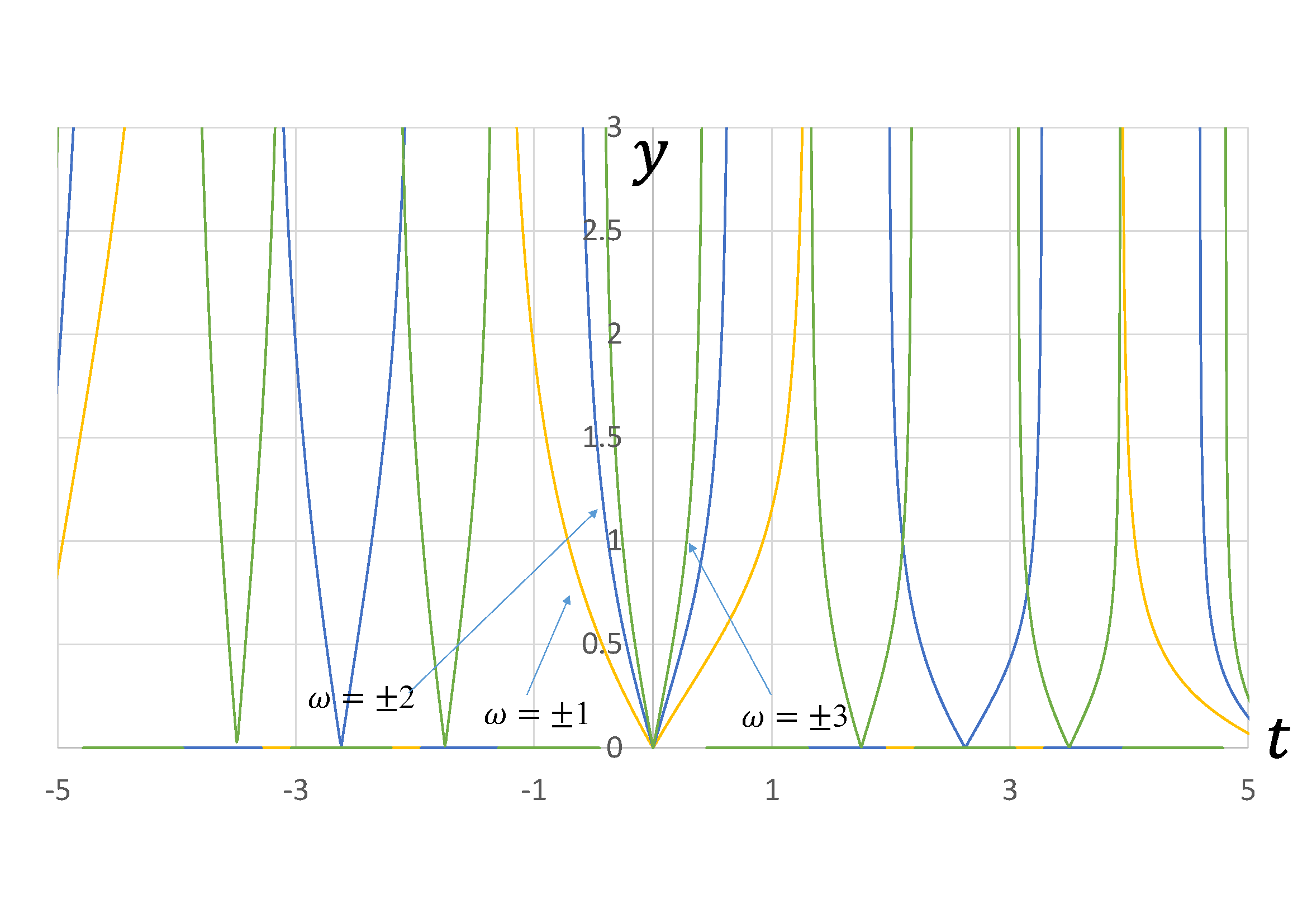}
\caption{Extended exact solution 14 with respect to variation of the parameter $\omega$ related to the period ($A=1$)}
\label{fig855}      
\end{center}
\end{figure*}

\begin{figure*}[tb]
\begin{center}
\includegraphics[width=0.70 \textwidth]{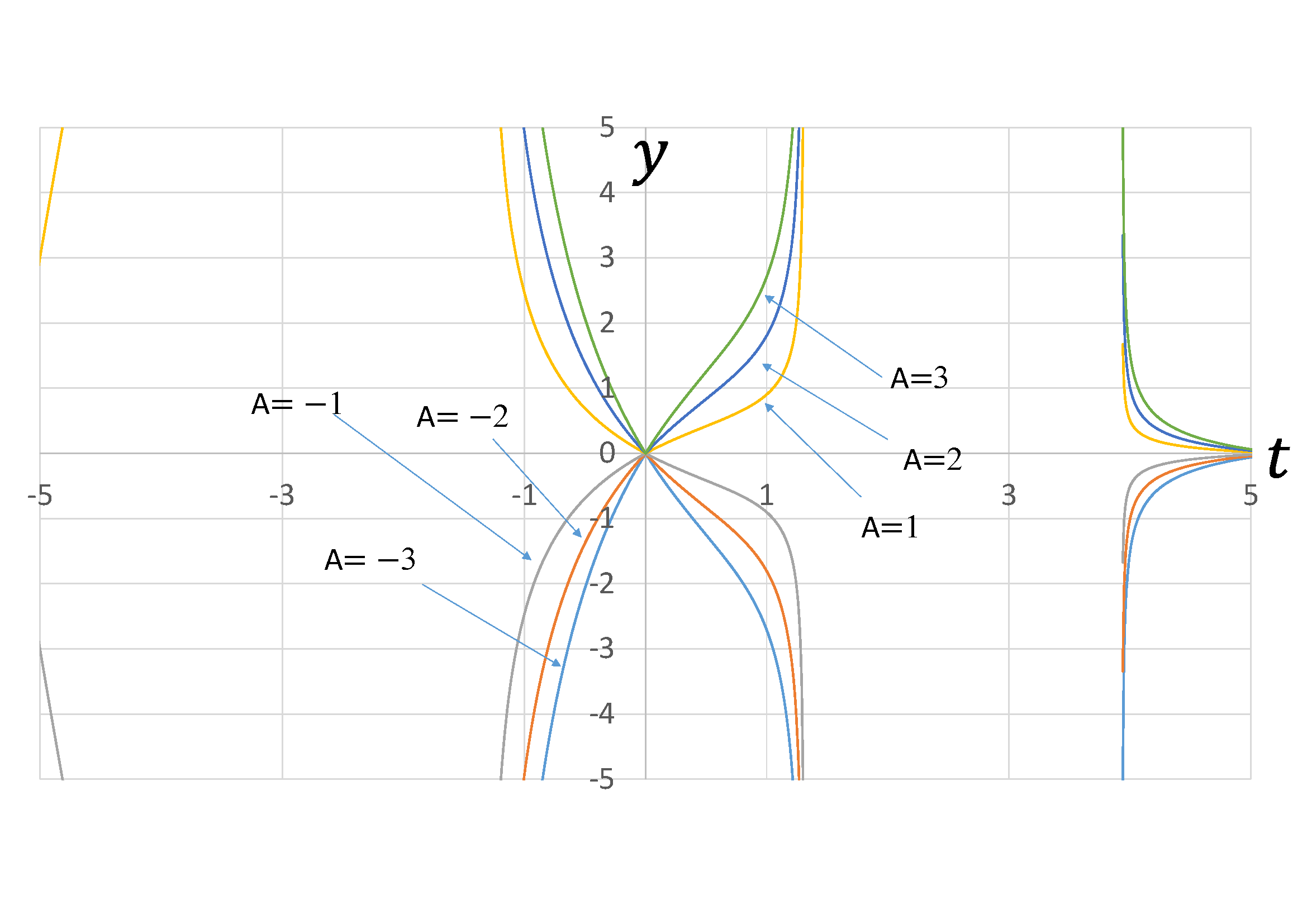}
\caption{Extended exact solution 14 with respect to variation of the parameter $A$ related to the amplitude($\omega=1$)}
\label{fig856}      
\end{center}
\end{figure*}

The extended exact solution satisfies the Duffing equation including the first derivative term as shown in Eq. (\ref{6.1}). As shown in Eqs (\ref{6.6})-(\ref{6.18}), extended exact solutions (Eqs (\ref{4.1})-(\ref{5.9})) are created by multiplying the exact solutions by the exponential function. Therefore, the amplitude of the wave gradually increases or decreases according to the time axis.
The present approach based on the leaf functions can be used to generalized models of dampers. Duffing equation that generalizes the classical harmonic vibration are presented in the Eqs (\ref{4.2})-(\ref{5.9}). By using a suitable transformation and an associated undamped harmonic vibration with similar frequency, we can obtain a system in the form of a generalized Duffing model which included an additional first derivative. The Duffing system with a damping term described by an exponent term are examined. The analysis results of the extended exact solution obtained by extending the undamped system or the divergent system are shown in Fig. \ref{fig829} - Fig. \ref{fig856}.
The vertical axis  and the horizontal axis represents the variable $x(t)$ and the variable $t$, respectively. In order to perform a numerical analysis concretely, the damped Duffing equation (\ref{6.6})-(\ref{7.9}) are numerically analyzed with the conditions $ \beta(t)=\frac{1}{2}$ and $c=0$. The amplitude coefficients $A$ in the Duffing equation term and the initial condition are varied  as shown in Fig. \ref{fig829} - Fig. \ref{fig856}. Only the amplitude of the wave can be controlled while keeping the wave period $\omega$. On the other hand, only the period $\omega$ of the wave can be controlled while keeping the amplitude $A$. Using the leaf functions, the term coefficents of the cubic-quint damped duffing equaiotn can be split in two elements both the the amplitude $A$ and the period $\omega$.
 In this analysis result, $\beta(t)$ is a constant. However, it is also possible to analyze with the function $\beta(t)$ as a polynomial of the variable $t$. This makes it possible to create waves of various damping systems or divergent systems with leaf functions.


\maxdeadcycles=1000
\clearpage

\section{Conclusion}
In a previous study, the exact solutions of the cubic Duffing equation were presented by using the leaf function.
However, the calculations of the exact solutions are complicated because they are described by  using the leaf function for the phase of the trigonometric function. In this study, the exact solutions of the cubic Duffing equation are presented by using only leaf functions without using trigonometric functions.
Furthermore, the exact solutions are derived for the cubic--quintic Duffing equation, which has strong nonlinearity. 

The conclusions are as follows:

\begin{itemize}
\item The exact solutions of the cubic Duffing equation and cubic--quintic Duffing equation are presented by using leaf functions. Some of these combinations of factors based on factorization of the first derivatives can be applied as the exact solutions of the cubic-quintic Duffing equation.
\item The exact solutions have two parameters: parameter $A$ related to the amplitude and parameter $\omega$ related to the period. By using parameter $A$, the amplitude of the wave according to the cubic-quintic Duffing equation can be adjusted while the period of the wave is constant.  By using parameter $\omega$, the wave period can be adjusted while the wave amplitude is constant. The relations between these parameters and the coefficient of the term of the Duffing equation are clarified.
\item  Based on the exact solution that satisfies the Duffing equation, the extended exact solutions that satisfy the damping system Duffing equation and the divergence system Duffing equation are presented. By applying the equation transformation  with respect to the undamped Duffing equation, it is possible to derive the extended Duffing equation consisting of the second derivative, the first derivative and the polynomial. The exact solutions that satisfies the extended Duffing equation can be derived by both the exponential function and the leaf function. These exact solutions becomes very simple. The wave based on the extended Duffing equation is visualized by numerical analysis. The wave is damped while keeping the period of the wave by the Duffing equation. Since the waves by the numerical analysis completely match the wave by the exact solutions, the numerical analysis are consistent with the exact solutions.

\end{itemize}

\appendix
\section{ Symbols }
The parameter $m$ represents an integer. The symbols $\eta_n$, $\pi_n$ and $\zeta_n$ represent constants (see Appendices P, Q and R), respectively. The sign $\pm$ of the derivative $\pm \sqrt{1-(\mathrm{sleaf}_2(t))^4}$ of the leaf function sleaf$_2(t)$ depends on the domain of the phase $l$ \cite{Kaz_sl}. The leaf function cleaf$_2(t)$ and the hyperbolic leaf function cleafh$_2(t)$ are the same \cite{Kaz_cl} \cite{Kaz_clh}. Therefore, the differentiations are derived by dividing them into regions.

The symbol $\lim_{t  \downarrow a } f(t) $   represents the right-side limit. As the variable $t$ decreases, it approaches point $a$ from right to left on the $t$ axis of the graph. The symbol $\lim_{t \uparrow a }$ represents the left-side limit. As the variable $t$ increases, it approaches point $a$ from left to right on the $t$ axis of the graph. 

\section{ Ordinary differential equation of exact solution 1 }
\renewcommand{\theequation}{B.\arabic{equation}}
\setcounter{equation}{0}

Here we prove that exact solution 1 satisfies the Duffing equation.
The first derivatives of exact solution 1 are obtained as follows:

(i) In the case where $(2m-1) \pi_2 \leqq \omega t \leqq 2m \pi_2$, 
\begin{equation}
\begin{split}
\frac{\mathrm{d} x(t) }{\mathrm{d} t} 
&= A
\frac{ 2 \mathrm{cleaf}_2(\omega t) 
\left(\omega \sqrt{1-(\mathrm{cleaf}_2(\omega t))^4 } \right) }
{ 2 \sqrt{1+(\mathrm{cleaf}_2(\omega t))^2} } \\
&= A \omega \mathrm{cleaf}_2(\omega t) \sqrt{1-(\mathrm{cleaf}_2(\omega t))^2 }.
\label{B1}
\end{split}
\end{equation}

(ii) In the case where $2m \pi_2 \leqq \omega t \leqq (2m+1) \pi_2 $, 
\begin{equation}
\begin{split}
\frac{\mathrm{d} x(t) }{\mathrm{d} t} 
&= A
\frac{ 2\mathrm{cleaf}_2(\omega t) 
\left( -\omega\sqrt{1-(\mathrm{cleaf}_2(\omega t))^4 } \right) }
{ 2 \sqrt{1+(\mathrm{cleaf}_2(\omega t))^2} } \\
&= -A \omega \mathrm{cleaf}_2(\omega t) \sqrt{1-(\mathrm{cleaf}_2(\omega t))^2 }.
\label{B2}
\end{split}
\end{equation}

The second derivatives of exact solution 1 are obtained as follows:

(i) In the case where $(2m-1) \pi_2 \leqq \omega t \leqq 2m \pi_2$, 
\begin{equation}
\begin{split}
\frac{\mathrm{d}^2 x(t) }{\mathrm{d} t^2}
&=    A \omega^2 \left( \sqrt{1-(\mathrm{cleaf}_2(\omega t))^4 } \right) \sqrt{1-(\mathrm{cleaf}_2(\omega t))^2}
\\
&+ A \omega \mathrm{cleaf}_2(\omega t) \left( -\omega \mathrm{cleaf}_2(\omega t) \sqrt{1+(\mathrm{cleaf}_2(\omega t))^2}  \right) 
\\
&= A \omega^2
\left\{ 1 - 2 (\mathrm{cleaf}_2(\omega t) )^2 \right\}
\sqrt{1 + (\mathrm{cleaf}_2(\omega t))^2 } . 
\label{B3}
\end{split}
\end{equation}

(ii) In the case where $2m \pi_2 \leqq \omega t \leqq (2m+1) \pi_2 $, 
\begin{equation}
\begin{split}
\frac{\mathrm{d}^2 x(t) }{\mathrm{d} t^2}
&=   - A \omega^2 \left( -\sqrt{1-(\mathrm{cleaf}_2(\omega t))^4 } \right) \sqrt{1-(\mathrm{cleaf}_2(\omega t))^2} 
\\
&- A \omega^2 \mathrm{cleaf}_2(\omega t) \left(  \mathrm{cleaf}_2(\omega t) \sqrt{1+(\mathrm{cleaf}_2(\omega t))^2}  \right) 
\\
&= A \omega^2
\left\{ 1 - 2 (\mathrm{cleaf}_2(\omega t) )^2 \right\}
\sqrt{1 + (\mathrm{cleaf}_2(\omega t))^2 }.  
\label{B4}
\end{split}
\end{equation}

By substituting Eq. (\ref{4.1}) into Eq. (\ref{B3}) (or Eq. (\ref{B4})), Eq. (\ref{4.2}) of the Duffing equation is derived. The initial conditions are as follows:
\begin{equation}
 x(0)=A \sqrt{1+ (\mathrm{cleaf}_2(0))^2 }=A \sqrt{1+ (1)^2 }=\sqrt{2 } A, \label{B5}
\end{equation}
\begin{equation}
\frac{\mathrm{d} x(0) }{\mathrm{d} t} 
= \pm A \omega \mathrm{cleaf}_2(0) \sqrt{1-(\mathrm{cleaf}_2(0))^2 }
= \pm A \omega \sqrt{1-(1)^2 }=0.
\label{B6}
\end{equation}

\section{ Ordinary differential equation of exact solution 2 }
\renewcommand{\theequation}{C.\arabic{equation}}
\setcounter{equation}{0}

Here we prove that exact solution 2 satisfies the Duffing equation
The first derivatives of exact solution 1 are obtained as follows:

(i) In the case where $(2m-1) \pi_2 \leqq \omega t \leqq 2m \pi_2 $, 
\begin{equation}
\begin{split}
\frac{\mathrm{d} x(t) }{\mathrm{d} t} 
&= A
\frac{ -2 \mathrm{cleaf}_2(\omega t) 
\left(\omega \sqrt{1-(\mathrm{cleaf}_2(\omega t))^4 } \right) }
{ 2 \sqrt{1-(\mathrm{cleaf}_2(\omega t))^2} } \\
&= -A \omega \mathrm{cleaf}_2(\omega t) \sqrt{1+(\mathrm{cleaf}_2(\omega t))^2 }.
\label{C1}
\end{split}
\end{equation}

(ii) In the case where $2m \pi_2 \leqq \omega t \leqq (2m+1) \pi_2$, 
\begin{equation}
\begin{split}
\frac{\mathrm{d} x(t) }{\mathrm{d} t} 
&= A
\frac{ -2 \mathrm{cleaf}_2(\omega t) 
\left( -\omega\sqrt{1-(\mathrm{cleaf}_2(\omega t))^4 } \right) }
{ 2 \sqrt{1-(\mathrm{cleaf}_2(\omega t))^2} } \\
&= A \omega \mathrm{cleaf}_2(\omega t) \sqrt{1+(\mathrm{cleaf}_2(\omega t))^2 .}
\label{C2}
\end{split}
\end{equation}

The second derivatives of exact solution 2 are obtained as follows:

(i) In the case where $(2m-1) \pi_2 \leqq \omega t \leqq 2m \pi_2 $, 
\begin{equation}
\begin{split}
\frac{\mathrm{d}^2 x(t) }{\mathrm{d} t^2}
&=    -A \omega \left( \omega \sqrt{1-(\mathrm{cleaf}_2(\omega t))^4 } \right) \sqrt{1+(\mathrm{cleaf}_2(\omega t))^2} 
\\
&- A \omega \mathrm{cleaf}_2(\omega t) \left( \omega  \mathrm{cleaf}_2(\omega t) \sqrt{1-(\mathrm{cleaf}_2(\omega t))^2}  \right) 
\\
&= - A \omega^2
\left\{ 1 + 2 (\mathrm{cleaf}_2(\omega t) )^2 \right\}
\sqrt{1 - (\mathrm{cleaf}_2(\omega t))^2 } . 
\label{C3}
\end{split}
\end{equation}

(ii) In the case where $2m \pi_2 \leqq \omega t \leqq (2m+1) \pi_2$, 
\begin{equation}
\begin{split}
\frac{\mathrm{d}^2 x(t) }{\mathrm{d} t^2}
&=    A \omega^2 \left( -\sqrt{1-(\mathrm{cleaf}_2(\omega t))^4 } \right) \sqrt{1+(\mathrm{cleaf}_2(\omega t))^2}
\\
&+ A \omega^2 \mathrm{cleaf}_2(\omega t) \left( - \mathrm{cleaf}_2(\omega t) \sqrt{1-(\mathrm{cleaf}_2(\omega t))^2}  \right) 
\\
&= - A \omega^2
\left\{ 1 + 2 (\mathrm{cleaf}_2(\omega t) )^2 \right\}
\sqrt{1 - (\mathrm{cleaf}_2(\omega t))^2 }  .
\label{C4}
\end{split}
\end{equation}

By substituting Eq. (\ref{4.3}) into Eq. (\ref{C3}) (or Eq. (\ref{C4})), Eq. (\ref{4.4}) of the Duffing equation is obtained. The initial condition of the differential equation, Eq. (\ref{4.4}), is as follows:

\begin{equation}
 x(0)=A \sqrt{1- (\mathrm{cleaf}_2(0))^2 }=A \sqrt{1- (1)^2 }=0. \label{C5}
\end{equation}

Exact solution 2 is a continuous function. However,  it cannot be partially differentiated. 
Exact solution 2 is not differentiable at $t=m \frac{\pi_2}{\omega}$:
\begin{equation}
\begin{split}
\lim_{t \uparrow m\frac{\pi_2}{\omega} } \frac{\mathrm{d} x(t) }{\mathrm{d} t} 
=- \sqrt{2} A \omega,
\label{C6}
\end{split}
\end{equation}

\begin{equation}
\begin{split}
\lim_{t \downarrow m \frac{\pi_2}{\omega} } \frac{\mathrm{d} x(t) }{\mathrm{d} t} 
=  \sqrt{2} A \omega.
\label{C7}
\end{split}
\end{equation}

\section{ Ordinary differential equation of exact solution 3 }
\renewcommand{\theequation}{D.\arabic{equation}}
\setcounter{equation}{0}
Here we prove that exact solution 3 satisfies the Duffing equation. The first derivatives of exact solution 3 are obtained as follows:

(i) In the case where $\frac{4m-1}{2} \pi_2 \leqq \omega t \leqq \frac{4m+1}{2} \pi_2$,
\begin{equation}
\begin{split}
\frac{\mathrm{d} x(t) }{\mathrm{d} t} 
&= A
\frac{ 2\mathrm{sleaf}_2(\omega t) 
\left( \omega \sqrt{1-(\mathrm{sleaf}_2(\omega t))^4 } \right) }
{ 2 \sqrt{1+(\mathrm{sleaf}_2(\omega t))^2} } \\
&= A \omega \mathrm{sleaf}_2(\omega t) \sqrt{1-(\mathrm{sleaf}_2(\omega t))^2 }.
\label{D1}
\end{split}
\end{equation}

(ii) In the case where $\frac{4m+1}{2} \pi_2 \leqq \omega t \leqq \frac{4m+3}{2} \pi_2$,
\begin{equation}
\begin{split}
\frac{\mathrm{d} x(t) }{\mathrm{d} t} 
&= A
\frac{ 2\mathrm{sleaf}_2(\omega t) 
\left( -\omega \sqrt{1-(\mathrm{sleaf}_2(\omega t))^4 } \right) }
{ 2 \sqrt{1+(\mathrm{sleaf}_2(\omega t))^2} } \\
&= - A \omega \mathrm{sleaf}_2(\omega t) \sqrt{1-(\mathrm{sleaf}_2(\omega t))^2 }.
\label{D2}
\end{split}
\end{equation}

The second derivatives of exact solution 3 are obtained as follows:

(i) In the case where $\frac{4m-1}{2} \pi_2 \leqq \omega t \leqq \frac{4m+1}{2} \pi_2$,
\begin{equation}
\begin{split}
\frac{\mathrm{d}^2 x(t) }{\mathrm{d} t^2}
&=    A \omega^2 \left( \sqrt{1-(\mathrm{sleaf}_2(\omega t))^4 } \right) \sqrt{1-(\mathrm{sleaf}_2(\omega t))^2}
\\
& +A \omega^2 \mathrm{sleaf}_2(\omega t) \left( - \mathrm{sleaf}_2(\omega t) \sqrt{1+(\mathrm{sleaf}_2(\omega t))^2}  \right) 
\\
&= A \omega^2
\left\{ 1-2 (\mathrm{sleaf}_2(\omega t) )^2 \right\}
\sqrt{1 + (\mathrm{sleaf}_2(\omega t))^2 } . 
\label{D3}
\end{split}
\end{equation}

(ii) In the case where $\frac{4m+1}{2} \pi_2 \leqq \omega t \leqq \frac{4m+3}{2} \pi_2$,
\begin{equation}
\begin{split}
\frac{\mathrm{d}^2 x(t) }{\mathrm{d} t^2}
&=    - A \omega \left( -\omega \sqrt{1-(\mathrm{sleaf}_2(\omega t))^4 } \right) \sqrt{1-(\mathrm{sleaf}_2(\omega t))^2} 
\\
& -A \omega \mathrm{sleaf}_2(\omega t) \left( \omega \mathrm{sleaf}_2(\omega t) \sqrt{1+(\mathrm{sleaf}_2(\omega t))^2}  \right) 
\\
&= A \omega^2
\left\{ 1-2 (\mathrm{sleaf}_2(\omega t) )^2 \right\}
\sqrt{1 + (\mathrm{sleaf}_2(\omega t))^2 }  .
\label{D4}
\end{split}
\end{equation}

By substituting Eq. (\ref{4.5}) into Eq. (\ref{D3}) (or Eq. (\ref{D4})), Eq. (\ref{4.6}) of the Duffing equation is obtained. The initial conditions of the differential equation, Eq. (\ref{4.6}), are as follows:
\begin{equation}
 x(0)=A \sqrt{1+ (\mathrm{sleaf}_2(0))^2 }=A \sqrt{1+ (0)^2 }= A \label{D5}, 
\end{equation}

\begin{equation}
\frac{\mathrm{d} x(0) }{\mathrm{d} t} 
= \pm A \omega \mathrm{sleaf}_2(0) \sqrt{1-(\mathrm{sleaf}_2(0))^2 }
=\pm A \omega \cdot 0 \cdot \sqrt{1-(0)^2 }=0.
\label{D6}
\end{equation}

\section{ Ordinary differential equation of exact solution 4 }
\renewcommand{\theequation}{E.\arabic{equation}}
\setcounter{equation}{0}
Here we prove that exact solution 4 satisfies the Duffing equation
The first derivatives of exact solution 4 are obtained as follows:

(i) In the case where $\frac{4m-1}{2} \pi_2 \leqq \omega t \leqq \frac{4m+1}{2} \pi_2$,
\begin{equation}
\begin{split}
\frac{\mathrm{d} x(t) }{\mathrm{d} t} 
&= A
\frac{ 2 \mathrm{sleaf}_2(\omega t) 
\left(\omega \sqrt{1-(\mathrm{sleaf}_2(\omega t))^4 } \right) }
{ 2 \sqrt{1-(\mathrm{sleaf}_2(\omega t))^2} } \\
&= A \omega \mathrm{sleaf}_2(\omega t) \sqrt{1+(\mathrm{sleaf}_2(\omega t))^2 }.
\label{E1}
\end{split}
\end{equation}

(ii) In the case where $\frac{4m+1}{2} \pi_2 \leqq \omega t \leqq \frac{4m+3}{2} \pi_2 $,
\begin{equation}
\begin{split}
\frac{\mathrm{d} x(t) }{\mathrm{d} t} 
&= A
\frac{ 2 \mathrm{sleaf}_2(\omega t) 
\left( -\omega\sqrt{1-(\mathrm{sleaf}_2(\omega t))^4 } \right) }
{ 2 \sqrt{1-(\mathrm{sleaf}_2(\omega t))^2} } \\
&= - A \omega \mathrm{sleaf}_2(\omega t) \sqrt{1+(\mathrm{sleaf}_2(\omega t))^2 }.
\label{E2}
\end{split}
\end{equation}

The second derivatives of exact solution 4 are obtained as follows:

(i) In the case where $\frac{4m-1}{2} \pi_2 \leqq \omega t \leqq \frac{4m+1}{2} \pi_2$,
\begin{equation}
\begin{split}
\frac{\mathrm{d}^2 x(t) }{\mathrm{d} t^2}
&=    A \omega^2 \left( \sqrt{1-(\mathrm{sleaf}_2(\omega t))^4 } \right) \sqrt{1+(\mathrm{sleaf}_2(\omega t))^2}
\\
& +A \omega^2 \mathrm{sleaf}_2(\omega t) \left(  \mathrm{sleaf}_2(\omega t) \sqrt{1-(\mathrm{sleaf}_2(\omega t))^2}  \right) 
\\
&= A \omega^2
\left\{ 1 + 2 (\mathrm{sleaf}_2(\omega t) )^2 \right\}
\sqrt{1 - (\mathrm{sleaf}_2(\omega t))^2 }  .
\label{E3}
\end{split}
\end{equation}

(ii) In the case where $\frac{4m+1}{2} \pi_2 \leqq \omega t \leqq \frac{4m+3}{2} \pi_2 $,
\begin{equation}
\begin{split}
\frac{\mathrm{d}^2 x(t) }{\mathrm{d} t^2}
&=    -A \omega \left( -\omega \sqrt{1-(\mathrm{sleaf}_2(\omega t))^4 } \right) \sqrt{1+(\mathrm{sleaf}_2(\omega t))^2} 
\\
& -A \omega \mathrm{sleaf}_2(\omega t) \left(-\omega  \mathrm{sleaf}_2(\omega t) \sqrt{1-(\mathrm{sleaf}_2(\omega t))^2}  \right) 
\\
&= A \omega^2
\left\{ 1 + 2 (\mathrm{sleaf}_2(\omega t) )^2 \right\}
\sqrt{1-(\mathrm{sleaf}_2(\omega t))^2 }  .
\label{E4}
\end{split}
\end{equation}

By substituting Eq. (\ref{4.6}) into Eq. (\ref{E3}) (or Eq. (\ref{E4})), Eq. (\ref{4.4}) of the Duffing equation is obtained. The initial condition of the differential equation is as follows:
\begin{equation}
 x(0)=A \sqrt{1- (\mathrm{sleaf}_2(0))^2 }=A \sqrt{1- (0)^2 }=A .\label{E5}
\end{equation}

Exact solution 4 is a continuous function. However,  it cannot be partially differentiated. 
Exact solution 4 is not differentiable at $t=\frac{2m-1}{2} \frac{\pi_2}{\omega}$:

\begin{equation}
\begin{split}
\lim_{t \uparrow \frac{2m-1}{2} \frac{\pi_2}{\omega} } \frac{\mathrm{d} x(t) }{\mathrm{d} t} = - \sqrt{2} A \omega,
\label{E6}
\end{split}
\end{equation}

\begin{equation}
\begin{split}
\lim_{t \downarrow \frac{2m-1}{2} \frac{\pi_2}{\omega} } \frac{\mathrm{d} x(t) }{\mathrm{d} t} = \sqrt{2} A \omega.
\label{E7}
\end{split}
\end{equation}

\section{ Ordinary differential equation of exact solution 5 }
\renewcommand{\theequation}{F.\arabic{equation}}
\setcounter{equation}{0}

Here we prove that exact solution 5 satisfies the Duffing equation.
The first derivative of exact solution 5 is obtained  by Eq. (\ref{B1}) (or Eq. (\ref{B2})) and Eq. (\ref{D1}) (or Eq. (\ref{D2})).

(i) In the case where $\frac{4m-1}{2} \pi_2 \leqq \omega t \leqq \frac{4m}{2} \pi_2$,
\begin{equation}
\begin{split}
\frac{\mathrm{d} x(t) }{\mathrm{d} t} 
= A \omega \mathrm{sleaf}_2(\omega t) 
\sqrt{1-(\mathrm{sleaf}_2(\omega t))^2 }
+A \omega \mathrm{cleaf}_2(\omega t)
 \sqrt{1-(\mathrm{cleaf}_2(\omega t))^2}.
\label{F1}
\end{split}
\end{equation}

(ii) In the case where $ \frac{4m}{2}  \pi_2 \leqq \omega t \leqq  \frac{4m+1}{2} \pi_2$,
\begin{equation}
\begin{split}
\frac{\mathrm{d} x(t) }{\mathrm{d} t} 
= A \omega \mathrm{sleaf}_2(\omega t) 
\sqrt{1-(\mathrm{sleaf}_2(\omega t))^2 }
-A \omega \mathrm{cleaf}_2(\omega t)
 \sqrt{1-(\mathrm{cleaf}_2(\omega t))^2}.
\label{F2}
\end{split}
\end{equation}

(iii) In the case where $ \frac{4m+1}{2} \pi_2 \leqq \omega t \leqq  \frac{4m+2}{2} \pi_2$,
\begin{equation}
\begin{split}
\frac{\mathrm{d} x(t) }{\mathrm{d} t} 
= -A \omega \mathrm{sleaf}_2(\omega t) 
\sqrt{1-(\mathrm{sleaf}_2(\omega t))^2 }
-A \omega \mathrm{cleaf}_2(\omega t)
 \sqrt{1-(\mathrm{cleaf}_2(\omega t))^2}.
\label{F3}
\end{split}
\end{equation}

(iv) In the case where $ \frac{4m+2}{2} \pi_2 \leqq \omega t \leqq  \frac{4m+3}{2} \pi_2$,
\begin{equation}
\begin{split}
\frac{\mathrm{d} x(t) }{\mathrm{d} t} 
= -A \omega \mathrm{sleaf}_2(\omega t) 
\sqrt{1-(\mathrm{sleaf}_2(\omega t))^2 }
+A \omega \mathrm{cleaf}_2(\omega t)
 \sqrt{1-(\mathrm{cleaf}_2(\omega t))^2}.
\label{F4}
\end{split}
\end{equation}

With respect to all domains in (i)-(iv), the second derivative of exact solution 5 is obtained as follows:
\begin{equation}
\begin{split}
\frac{\mathrm{d}^2 x(t) }{\mathrm{d} t^2}
&=   A \omega^2
\left\{ 1 - 2 (\mathrm{sleaf}_2(\omega t) )^2 \right\}
\sqrt{1 + (\mathrm{sleaf}_2(\omega t))^2 }  \\
&+A \omega^2
\left\{ 1 - 2 (\mathrm{cleaf}_2(\omega t) )^2 \right\}
\sqrt{1 + (\mathrm{cleaf}_2(\omega t))^2 } . 
\label{F5}
\end{split}
\end{equation}

By using Eq. (\ref{S1}), the relational expression is obtained by the following transformation of the equation: 
\begin{equation}
\begin{split}
& \left\{ 1 + (\mathrm{sleaf}_2(\omega t) )^2 \right\}
\sqrt{1 + (\mathrm{cleaf}_2(\omega t))^2 }  \\
&=\frac{1+(\mathrm{cleaf}_2(\omega t))^2+
1-(\mathrm{cleaf}_2(\omega t))^2}
{1+(\mathrm{cleaf}_2(\omega t))^2}
\sqrt{1 + (\mathrm{cleaf}_2(\omega t))^2 } \\
&=\frac{2}{ \sqrt{1 + (\mathrm{cleaf}_2(\omega t))^2 } }
=\sqrt{2}\sqrt{1 + (\mathrm{sleaf}_2(\omega t))^2}.
\label{F6}
\end{split}
\end{equation}

Similarly, the following relational expression is obtained:
\begin{equation}
\begin{split}
& \left\{ 1 + (\mathrm{cleaf}_2(\omega t) )^2 \right\}
\sqrt{1 + (\mathrm{sleaf}_2(\omega t))^2 }  \\
&=\frac{1+(\mathrm{sleaf}_2(\omega t))^2+
1-(\mathrm{sleaf}_2(\omega t))^2}
{1+(\mathrm{sleaf}_2(\omega t))^2}
\sqrt{1 + (\mathrm{sleaf}_2(\omega t))^2 } \\
&=\frac{2}{ \sqrt{1 + (\mathrm{sleaf}_2(\omega t))^2 } }
=\sqrt{2}\sqrt{1 + (\mathrm{cleaf}_2(\omega t))^2}.
\label{F7}
\end{split}
\end{equation}

Next, the following equation is derived:
\begin{equation}
\begin{split}
&(x(t))^3 =A^3 
\left\{ 1 + (\mathrm{sleaf}_2(\omega t) )^2 \right\} 
\sqrt{1 + (\mathrm{sleaf}_2(\omega t))^2 }  \\
&+3A^3 
\left\{ 1 + (\mathrm{cleaf}_2(\omega t) )^2 \right\} 
\sqrt{1 + (\mathrm{sleaf}_2(\omega t))^2 }  \\
&+3A^3 
\left\{ 1 + (\mathrm{sleaf}_2(\omega t) )^2 \right\} 
\sqrt{1 + (\mathrm{cleaf}_2(\omega t))^2 }  \\
&+A^3 
\left\{ 1 + (\mathrm{cleaf}_2(\omega t) )^2 \right\} 
\sqrt{1 + (\mathrm{cleaf}_2(\omega t))^2 }  \\
&=A^3 
\left\{ 1 + (\mathrm{sleaf}_2(\omega t) )^2 \right\} 
\sqrt{1 + (\mathrm{sleaf}_2(\omega t))^2 }  \\
&+3A^3 \sqrt{2}
\sqrt{1 + (\mathrm{cleaf}_2(\omega t))^2 }  
+3A^3 \sqrt{2}
\sqrt{1 + (\mathrm{sleaf}_2(\omega t))^2 }  \\
&+A^3 
\left\{ 1 + (\mathrm{cleaf}_2(\omega t) )^2 \right\} 
\sqrt{1 + (\mathrm{cleaf}_2(\omega t))^2 }  \\
&=A^3  (\mathrm{sleaf}_2(\omega t) )^2 
\sqrt{1 + (\mathrm{sleaf}_2(\omega t))^2 } 
+A^3  (\mathrm{cleaf}_2(\omega t) )^2 
\sqrt{1 + (\mathrm{cleaf}_2(\omega t))^2 } \\
&+A^3 (3 \sqrt{2}+1)
\sqrt{1 + (\mathrm{cleaf}_2(\omega t))^2 }  
+A^3 (3 \sqrt{2}+1)
\sqrt{1 + (\mathrm{sleaf}_2(\omega t))^2 } . 
\label{F8}
\end{split}
\end{equation}

The above  equation is transformed as follows: 
\begin{equation}
\begin{split}
&(\mathrm{sleaf}_2(\omega t) )^2 
\sqrt{1 + (\mathrm{sleaf}_2(\omega t))^2 }
+
(\mathrm{cleaf}_2(\omega t) )^2
\sqrt{1 + (\mathrm{cleaf}_2(\omega t))^2 } \\
&=\frac{(x(t))^3}{A^3}-(3 \sqrt{2}+1) \frac{x(t)}{A}.
\label{F9}
\end{split}
\end{equation}

Eq. (\ref{4.8}) of the Duffing equation is obtained. The initial conditions are as follows:
\begin{equation}
 x(0)=A \sqrt{1+ (\mathrm{sleaf}_2(0))^2 }
+A \sqrt{1+ (\mathrm{cleaf}_2(0))^2 }
=(1+\sqrt{2})A, \label{F10}
\end{equation}

\begin{equation}
\begin{split}
\frac{\mathrm{d} x(0) }{\mathrm{d} t}
=   0.
\label{F11}
\end{split}
\end{equation}

\section{ Ordinary differential equation of exact solution 6 }
\renewcommand{\theequation}{G.\arabic{equation}}
\setcounter{equation}{0}

Here we prove that exact solution 6 satisfies the Duffing equation. The first derivative of exact solution 6 is obtained  by Eq. (\ref{B1}) (or Eq. (\ref{B2})) and Eq. (\ref{D1}) (or Eq. (\ref{D2})).

(i) In the case where $\frac{4m-1}{2} \pi_2 \leqq \omega t \leqq \frac{4m}{2} \pi_2$,
\begin{equation}
\begin{split}
\frac{\mathrm{d} x(t) }{\mathrm{d} t} 
= - A \omega \mathrm{sleaf}_2(\omega t) 
\sqrt{1-(\mathrm{sleaf}_2(\omega t))^2 }
+A \omega \mathrm{cleaf}_2(\omega t)
 \sqrt{1-(\mathrm{cleaf}_2(\omega t))^2}.
\label{G1}
\end{split}
\end{equation}

(ii) In the case where $ \frac{4m}{2}  \pi_2 \leqq \omega t \leqq  \frac{4m+1}{2} \pi_2$,
\begin{equation}
\begin{split}
\frac{\mathrm{d} x(t) }{\mathrm{d} t} 
= -A \omega \mathrm{sleaf}_2(\omega t) 
\sqrt{1-(\mathrm{sleaf}_2(\omega t))^2 }
-A \omega \mathrm{cleaf}_2(\omega t)
 \sqrt{1-(\mathrm{cleaf}_2(\omega t))^2}.
\label{G2}
\end{split}
\end{equation}

(iii) In the case where $ \frac{4m+1}{2} \pi_2 \leqq \omega t \leqq  \frac{4m+2}{2} \pi_2$,
\begin{equation}
\begin{split}
\frac{\mathrm{d} x(t) }{\mathrm{d} t} 
= A \omega \mathrm{sleaf}_2(\omega t) 
\sqrt{1-(\mathrm{sleaf}_2(\omega t))^2 }
-A \omega \mathrm{cleaf}_2(\omega t)
 \sqrt{1-(\mathrm{cleaf}_2(\omega t))^2}.
\label{G3}
\end{split}
\end{equation}

(iv) In the case where $ \frac{4m+2}{2} \pi_2 \leqq \omega t \leqq  \frac{4m+3}{2} \pi_2$,
\begin{equation}
\begin{split}
\frac{\mathrm{d} x(t) }{\mathrm{d} t} 
= A \omega \mathrm{sleaf}_2(\omega t) 
\sqrt{1-(\mathrm{sleaf}_2(\omega t))^2 }
+A \omega \mathrm{cleaf}_2(\omega t)
 \sqrt{1-(\mathrm{cleaf}_2(\omega t))^2}.
\label{G4}
\end{split}
\end{equation}

With respect to all domains in (i)-(iv), the second derivative of exact solution 6 is obtained as follows:
\begin{equation}
\begin{split}
\frac{\mathrm{d}^2 x(t) }{\mathrm{d} t^2}
&=  - A \omega^2
\left\{ 1 - 2 (\mathrm{sleaf}_2(\omega t) )^2 \right\}
\sqrt{1 + (\mathrm{sleaf}_2(\omega t))^2 }  \\
&+A \omega^2
\left\{ 1 - 2 (\mathrm{cleaf}_2(\omega t) )^2 \right\}
\sqrt{1 + (\mathrm{cleaf}_2(\omega t))^2 }  .
\label{G5}
\end{split}
\end{equation}

Next, the following equation is derived:
\begin{equation}
\begin{split}
&(x(t))^3 =-A^3 
\left\{ 1 + (\mathrm{sleaf}_2(\omega t) )^2 \right\} 
\sqrt{1 + (\mathrm{sleaf}_2(\omega t))^2 }  \\
&-3A^3 
\left\{ 1 + (\mathrm{cleaf}_2(\omega t) )^2 \right\} 
\sqrt{1 + (\mathrm{sleaf}_2(\omega t))^2 }  \\
&+3A^3 
\left\{ 1 + (\mathrm{sleaf}_2(\omega t) )^2 \right\} 
\sqrt{1 + (\mathrm{cleaf}_2(\omega t))^2 }  \\
&+A^3 
\left\{ 1 + (\mathrm{cleaf}_2(\omega t) )^2 \right\} 
\sqrt{1 + (\mathrm{cleaf}_2(\omega t))^2 }  \\
&=-A^3 
\left\{ 1 + (\mathrm{sleaf}_2(\omega t) )^2 \right\} 
\sqrt{1 + (\mathrm{sleaf}_2(\omega t))^2 }  \\
&-3A^3 \sqrt{2}
\sqrt{1 + (\mathrm{cleaf}_2(\omega t))^2 }  
+3A^3 \sqrt{2}
\sqrt{1 + (\mathrm{sleaf}_2(\omega t))^2 }  \\
&+A^3 
\left\{ 1 + (\mathrm{cleaf}_2(\omega t) )^2 \right\} 
\sqrt{1 + (\mathrm{cleaf}_2(\omega t))^2 }  \\
&=-A^3  (\mathrm{sleaf}_2(\omega t) )^2 
\sqrt{1 + (\mathrm{sleaf}_2(\omega t))^2 } 
+A^3  (\mathrm{cleaf}_2(\omega t) )^2 
\sqrt{1 + (\mathrm{cleaf}_2(\omega t))^2 } \\
&-A^3 (3 \sqrt{2}-1)
\sqrt{1 + (\mathrm{cleaf}_2(\omega t))^2 }  
+A^3 (3 \sqrt{2}-1)
\sqrt{1 + (\mathrm{sleaf}_2(\omega t))^2 }.  
\label{G6}
\end{split}
\end{equation}

The above  equation is transformed as follows: 
\begin{equation}
\begin{split}
-&(\mathrm{sleaf}_2(\omega t) )^2 
\sqrt{1 + (\mathrm{sleaf}_2(\omega t))^2 }
+
(\mathrm{cleaf}_2(\omega t) )^2
\sqrt{1 + (\mathrm{cleaf}_2(\omega t))^2 } \\
&=\frac{(x(t))^3}{A^3} + (3 \sqrt{2}-1) \frac{x(t)}{A}.
\label{G7}
\end{split}
\end{equation}

Eq. (\ref{4.12}) of the Duffing equation is obtained. The initial conditions of Eq. (\ref{4.12}) are as follows:
\begin{equation}
 x(0)=-A \sqrt{1+ (\mathrm{sleaf}_2(0))^2 }
+A \sqrt{1+ (\mathrm{cleaf}_2(0))^2 }
=(\sqrt{2}-1)A, \label{G8}
\end{equation}

\begin{equation}
\begin{split}
\frac{\mathrm{d} x(0) }{\mathrm{d} t}
=0.
\label{G9}
\end{split}
\end{equation}

\section{ Ordinary differential equation of exact solution 7 }
\renewcommand{\theequation}{H.\arabic{equation}}
\setcounter{equation}{0}

Here we prove that exact solution 7 satisfies the Duffing equation
The first derivatives of exact solution 7 are obtained as follows:

(i) In the case where $4m \eta_2 \leqq \omega t < (4m+1) \eta_2$ and $(4m+1) \eta_2 < \omega t \leqq (4m+2) \eta_2$,
\begin{equation}
\begin{split}
\frac{\mathrm{d} x(t) }{\mathrm{d} t} 
&= A
\frac{ 2 \mathrm{cleafh}_2(\omega t) 
\left(\omega \sqrt{(\mathrm{cleafh}_2(\omega t))^4-1 } \right) }
{ 2 \sqrt{(\mathrm{cleafh}_2(\omega t))^2+1} } \\
&= A \omega \mathrm{cleafh}_2(\omega t) \sqrt{(\mathrm{cleafh}_2(\omega t))^2-1 }.
\label{H1}
\end{split}
\end{equation}

(ii) In the case where $(4m-2) \eta_2 \leqq \omega t < (4m-1) \eta_2$ and $(4m-1) \eta_2 < \omega t \leqq 4m \eta_2$, 
\begin{equation}
\begin{split}
\frac{\mathrm{d} x(t) }{\mathrm{d} t} 
&= A
\frac{ 2 \mathrm{cleafh}_2(\omega t) 
\left(-\omega \sqrt{(\mathrm{cleafh}_2(\omega t))^4-1 } \right) }
{ 2 \sqrt{(\mathrm{cleafh}_2(\omega t))^2+1} } \\
&=- A \omega \mathrm{cleafh}_2(\omega t) \sqrt{(\mathrm{cleafh}_2(\omega t))^2-1 }.
\label{H2}
\end{split}
\end{equation}

With respect to all domains in (i) and (ii), 
the second derivative of exact solution 7 is obtained as follows:
\begin{equation}
\begin{split}
\frac{\mathrm{d}^2 x(t) }{\mathrm{d} t^2}
&= A \omega^2
\left\{  2 (\mathrm{cleafh}_2(\omega t) )^2-1 \right\}
\sqrt{ (\mathrm{cleafh}_2(\omega t))^2+1 }  
\label{H3}.
\end{split}
\end{equation}

By substituting Eq. (\ref{4.11}) into Eq. (\ref{H3}), Eq. (\ref{4.11}) of the Duffing equation is obtained. The initial conditions are as follows:

\begin{equation}
 x(0)=A \sqrt{ (\mathrm{cleafh}_2(0))^2+1 }
=\sqrt{2}A \label{H4},
\end{equation}

\begin{equation}
\begin{split}
\frac{\mathrm{d} x(0) }{\mathrm{d} t}=0
\label{H5}.
\end{split}
\end{equation}

\section{ Ordinary differential equation of exact solution 8 }
\renewcommand{\theequation}{I.\arabic{equation}}
\setcounter{equation}{0}

Here we prove that exact solution 8 satisfies the Duffing equation
The first derivatives of  exact solution 8 are obtained as follows:

(i) In the case where $4m \eta_2 \leqq \omega t < (4m+1) \eta_2$ and $(4m+1) \eta_2 < \omega t \leqq (4m+2) \eta_2$,
\begin{equation}
\begin{split}
\frac{\mathrm{d} x(t) }{\mathrm{d} t} 
&= A
\frac{ 2\mathrm{cleafh}_2(\omega t) 
\left(\omega  \sqrt{(\mathrm{cleafh}_2(\omega t))^4-1 } \right) }
{ 2 \sqrt{(\mathrm{cleafh}_2(\omega t))^2-1} } \\
&= A \omega \mathrm{cleafh}_2(\omega t) \sqrt{(\mathrm{cleafh}_2(\omega t))^2+1 }.
\label{I1}
\end{split}
\end{equation}

(ii) In the case where $(4m-2) \eta_2 \leqq \omega t < (4m-1) \eta_2$ and $(4m-1) \eta_2 < \omega t \leqq 4m \eta_2$, 
\begin{equation}
\begin{split}
\frac{\mathrm{d} x(t) }{\mathrm{d} t} 
&= A
\frac{ 2 \mathrm{cleafh}_2(\omega t) 
\left( -\omega \sqrt{(\mathrm{cleafh}_2(\omega t))^4-1 } \right) }
{ 2 \sqrt{(\mathrm{cleafh}_2(\omega t))^2-1} } \\
&= -A \omega \mathrm{cleafh}_2(\omega t) \sqrt{(\mathrm{cleafh}_2(\omega t))^2+1 }.
\label{I2}
\end{split}
\end{equation}

With respect to all domains in (i) and (ii), 
the second derivative of the exact solution 8 is obtained as follows:
\begin{equation}
\begin{split}
\frac{\mathrm{d}^2 x(t) }{\mathrm{d} t^2}
= A \omega^2
\left\{  2 (\mathrm{cleafh}_2(\omega t) )^2+1 \right\}
\sqrt{ (\mathrm{cleafh}_2(\omega t))^2-1 }  
\label{I3}.
\end{split}
\end{equation}

By substituting Eq. (\ref{4.12}) into Eq. (\ref{I3}), Eq. (\ref{4.13}) of the Duffing equation is obtained.
The initial condition is as follows:
\begin{equation}
 x(0)=A \sqrt{ (\mathrm{cleafh}_2(0))^2-1 }
=0 \label{I4}.
\end{equation}

Exact solution 8 is a continuous function. However,  it cannot be partially differentiated. Exact solution 8 is not differentiable at $t=\frac{2m\eta_2}{\omega}$:
\begin{equation}
\begin{split}
\lim_{t \uparrow  \frac{2m\eta_2}{\omega} } \frac{\mathrm{d} x(t) }{\mathrm{d} t} 
=  - \sqrt{2} A \omega,
\label{I5}
\end{split}
\end{equation}

\begin{equation}
\begin{split}
\lim_{t \downarrow  \frac{2m\eta_2}{\omega} } \frac{\mathrm{d} x(t) }{\mathrm{d} t} 
=  \sqrt{2} A \omega.
\label{I6}
\end{split}
\end{equation}

\section{ Ordinary differential equation of exact solution 9 }
\renewcommand{\theequation}{J.\arabic{equation}}
\setcounter{equation}{0}

Here we prove that exact solution 9 satisfies the Duffing equation. The first derivatives of exact solution 9 are obtained as follows: 

(i) In the case where $\frac{4m-1}{2} \pi_2 \leqq \omega t \leqq \frac{4m+1}{2} \pi_2$,
\begin{equation}
\begin{split}
\frac{\mathrm{d} x(t) }{\mathrm{d} t} 
&= A \omega
\frac{  \sqrt{1-(\mathrm{sleaf}_2(\omega t))^4 }  }
{ 2 \sqrt{1+ \mathrm{sleaf}_2(\omega t)} } \\
&= \frac{A \omega}{2} \sqrt{1- \mathrm{sleaf}_2(\omega t)
+(\mathrm{sleaf}_2(\omega t))^2
-(\mathrm{sleaf}_2(\omega t))^3  }
\label{J1}.
\end{split}
\end{equation}

(ii) In the case where $\frac{4m+1}{2} \pi_2 \leqq \omega t \leqq \frac{4m+3}{2} \pi_2$,
\begin{equation}
\begin{split}
\frac{\mathrm{d} x(t) }{\mathrm{d} t} 
&= A \omega
\frac{  -\sqrt{1-(\mathrm{sleaf}_2(\omega t))^4 }  }
{ 2 \sqrt{1+ \mathrm{sleaf}_2(\omega t)} } \\
&= -\frac{A \omega}{2} \sqrt{1- \mathrm{sleaf}_2(\omega t)
+(\mathrm{sleaf}_2(\omega t))^2
-(\mathrm{sleaf}_2(\omega t))^3  }.
\label{J2}
\end{split}
\end{equation}

The second derivatives of exact solution 9 are obtained as follows:

(i) In the case where $\frac{4m-1}{2} \pi_2 \leqq \omega t \leqq \frac{4m+1}{2} \pi_2$,
\begin{equation}
\begin{split}
\frac{\mathrm{d}^2 x(t) }{\mathrm{d} t^2}
&=    \frac{A \omega^2}{4} 
\frac{  \{ -1 +2 \mathrm{sleaf}_2(\omega t) 
-3 (\mathrm{sleaf}_2(\omega t))^2 \} 
\sqrt{1-(\mathrm{sleaf}_2(\omega t))^4 }  }
{ \sqrt{1- \mathrm{sleaf}_2(\omega t)
+(\mathrm{sleaf}_2(\omega t))^2
-(\mathrm{sleaf}_2(\omega t))^3} } \\
&=    \frac{A \omega^2}{4}  
 \{ -1 +2 \mathrm{sleaf}_2(\omega t) 
-3 (\mathrm{sleaf}_2(\omega t))^2 \} 
\sqrt{1+\mathrm{sleaf}_2(\omega t) } 
\\
&=    \frac{A \omega^2}{4}  
 \left\{ -1 +2 \left( \frac{x(t)^2}{A^2}-1 \right)  -3 \left( \frac{x(t)^2}{A^2}-1 \right)^2 \right\} \frac{x(t)}{A} \\
&=   -\frac{3}{2} \omega^2 x(t) + 2\frac{\omega^2}{A^2} x(t)^3
-\frac{3}{4} \frac{\omega^2}{A^4} x(t)^5.
\label{J3}
\end{split}
\end{equation}

(ii) In the case where $\frac{4m+1}{2} \pi_2 \leqq \omega t \leqq \frac{4m+3}{2} \pi_2$,
\begin{equation}
\begin{split}
\frac{\mathrm{d}^2 x(t) }{\mathrm{d} t^2}
&=    -\frac{A \omega^2}{4} 
\frac{  \{ -1 +2 \mathrm{sleaf}_2(\omega t) 
-3 (\mathrm{sleaf}_2(\omega t))^2 \} 
(-\sqrt{1-(\mathrm{sleaf}_2(\omega t))^4 } )  }
{ \sqrt{1- \mathrm{sleaf}_2(\omega t)
+(\mathrm{sleaf}_2(\omega t))^2
-(\mathrm{sleaf}_2(\omega t))^3} } \\
&=   -\frac{3}{2} \omega^2 x(t) + 2\frac{\omega^2}{A^2} x(t)^3
-\frac{3}{4} \frac{\omega^2}{A^4} x(t)^5
\label{J4}.
\end{split}
\end{equation}

Eq. (\ref{5.2}) of the Duffing equation is obtained. The initial condition of the differential equation, Eq. (\ref{5.2}), is as follows:
\begin{equation}
 x(0)=A \sqrt{1+ \mathrm{sleaf}_2(0) }=A \sqrt{1+ 0 }= A. \label{J5}
\end{equation}

Exact solution 9 is a continuous function. However,  it cannot be partially differentiated. Exact solution 9 is not differentiable at $t=\frac{4m-1}{2} \frac{\pi_2}{\omega}$:
\begin{equation}
\begin{split}
&\lim_{t \uparrow \frac{4m-1}{2} \frac{\pi_2}{\omega} }  \frac{\mathrm{d} x(t) }{\mathrm{d} t} 
= -A \omega,
\label{J6}
\end{split}
\end{equation}

\begin{equation}
\begin{split}
&\lim_{t \downarrow \frac{4m-1}{2} \frac{\pi_2}{\omega} }  \frac{\mathrm{d} x(t) }{\mathrm{d} t} 
= A \omega.
\label{J7}
\end{split}
\end{equation}

\section{ Ordinary differential equation of exact solution 10 }
\renewcommand{\theequation}{K.\arabic{equation}}
\setcounter{equation}{0}

Here we prove that exact solution 10 satisfies the Duffing equation. The first derivatives of exact solution 10 are obtained as follows: 

(i) In the case where $\frac{4m-1}{2} \pi_2 \leqq \omega t \leqq \frac{4m+1}{2} \pi_2$,
\begin{equation}
\begin{split}
\frac{\mathrm{d} x(t) }{\mathrm{d} t} 
&= -A \omega
\frac{  \sqrt{1-(\mathrm{sleaf}_2(\omega t))^4 }  }
{ 2 \sqrt{1- \mathrm{sleaf}_2(\omega t)} } \\
&= -\frac{A \omega}{2} \sqrt{1+ \mathrm{sleaf}_2(\omega t)
+(\mathrm{sleaf}_2(\omega t))^2
+(\mathrm{sleaf}_2(\omega t))^3  }.
\label{K1}
\end{split}
\end{equation}

(ii) In the case where $\frac{4m+1}{2} \pi_2 \leqq \omega t \leqq \frac{4m+3}{2} \pi_2$,
\begin{equation}
\begin{split}
\frac{\mathrm{d} x(t) }{\mathrm{d} t} 
&= -A \omega
\frac{  - \sqrt{1-(\mathrm{sleaf}_2(\omega t))^4 }  }
{ 2 \sqrt{1- \mathrm{sleaf}_2(\omega t)} } \\
&= \frac{A \omega}{2} \sqrt{1+ \mathrm{sleaf}_2(\omega t)
+(\mathrm{sleaf}_2(\omega t))^2
+(\mathrm{sleaf}_2(\omega t))^3  }.
\label{K2}
\end{split}
\end{equation}

The second derivatives of exact solution 10 are obtained as follows:

(i) In the case where $\frac{4m-1}{2} \pi_2 \leqq \omega t \leqq \frac{4m+1}{2} \pi_2$,
\begin{equation}
\begin{split}
\frac{\mathrm{d}^2 x(t) }{\mathrm{d} t^2}
&=    -\frac{A \omega^2}{4} 
\frac{  \{ 1 +2 \mathrm{sleaf}_2(\omega t) 
+3 (\mathrm{sleaf}_2(\omega t))^2 \} 
\sqrt{1-(\mathrm{sleaf}_2(\omega t))^4 }  }
{ \sqrt{1+ \mathrm{sleaf}_2(\omega t)
+(\mathrm{sleaf}_2(\omega t))^2
+(\mathrm{sleaf}_2(\omega t))^3} } \\
&=    -\frac{A \omega^2}{4}  
 \{ 1 +2 \mathrm{sleaf}_2(\omega t) 
+3 (\mathrm{sleaf}_2(\omega t))^2 \} 
\sqrt{1-\mathrm{sleaf}_2(\omega t) } 
\\
&=    -\frac{A \omega^2}{4}  
 \left\{ 1 +2 \left( 1 - \frac{x(t)^2}{A^2} \right)  +3 \left( 1-\frac{x(t)^2}{A^2} \right)^2 \right\} \frac{x(t)}{A} \\
&=  - \frac{3}{2} \omega^2 x(t) + 2\frac{\omega^2}{A^2} x(t)^3
-\frac{3}{4} \frac{\omega^2}{A^4} x(t)^5.
\label{K3}
\end{split}
\end{equation}

(ii) In the case where $\frac{4m+1}{2} \pi_2 \leqq \omega t \leqq \frac{4m+3}{2} \pi_2$,
\begin{equation}
\begin{split}
\frac{\mathrm{d}^2 x(t) }{\mathrm{d} t^2}
&=    \frac{A \omega^2}{4} 
\frac{  \{ 1 +2 \mathrm{sleaf}_2(\omega t) 
+3 (\mathrm{sleaf}_2(\omega t))^2 \} 
(-\sqrt{1-(\mathrm{sleaf}_2(\omega t))^4 } ) }
{ \sqrt{1+ \mathrm{sleaf}_2(\omega t)
+(\mathrm{sleaf}_2(\omega t))^2
+(\mathrm{sleaf}_2(\omega t))^3} } \\
&=  - \frac{3}{2} \omega^2 x(t) + 2\frac{\omega^2}{A^2} x(t)^3
-\frac{3}{4} \frac{\omega^2}{A^4} x(t)^5.
\label{K4}
\end{split}
\end{equation}

Eq. (\ref{5.2}) of the Duffing equation is obtained. The initial condition of the differential equation, Eq. (\ref{5.2}), is as follows:
\begin{equation}
 x(0)=A \sqrt{1- \mathrm{sleaf}_2(0) }=A \sqrt{1- 0 }= A. \label{K5}
\end{equation}

Exact solution 10 is a continuous function. However,  it cannot be partially differentiated. Exact solution 10 is not differentiable at $t=\frac{4m+1}{2} \frac{\pi_2}{\omega} $: 
\begin{equation}
\begin{split}
&\lim_{t \uparrow \frac{4m+1}{2}  \frac{\pi_2}{\omega} }  \frac{\mathrm{d} x(t) }{\mathrm{d} t} 
=  -A \omega,
\label{K6}
\end{split}
\end{equation}

\begin{equation}
\begin{split}
&\lim_{t \downarrow \frac{4m+1}{2} \frac{\pi_2}{\omega} }  \frac{\mathrm{d} x(t) }{\mathrm{d} t} 
=  A \omega.
\label{K7}
\end{split}
\end{equation}

\section{ Ordinary differential equation of exact solution 11 }
\renewcommand{\theequation}{L.\arabic{equation}}
\setcounter{equation}{0}

Here we prove that exact solution 11 satisfies the Duffing equation. The first derivatives of exact solution 11 are obtained as follows: 

(i) In the case where $2m \pi_2 \leqq \omega t \leqq (2m +1) \pi_2$,
\begin{equation}
\begin{split}
\frac{\mathrm{d} x(t) }{\mathrm{d} t} 
&= A \omega
\frac{  -\sqrt{1-(\mathrm{cleaf}_2(\omega t))^4 }  }
{ 2 \sqrt{1+ \mathrm{cleaf}_2(\omega t)} } \\
&= -\frac{A \omega}{2} \sqrt{1- \mathrm{cleaf}_2(\omega t)
+(\mathrm{cleaf}_2(\omega t))^2
-(\mathrm{cleaf}_2(\omega t))^3  }.
\label{L1}
\end{split}
\end{equation}

(ii) In the case where $(2m+1) \pi_2 \leqq \omega t \leqq (2m+2) \pi_2$, 
\begin{equation}
\begin{split}
\frac{\mathrm{d} x(t) }{\mathrm{d} t} 
&= A \omega
\frac{  \sqrt{1-(\mathrm{cleaf}_2(\omega t))^4 }  }
{ 2 \sqrt{1+ \mathrm{cleaf}_2(\omega t)} } \\
&= \frac{A \omega}{2} \sqrt{1- \mathrm{cleaf}_2(\omega t)
+(\mathrm{cleaf}_2(\omega t))^2
-(\mathrm{cleaf}_2(\omega t))^3  }.
\label{L2}
\end{split}
\end{equation}

The second derivatives of exact solution 11 are obtained as follows:

(i) In the case where $2m \pi_2 \leqq \omega t \leqq (2m +1) \pi_2$,
\begin{equation}
\begin{split}
\frac{\mathrm{d}^2 x(t) }{\mathrm{d} t^2}
&=    -\frac{A \omega}{2} 
\frac{  \{ -1 +2 \mathrm{cleaf}_2(\omega t) 
-3 (\mathrm{cleaf}_2(\omega t))^2 \} 
( -\omega \sqrt{1-(\mathrm{cleaf}_2(\omega t))^4 } )  }
{ 2 \sqrt{1- \mathrm{cleaf}_2(\omega t)
+(\mathrm{cleaf}_2(\omega t))^2
-(\mathrm{cleaf}_2(\omega t))^3} } \\
&=    \frac{A \omega^2}{4}  
 \{ -1 +2 \mathrm{cleaf}_2(\omega t) 
-3 (\mathrm{cleaf}_2(\omega t))^2 \} 
\sqrt{1+\mathrm{cleaf}_2(\omega t) } 
\\
&=    \frac{A \omega^2}{4}  
 \left\{ -1 +2 \left( \frac{x(t)^2}{A^2}-1 \right)  -3 \left( \frac{x(t)^2}{A^2}-1 \right)^2 \right\} \frac{x(t)}{A} \\
&=   -\frac{3}{2} \omega^2 x(t) + 2\frac{\omega^2}{A^2} x(t)^3
-\frac{3}{4} \frac{\omega^2}{A^4} x(t)^5.
\label{L3}
\end{split}
\end{equation}

(ii) In the case where $(2m+1) \pi_2 \leqq \omega t \leqq (2m+2) \pi_2$, 
\begin{equation}
\begin{split}
\frac{\mathrm{d}^2 x(t) }{\mathrm{d} t^2}
&=   \frac{A \omega}{2} 
\frac{  \{ -1 +2 \mathrm{cleaf}_2(\omega t) 
-3 (\mathrm{cleaf}_2(\omega t))^2 \} 
( \omega \sqrt{1-(\mathrm{cleaf}_2(\omega t))^4 } )  }
{ 2 \sqrt{1- \mathrm{cleaf}_2(\omega t)
+(\mathrm{cleaf}_2(\omega t))^2
-(\mathrm{cleaf}_2(\omega t))^3} } \\
&=   -\frac{3}{2} \omega^2 x(t) + 2\frac{\omega^2}{A^2} x(t)^3
-\frac{3}{4} \frac{\omega^2}{A^4} x(t)^5.
\label{L4}
\end{split}
\end{equation}

Eq. (\ref{5.2}) of the Duffing equation is obtained. The initial condition of the differential equation, Eq. (\ref{5.2}), is as follows:
\begin{equation}
 x(0)=A \sqrt{1+ \mathrm{cleaf}_2(0) }=A \sqrt{1+ 1 }= \sqrt{2} A. \label{L5}
\end{equation}

Exact solution 11 is a continuous function. However,  it cannot be partially differentiated. 
Exact solution 11 is not differentiable at $t=(2m-1) \frac{\pi_2}{\omega}$:
\begin{equation}
\begin{split}
&\lim_{t \uparrow (2m-1) \frac{\pi_2}{\omega} }  \frac{\mathrm{d} x(t) }{\mathrm{d} t} 
=  - A \omega,
\label{L6}
\end{split}
\end{equation}

\begin{equation}
\begin{split}
&\lim_{t \downarrow (2m-1) \frac{\pi_2}{\omega} }  \frac{\mathrm{d} x(t) }{\mathrm{d} t} 
=   A \omega.
\label{L7}
\end{split}
\end{equation}

\section{Ordinary differential equation of exact solution 12  }
\renewcommand{\theequation}{M.\arabic{equation}}
\setcounter{equation}{0}

Here we prove that exact solution 12 satisfies the Duffing equation. The first derivatives of exact solution 12 are obtained as follows: 

(i) In the case where $2m \pi_2 \leqq \omega t \leqq (2m +1) \pi_2$,
\begin{equation}
\begin{split}
\frac{\mathrm{d} x(t) }{\mathrm{d} t} 
&= A \omega
\frac{  -(-\sqrt{1-(\mathrm{cleaf}_2(\omega t))^4 } ) }
{ 2 \sqrt{1- \mathrm{cleaf}_2(\omega t)} } \\
&= \frac{A \omega}{2} \sqrt{1+ \mathrm{cleaf}_2(\omega t)
+(\mathrm{cleaf}_2(\omega t))^2
+(\mathrm{cleaf}_2(\omega t))^3  }.
\label{M1}
\end{split}
\end{equation}

(ii) In the case where $(2m +1) \pi_2 \leqq \omega t \leqq (2m +2) \pi_2$, 
\begin{equation}
\begin{split}
\frac{\mathrm{d} x(t) }{\mathrm{d} t} 
&= A \omega
\frac{  -\sqrt{1-(\mathrm{cleaf}_2(\omega t))^4 }  }
{ 2 \sqrt{1- \mathrm{cleaf}_2(\omega t)} } \\
&= - \frac{A \omega}{2} \sqrt{1+ \mathrm{cleaf}_2(\omega t)
+(\mathrm{cleaf}_2(\omega t))^2
+(\mathrm{cleaf}_2(\omega t))^3  }.
\label{M2}
\end{split}
\end{equation}

The second derivatives of exact solution 12 are obtained as follows:

(i) In the case where $2m \pi_2 \leqq \omega t \leqq (2m +1) \pi_2$,
\begin{equation}
\begin{split}
\frac{\mathrm{d}^2 x(t) }{\mathrm{d} t^2}
&=    \frac{A \omega}{2} 
\frac{  \{ 1 +2 \mathrm{cleaf}_2(\omega t) 
+3 (\mathrm{cleaf}_2(\omega t))^2 \} 
(-\omega \sqrt{1-(\mathrm{cleaf}_2(\omega t))^4 } )  }
{ 2 \sqrt{1+ \mathrm{cleaf}_2(\omega t)
+(\mathrm{cleaf}_2(\omega t))^2
+(\mathrm{cleaf}_2(\omega t))^3} } \\
&=    - \frac{A \omega^2}{4}  
 \{ 1 +2 \mathrm{cleaf}_2(\omega t) 
+3 (\mathrm{cleaf}_2(\omega t))^2 \} 
\sqrt{1-\mathrm{cleaf}_2(\omega t) } 
\\
&=    - \frac{A \omega^2}{4}  
 \left\{ 1 +2 \left( 1 - \frac{x(t)^2}{A^2} \right)  +3 \left( 1-\frac{x(t)^2}{A^2} \right)^2 \right\} \frac{x(t)}{A} \\
&=   -\frac{3}{2} \omega^2 x(t) + 2\frac{\omega^2}{A^2} x(t)^3
-\frac{3}{4} \frac{\omega^2}{A^4} x(t)^5.
\label{M3}
\end{split}
\end{equation}

(ii) In the case where $(2m +1) \pi_2 \leqq \omega t \leqq (2m +2) \pi_2$, 
\begin{equation}
\begin{split}
\frac{\mathrm{d}^2 x(t) }{\mathrm{d} t^2}
&=   - \frac{A \omega}{2} 
\frac{  \{ 1 +2 \mathrm{cleaf}_2(\omega t) 
+3 (\mathrm{cleaf}_2(\omega t))^2 \} 
(\omega \sqrt{1-(\mathrm{cleaf}_2(\omega t))^4 } )  }
{ 2 \sqrt{1+ \mathrm{cleaf}_2(\omega t)
+(\mathrm{cleaf}_2(\omega t))^2
+(\mathrm{cleaf}_2(\omega t))^3} } \\
&=   -\frac{3}{2} \omega^2 x(t) + 2\frac{\omega^2}{A^2} x(t)^3
-\frac{3}{4} \frac{\omega^2}{A^4} x(t)^5.
\label{M4}
\end{split}
\end{equation}

Eq. (\ref{5.2}) of the Duffing equation is obtained. The initial condition of the differential equation, Eq. (\ref{5.2}), is as follows:
\begin{equation}
 x(0)=A \sqrt{1 - \mathrm{cleaf}_2(0) }=A \sqrt{1 - 1 }= 0. \label{M5}
\end{equation}

Exact solution 12 is a continuous function. However,  it cannot be partially differentiated. Exact solution 12 is not differentiable at $t=2m \frac{\pi_2}{\omega}$:
\begin{equation}
\begin{split}
&\lim_{t \uparrow 2m \frac{\pi_2}{\omega} }  \frac{\mathrm{d} x(t) }{\mathrm{d} t} 
=  - A \omega,
\label{M6}
\end{split}
\end{equation}

\begin{equation}
\begin{split}
&\lim_{t \downarrow 2m \frac{\pi_2}{\omega} }  \frac{\mathrm{d} x(t) }{\mathrm{d} t} 
=   A \omega.
\label{M7}
\end{split}
\end{equation}

\section{ Ordinary differential equation of exact solution 13 }
\renewcommand{\theequation}{N.\arabic{equation}}
\setcounter{equation}{0}

Here we prove that exact solution 13 satisfies the Duffing equation. The first derivatives of exact solution 13 are obtained as follows: 

(i) In the case where $4m \eta_2 \leqq \omega t < (4m+1) \eta_2$ and $(4m+1) \eta_2 < \omega t \leqq (4m+2) \eta_2$,
\begin{equation}
\begin{split}
\frac{\mathrm{d} x(t) }{\mathrm{d} t} 
&= A \omega
\frac{  \sqrt{(\mathrm{cleafh}_2(\omega t))^4-1 }  }
{ 2 \sqrt{\mathrm{cleafh}_2(\omega t)+1}   } \\
&= \frac{A \omega}{2} \sqrt{ (\mathrm{cleafh}_2(\omega t))^3-(\mathrm{cleafh}_2(\omega t))^2+ \mathrm{cleafh}_2(\omega t)-1 }.
\label{N1}
\end{split}
\end{equation}

(ii) In the case where $(4m-2) \eta_2 \leqq \omega t < (4m-1) \eta_2$ and $(4m-1) \eta_2 < \omega t \leqq 4m \eta_2$, 
\begin{equation}
\begin{split}
\frac{\mathrm{d} x(t) }{\mathrm{d} t} 
&= A 
\frac{ - \omega \sqrt{(\mathrm{cleafh}_2(\omega t))^4-1 }  }
{ 2 \sqrt{\mathrm{cleafh}_2(\omega t)+1}   } \\
&= -\frac{A \omega}{2} \sqrt{ (\mathrm{cleafh}_2(\omega t))^3-(\mathrm{cleafh}_2(\omega t))^2+ \mathrm{cleafh}_2(\omega t)-1 }.
\label{N2}
\end{split}
\end{equation}

With respect to all domains in (i) and (ii),
the second derivative of exact solution 13 is obtained as follows:
\begin{equation}
\begin{split}
\frac{\mathrm{d}^2 x(t) }{\mathrm{d} t^2}
=    \frac{3}{2} \omega^2 x(t) - 2\frac{\omega^2}{A^2} x(t)^3
+\frac{3}{4} \frac{\omega^2}{A^4} x(t)^5.
\label{N3}
\end{split}
\end{equation}

Eq. (\ref{5.7}) of the Duffing equation is obtained.  The initial conditions of the differential equation  are as follows:
\begin{equation}
 x(0)=A \sqrt{ \mathrm{cleafh}_2(0)+1 }=A \sqrt{1+ 1 }= \sqrt{2} A, \label{N4}
\end{equation}

\begin{equation}
\begin{split}
\frac{\mathrm{d} x(0) }{\mathrm{d} t}=0
\label{N5}
\end{split}
\end{equation}

\section{ Ordinary differential equation of exact solution 14 }
\renewcommand{\theequation}{O.\arabic{equation}}
\setcounter{equation}{0}

Here we prove that exact solution 14 satisfies the Duffing equation. The first derivatives of exact solution 14 are obtained as follows: 

(i) In the case where $4m \eta_2 \leqq \omega t < (4m+1) \eta_2$ and $(4m+1) \eta_2 < \omega t \leqq (4m+2) \eta_2$,
\begin{equation}
\begin{split}
\frac{\mathrm{d} x(t) }{\mathrm{d} t} 
&= A \omega
\frac{  \sqrt{(\mathrm{cleafh}_2(\omega t))^4-1 }  }
{ 2 \sqrt{\mathrm{cleafh}_2(\omega t)-1}   } \\
&= \frac{A \omega}{2} \sqrt{ (\mathrm{cleafh}_2(\omega t))^3 + (\mathrm{cleafh}_2(\omega t))^2+ \mathrm{cleafh}_2(\omega t)+1 }. 
\label{O1}
\end{split}
\end{equation}

(ii) In the case where $(4m-2) \eta_2 \leqq \omega t < (4m-1) \eta_2$ and $(4m-1) \eta_2 < \omega t \leqq 4m \eta_2$, 
\begin{equation}
\begin{split}
\frac{\mathrm{d} x(t) }{\mathrm{d} t} 
&= A \omega
\frac{  -\sqrt{(\mathrm{cleafh}_2(\omega t))^4-1 }  }
{ 2 \sqrt{\mathrm{cleafh}_2(\omega t)-1}   } \\
&= -\frac{A \omega}{2} \sqrt{ (\mathrm{cleafh}_2(\omega t))^3 + (\mathrm{cleafh}_2(\omega t))^2+ \mathrm{cleafh}_2(\omega t)+1 }. 
\label{O2}
\end{split}
\end{equation}

With respect to all domains in (i) and (ii),
the second derivative of exact solution 14 is obtained as follows:

\begin{equation}
\begin{split}
\frac{\mathrm{d}^2 x(t) }{\mathrm{d} t^2}
=    \frac{3}{2} \omega^2 x(t) + 2\frac{\omega^2}{A^2} x(t)^3
+\frac{3}{4} \frac{\omega^2}{A^4} x(t)^5.
\label{O3}
\end{split}
\end{equation}

Eq. (\ref{5.9}) of the Duffing equation is obtained.
The initial condition of the differential equation is as follows:
\begin{equation}
 x(0)=A \sqrt{ \mathrm{cleafh}_2(0)-1 }=A \sqrt{1 - 1 }= 0. \label{O4}
\end{equation}

Exact solution 14 is a continuous function. However,  it cannot be partially differentiated. Exact solution 14 is not differentiable at $t=4m \frac{\eta_2}{\omega}$:
\begin{equation}
\begin{split}
&\lim_{t \uparrow 4m \frac{\eta_2}{\omega} }  \frac{\mathrm{d} x(t) }{\mathrm{d} t} 
=  - A \omega,
\label{O5}
\end{split}
\end{equation}

\begin{equation}
\begin{split}
&\lim_{t \downarrow 4m \frac{\eta_2}{\omega} }  \frac{\mathrm{d} x(t) }{\mathrm{d} t} 
=  A \omega.
\label{O6}
\end{split}
\end{equation}

\section{ Limit of hyperbolic leaf function: $\mathrm{cleafh}_n(t)$  }
\renewcommand{\theequation}{P.\arabic{equation}}
\setcounter{equation}{0}

The limit of the hyperbolic leaf function $\mathrm{cleafh}_n(t)$ is calculated using the following equation \cite{Kaz_clh}:
\begin{equation}
\begin{split}
\eta_n=\int_1^\infty \frac{1}{\sqrt{u^{2n}-1}} du.
\label{P1}
\end{split}
\end{equation}

The numerical results of $\eta_n$ (for $n$=2,3 $\cdots$)  are summarized in Table \ref{tabP1}. 

\begin{table}
\begin{center}
\caption{ Limits $\eta_n$ of hyperbolic leaf functions; cleafh$_n(t)$}
\label{tabP1}
\begin{tabular}{cc}
\hline\noalign{\smallskip}
$n$ &  $ \eta_n$    \\
\noalign{\smallskip}\hline\noalign{\smallskip}
$1$ &  $ N/A $    \\
$2$ &  $\eta_2=1.311 \cdots $    \\
$3$ &  $\eta_3=0.701 \cdots $    \\
$\cdots$ &  $\cdots$    \\
\noalign{\smallskip}\hline
\end{tabular}
\end{center}
\end{table}

\section{ Period of leaf function }
\renewcommand{\theequation}{Q.\arabic{equation}}
\setcounter{equation}{0}

The constant $\pi_n$ is calculated using the following equation \cite{Kaz_sl} \cite{Kaz_cl}:
\begin{equation}
\begin{split}
\pi_n=2 \int_0^1 \frac{1}{\sqrt{1-u^{2n}}} du.
\label{Q1}
\end{split}
\end{equation}

The constant values 2$\pi_n$ represent one periodicity with respect to the arbitrary parameter $n$. The numerical results of $\pi_n$ (for $n$=1,2,3 $\cdots$) are summarized in Table \ref{tabQ1}.

\begin{table}
\begin{center}
\caption{ Constants $\pi_n$}
\label{tabQ1}
\begin{tabular}{cc}
\hline\noalign{\smallskip}
$n$ &  $ \pi_n$    \\
\noalign{\smallskip}\hline\noalign{\smallskip}
$1$ &  $\pi_1=3.141 \cdots $    \\
$2$ &  $\pi_2=2.622 \cdots $    \\
$3$ &  $\pi_3=2.429 \cdots $    \\
$\cdots$ &  $\cdots$    \\
\noalign{\smallskip}\hline
\end{tabular}
\end{center}
\end{table}

\section{ Limit of hyperbolic leaf function: $\mathrm{sleafh}_n(t)$  }
\renewcommand{\theequation}{R.\arabic{equation}}
\setcounter{equation}{0}

The constant $\zeta_n$ is calculated using the following equation \cite{Kaz_slh}:
\begin{equation}
\begin{split}
\zeta_n=\int_0^\infty \frac{1}{\sqrt{1+u^{2n}}} du.
\label{R1}
\end{split}
\end{equation}

The numerical results of $\zeta_n$ (for $n$=2,3 $\cdots$) are summarized in Table \ref{tabR1}.

\begin{table}
\begin{center}
\caption{ Limits $\zeta_n$ of hyperbolic leaf functions; sleafh$_n(t)$}
\label{tabR1}
\begin{tabular}{cc}
\hline\noalign{\smallskip}
$n$ &  $ \zeta_n$    \\
\noalign{\smallskip}\hline\noalign{\smallskip}
$1$ &  $ N/A $    \\
$2$ &  $\zeta_2=1.854 \cdots $    \\
$3$ &  $\zeta_3=1.402 \cdots $    \\
$\cdots$ &  $\cdots$    \\
\noalign{\smallskip}\hline
\end{tabular}
\end{center}
\end{table}

\section{ Relationship between $\mathrm{sleaf}_2(t)$ and  $\mathrm{cleaf}_2(t)$  }
\renewcommand{\theequation}{S.\arabic{equation}}
\setcounter{equation}{0}
The relation equations with the basis $n=2$ are described. The relation equation between the leaf function $\mathrm{sleaf}_2(t)$ and the leaf function $\mathrm{cleaf}_2(t)$ is as follows \cite{Kaz_cl}:
\begin{equation}
 (\mathrm{sleaf}_2(t))^2+ (\mathrm{cleaf}_2(t))^2+ (\mathrm{sleaf}_2(t))^2 \cdot (\mathrm{cleaf}_2(t))^2=1 .\label{S1}
\end{equation}


\end{document}